\documentclass[12pt]{article}
\usepackage{amssymb, amsmath}
\usepackage[dvips]{graphics}

\usepackage{bbm}
\usepackage{enumitem}
\usepackage{setspace}
\usepackage[colorlinks=true]{hyperref}
\hypersetup{colorlinks   = true}
\hypersetup{linkcolor=blue}
\begin{document}

\newtheorem{The}{Theorem}[section]
\newtheorem{Lem}[The]{Lemma}
\newtheorem{Prop}[The]{Proposition}
\newtheorem{Cor}[The]{Corollary}
\newtheorem{Rem}[The]{Remark}
\newtheorem{Obs}[The]{Observation}
\newtheorem{SConj}[The]{Standard Conjecture}
\newtheorem{Titre}[The]{\!\!\!\! }
\newtheorem{Conj}[The]{Conjecture}
\newtheorem{Question}[The]{Question}
\newtheorem{Prob}[The]{Problem}
\newtheorem{Def}[The]{Definition}
\newtheorem{Not}[The]{Notation}
\newtheorem{Claim}[The]{Claim}
\newtheorem{Conc}[The]{Conclusion}
\newtheorem{Ex}[The]{Example}
\newtheorem{Fact}[The]{Fact}
\newcommand{\C}{\mathbb{C}}
\newcommand{\R}{\mathbb{R}}
\newcommand{\N}{\mathbb{N}}
\newcommand{\Z}{\mathbb{Z}}
\newcommand{\Q}{\mathbb{Q}}
\newcommand{\Proj}{\mathbb{P}}

\begin{center}

{\Large\bf Volume and Self-Intersection of Differences of Two Nef Classes}

\end{center}

\begin{center}

{\large Dan Popovici}

\end{center}

\vspace{1ex}

\noindent {\small {\bf Abstract.} Let $\{\alpha\}$ and $\{\beta\}$ be nef cohomology classes of bidegree $(1,\,1)$ on a compact $n$-dimensional K\"ahler manifold $X$ such that the difference of intersection numbers $\{\alpha\}^n - n\,\{\alpha\}^{n-1}.\,\{\beta\}$ is positive. We solve in a number of special but rather inclusive cases the quantitative part of Demailly's Transcendental Morse Inequalities Conjecture for this context predicting the lower bound $\{\alpha\}^n-n\,\{\alpha\}^{n-1}.\,\{\beta\}$ for the volume of the difference class $\{\alpha-\beta\}$. We completely solved the qualitative part in an earlier work. We also give general lower bounds for the volume of $\{\alpha-\beta\}$ and show that the self-intersection number $\{\alpha-\beta\}^n$ is always bounded below by $\{\alpha\}^n-n\,\{\alpha\}^{n-1}.\,\{\beta\}$. We also describe and estimate the relative psef and nef thresholds of $\{\alpha\}$ with respect to $\{\beta\}$ and relate them to the volume of $\{\alpha-\beta\}$. Finally, broadening the scope beyond the K\"ahler realm, we propose a conjecture relating the balanced and the Gauduchon cones of $\partial\bar\partial$-manifolds which, if proved to hold, would imply the existence of a balanced metric on any $\partial\bar\partial$-manifold.}

\vspace{1ex}


\vspace{1ex}

\section {Introduction}

Let $X$ be a compact K\"ahler manifold with $\mbox{dim}_{\C}X=n$ and let $\{\alpha\}, \{\beta\}\in H^{1,\,1}_{BC}(X,\,\R)$ be {\bf nef} Bott-Chern cohomology classes such that

\begin{equation}\label{eqn:initial-pos-cond}\{\alpha\}^n - n\,\{\alpha\}^{n-1}.\,\{\beta\} >0.\end{equation}

\noindent A (possibly transcendental) class $\{\alpha\}\in H^{1,\,1}_{BC}(X,\,\R)$ being {\bf nef} means (cf.\!\! Definition $1.3$ in [Dem92]) that for some (hence all) fixed Hermitian metric $\omega$ on $X$ and for every $\varepsilon>0$, there exists a $C^{\infty}$ form $\alpha_{\varepsilon}\in\{\alpha\}$ such that $\alpha_{\varepsilon}\geq -\varepsilon\,\omega$.

We have proved in [Pop14, Theorem $1.1$] that the class $\{\alpha - \beta\}$ is {\bf big} (i.e. contains a K\"ahler current $T$). This solved the {\it qualitative part} of Demailly's {\it Transcendental Morse Inequalities Conjecture} for differences of two nef classes (cf. [BDPP13, Conjecture $10.1,\,(ii)$]) on compact K\"ahler (and even more general) manifolds. This special form of the conjecture was originally motivated by attempts at extending to transcendental classes and to compact K\"ahler (not necessarily projective) manifolds the cone duality theorem of Boucksom, Demailly, Paun and Peternell [BDPP13, Theorem $2.2.$] that plays a major role in the theory of classification of projective manifolds. Recall that $T$ being a {\it K\"ahler current} means that $T$ is a $d$-closed positive $(1,\,1)$-current with the property that for some (hence all) fixed Hermitian metric $\omega$ on $X$, there exists $\varepsilon>0$ such that $T\geq\varepsilon\,\omega$ on $X$. Nefness and bigness are quite different positivity properties for real (possibly transcendental) $(1,\,1)$-classes and the by-now standard definitions just recalled extend classical algebraic definitions for integral classes on projective manifolds.

 In this paper we give a partial answer to the {\it quantitative part} of Demailly's {\it Transcendental Morse Inequalities Conjecture} for differences of two nef classes:

\begin{Conj}([BDPP13, Conjecture $10.1,\,(ii)$])\label{Conj:volume_lbound} Let $\{\alpha\}, \{\beta\}\in H^{1,\,1}_{BC}(X,\,\R)$ be {\bf nef} classes satisfying condition (\ref{eqn:initial-pos-cond}) on a compact K\"ahler manifold $X$ with $\mbox{dim}_{\C}X=n$. Then

\begin{equation}\label{eqn:volume_lbound}\mbox{Vol}\,\,(\{\alpha - \beta\}) \geq \{\alpha\}^n - n\,\{\alpha\}^{n-1}.\,\{\beta\}.\end{equation}

\end{Conj}

This is stated for arbitrary (i.e. possibly non-K\"ahler) compact complex manifolds in [BDPP13], but the volume is currently only known to be meaningful when $X$ is of {\it class} ${\cal C}$, a case reducible to the K\"ahler case by modifications. Thus, we may assume without loss of generality that $X$ is K\"ahler.

 Recall that the volume is a way of gauging the ``amount'' of positivity of a class $\{\gamma\}\in H^{1,\,1}_{BC}(X,\,\R)$ when $X$ is K\"ahler (or merely of {\it class} ${\cal C}$) and was introduced in [Bou02, Definition $1.3$] as

\begin{equation}\label{eqn:volume_def}\mbox{Vol}\,\,(\{\gamma\}):=\sup\limits_{T\in\{\gamma\},\,T\geq 0}\int\limits_XT_{ac}^n\end{equation}

\noindent if $\{\gamma\}$ is {\it pseudo-effective (psef)}, i.e. if $\{\gamma\}$ contains a {\it positive} $(1,\,1)$-current $T\geq 0$, where $T_{ac}$ denotes the absolutely continuous part of $T$ in the Lebesgue decomposition of its coefficients (which are complex measures when $T\geq 0$). If the class $\{\gamma\}$ is not psef, then its volume is set to be zero. It was proved in [Bou02, Theorem $1.2$] that this volume (which is always a {\it finite} non-negative quantity thanks to the K\"ahler, or more generally {\it class} ${\cal C}$, assumption on $X$) coincides with the standard volume of a holomorphic line bundle $L$ if the class $\{\gamma\}$ is integral (i.e. the first Chern class of some $L$). Moreover, the class $\{\gamma\}$ is big (i.e. contains a K\"ahler current) if and only if its volume is {\it positive}, by [Bou02, Theorem $4.7$]. 

 Thus, under the assumptions of Conjecture \ref{Conj:volume_lbound}, the main result in [Pop14] ensures that $\mbox{Vol}\,\,(\{\alpha - \beta\})>0$. In other words, $\{\alpha - \beta\}$ is positive in the {\it big} sense. The special case when $\{\beta\}=0$ had been proved in [DP04, Theorem $2.12$] and had served there as the main ingredient in the proof of the numerical characterisation of the K\"ahler cone. (In particular, the proof of the more general statement in [Pop14] reproves in a much simpler way the main technical result in [DP04].) The thrust of Conjecture \ref{Conj:volume_lbound} is to estimate from below the ``amount'' of positivity of the class $\{\alpha - \beta\}$.

 A first group of results that we obtain in the present paper can be summed up in the following positive answer to Conjecture \ref{Conj:volume_lbound} under an extra assumption. Recall that for {\it nef} classes $\{\gamma\}$, the volume equals the top self-intersection $\{\gamma\}^n$ (cf. [Bou02, Theorem $4.1$]), but for arbitrary classes, any order may occur between these two quantities.

\begin{The}\label{The:volume_lbound_all} Let $X$ be compact K\"ahler manifold with $\mbox{dim}_{\C}X=n$ and let $\{\alpha\}, \{\beta\}\in H^{1,\,1}_{BC}(X,\,\R)$ be {\bf nef} classes such that $\{\alpha\}^n - n\,\{\alpha\}^{n-1}.\,\{\beta\} >0$. Suppose, moreover, that 

\begin{equation}\label{eqn:extra-assumption}\mbox{Vol}\,\,(\{\alpha - \beta\}) \geq \{\alpha-\beta\}^n.\end{equation} 

\noindent Then $\mbox{Vol}\,\,(\{\alpha - \beta\}) \geq \{\alpha\}^n - n\,\{\alpha\}^{n-1}.\,\{\beta\}.$

\end{The}

 Although there are examples when the volume of $\{\alpha - \beta\}$ is strictly less than the top self-intersection, the assumption (\ref{eqn:extra-assumption}), that we hope to be able to remove in future work, is satisfied in quite a number of cases, e.g. when the class $\{\alpha - \beta\}$ is nef (treated in section \ref{section:nef_difference}). 

 Actually, we prove in full generality in section \ref{section:intersection-numbers} the analogue of Conjecture \ref{Conj:volume_lbound} for the top self-intersection number $\{\alpha - \beta\}^n$ in place of the volume of $\{\alpha - \beta\}$.

\begin{The}\label{The:s-intersection_lbound} Let $X$ be a compact K\"ahler manifold with $\mbox{dim}_{\C}X=n$ and let $\{\alpha\}, \{\beta\}\in H^{1,\,1}_{BC}(X,\,\R)$ be {\bf nef} classes such that $\{\alpha\}^n - n\,\{\alpha\}^{n-1}.\,\{\beta\} >0$. Then 

$$\{\alpha - \beta\}^n \geq \{\alpha\}^n - n\,\{\alpha\}^{n-1}.\,\{\beta\}.$$

\end{The}

 It is clear that Theorem \ref{The:volume_lbound_all} follows immediately from Theorem \ref{The:s-intersection_lbound}. Since the nef cone is the closure of the K\"ahler cone, we may assume without loss of generality that the classes $\{\alpha\}$ and $\{\beta\}$ are actually K\"ahler. Corcerning the volume of $\{\alpha - \beta\}$ in the general case (i.e. without assumption (\ref{eqn:extra-assumption})), we prove a lower bound that is weaker than the expected lower bound (\ref{eqn:volume_lbound}) in a way that depends explicitly on how far the class $\{\alpha - \beta\}$ is from being nef. The {\it nefness defect} of $\{\alpha - \beta\}$ is defined explicitly and investigated in relation to the volume in subsections \ref{subsection:nef}, \ref{subsection:trace-volume} and \ref{subsection:nef/psef/vol_again}. We call it the {\it nef threshold} (a term that is already present in the literature) of $\{\alpha\}$ w.r.t. $\{\beta\}$ and discuss it together with the analogous {\it psef threshold} of $\{\alpha\}$ w.r.t. $\{\beta\}$ in section \ref{section:dir-trace}. In $\S.$\ref{subsection:nef/psef/vol_again}, we prove the following general lower bound for the volume of $\{\alpha - \beta\}$.

\begin{The}\label{The:volume_lbound_nefT} Let $X$ be compact K\"ahler manifold with $\mbox{dim}_{\C}X=n$ and let $\{\alpha\}, \{\beta\}\in H^{1,\,1}_{BC}(X,\,\R)$ be {\bf K\"ahler} classes such that $\{\alpha\}^n - n\,\{\alpha\}^{n-1}.\,\{\beta\} >0$. Let $s_0:=N^{(\beta)}(\alpha)>0$ be the nef threshold of $\{\alpha\}$ w.r.t. $\{\beta\}$. Then:

\vspace{1ex}

$(i)$\, if $s_0\geq 1$, the class $\{\alpha - \beta\}$ is nef and the optimal volume estimate (\ref{eqn:volume_lbound}) holds;

\vspace{1ex}

$(ii)$\, if $s_0< 1$, the class $\{\alpha - \beta\}$ is not nef and the next volume estimate holds:

\begin{equation}\label{eqn:volume_lbound_nefT}\mbox{Vol}\,\,(\{\alpha - \beta\}) \geq (\{\alpha\}^n - n\,\{\alpha\}^{n-1}.\,\{\beta\})\,\bigg(\frac{\{\alpha\}^n - n\,\{\alpha\}^{n-1}.\,\{\beta\}}{\{\alpha\}^n - ns_0\,\{\alpha\}^{n-1}.\,\{\beta\}}\bigg)^{n-1}.\end{equation}

\end{The}

 A very special case of this result was also observed independently in [Tos15] using the technique introduced in [Pop14].

Taking our cue from the estimates we obtain in section \ref{section:M-A} for the supremum of $t\geq 0$ such that the class $\{\alpha\} -t\,\{\beta\}$ is psef in the setting of Conjecture \ref{Conj:volume_lbound}, we define the {\bf psef} and the {\bf nef thresholds} of $\{\alpha\}$ w.r.t. $\{\beta\}$ as functions $P^{(\beta)},\,N^{(\beta)}:H^{1,\,1}_{BC}(X,\,\R)\rightarrow\R$ by

\vspace{1ex}

\hspace{3ex} $\displaystyle (i)\hspace{2ex} P^{(\beta)}(\alpha):=\inf\int\limits_X\alpha\wedge\gamma^{n-1} \hspace{2ex} \mbox{and} \hspace{2ex} (ii)\hspace{2ex} N^{(\beta)}(\alpha):=\inf\int\limits_Y\alpha\wedge\omega^{n-p-1},$

\vspace{1ex}

\noindent where in $(i)$ the infimum is taken over all the Gauduchon metrics $\gamma$ on $X$ normalised by 

\vspace{1ex}

\hspace{20ex} $\displaystyle [\beta]_{BC}.[\gamma^{n-1}]_A=\int\limits_X\beta\wedge\gamma^{n-1}=1,$

\vspace{1ex}

\noindent while in $(ii)$ the infimum is taken over all $p=0,1, \dots , n-1$, over all the irreducible analytic subsets $Y\subset X$ such that $\mbox{codim}\,Y=p$ and over all K\"ahler classes $\{\omega\}$ normalised by $\int_Y\beta\wedge\omega^{n-p-1}=1$. The class $\{\beta\}$ is supposed to be {\bf big} in the case of $P^{(\beta)}$ and {\bf K\"ahler} in the case of $N^{(\beta)}$. (The subscripts BC and A will stand throughout for the Bott-Chern, resp. Aeppli cohomologies.) In section \ref{section:dir-trace}, we prove the following formulae that justify the terminology and make it match existing notions in the literature:

\vspace{1ex}

\hspace{6ex}  $P^{(\beta)}(\alpha) = \sup\,\{t\in\R\,\,\slash\,\,\mbox{the class}\,\,\{\alpha\} - t\,\{\beta\}\,\,\mbox{is psef}\},$

\vspace{1ex}

\hspace{6ex}  $N^{(\beta)}(\alpha) = \sup\{s\in\R\,\,/\,\,\mbox{the class}\,\,\{\alpha\}-s\{\beta\}\,\,\mbox{is nef}\,\}.$

\vspace{1ex}

The psef/nef thresholds of $\{\alpha\}$ w.r.t. $\{\beta\}$ turn out to gauge quite effectively the amount of positivity that the class $\{\alpha\}$ has in the ``direction'' of the class $\{\beta\}$. We study their various properties in section \ref{section:dir-trace}, estimate them in terms of intersection numbers as

\vspace{1ex}

\hspace{20ex} $\displaystyle\frac{\{\alpha\}^n}{n\,\{\alpha\}^{n-1}.\{\beta\}} \leq P^{(\beta)}(\alpha) \leq \frac{\{\alpha\}^n}{\{\alpha\}^{n-1}.\{\beta\}}$

\vspace{1ex}

\noindent and by similar, more involved inequalities for $N^{(\beta)}(\alpha)$, and relate them to the volume of $\{\alpha - \beta\}$ as

\vspace{1ex}

\hspace{20ex} $\displaystyle\mbox{Vol}\,(\{\alpha-\beta\})\geq\bigg(1 - \frac{1}{P^{(\beta)}(\alpha)}\bigg)^n\,\{\alpha\}^n$

\vspace{1ex}

\noindent whenever the classes $\{\alpha\}$ and $\{\beta\}$ are K\"ahler.

 Using these thresholds, we prove Conjecture \ref{Conj:volume_lbound} in yet another special case: when the psef and the nef thresholds of $\{\alpha\}$ w.r.t. $\{\beta\}$ are sufficiently {\it close} to each other (cf. Proposition \ref{Prop:psef-nef-vol}). Of course, we always have: $N^{(\beta)}(\alpha)\leq P^{(\beta)}(\alpha)$.

\vspace{2ex}

 As in our earlier work [Pop14] and as in [Xia13] that preceded it, we will repeatedly make use of two ingredients. The first one is {\it Lamari's positivity criterion}.

\begin{Lem}([Lam99, Lemme 3.3])\label{Lem:Lamari} Let $\{\alpha\}\in H^{1,\,1}_{BC}(X,\,\R)$ be any real Bott-Chern cohomology class on an $n$-dimensional compact complex manifold $X$. The following two statements are equivalent. 

\vspace{1ex}

 $(i)$\, There exists a $(1,\,1)$-current $T$ in $\{\alpha\}$ such that $T\geq 0$ on $X$ (i.e. $\{\alpha\}$ is psef). 

\vspace{1ex}

 $(ii)$\, $\displaystyle\int\limits_X\alpha\wedge\gamma^{n-1}\geq 0$ for all Gauduchon metrics $\gamma$ on $X.$

\end{Lem}

 In fact, Lamari's result holds more generally for any (i.e. not necessarily $d$-closed) $C^{\infty}$ real $(1,\,1)$-form $\alpha$ on $X$, but we will not use this here. The second ingredient that we will often use is Yau's solution of the Calabi Conjecture.

\begin{The}([Yau78])\label{The:Yau} Let $X$ be a compact complex $n$-dimensional manifold endowed with a K\"ahler metric $\omega$. Let $dV>0$ be any $C^{\infty}$ positive volume form on $X$ such that $\int_X\omega^n = \int_XdV$. Then, there exists a unique K\"ahler metric $\widetilde{\omega}$ in the K\"ahler class $\{\omega\}$ such that $\widetilde{\omega}^n=dV$.

\end{The}

 There is a non-K\"ahler analogue of Yau's theorem by Tosatti and Weinkove [TW10] that will not be used in this work. Moreover, most of the techniques that follow are still meaningful or can be extended to the non-K\"ahler context. This is part of the reason why we believe that a future development of the matters dealt with in this paper may be possible in the more general setting of $\partial\bar\partial$-manifolds. The conjecture we propose in section \ref{section:nKapplication} is an apt illustration of this idea.

 We will make repeated use of the technique based on the Cauchy-Schwarz inequality for estimating from below certain integrals of traces of K\"ahler metrics introduced in [Pop14]. Moreover, there are mainly two new techniques that we introduce in the current paper: $(i)$\, the observation, proof and use of certain pointwise inequalities involving products of positive smooth forms (cf. Appendix) reminiscent of the Hovanskii-Teissier inequalities and generalising Lemma $3.1$ in [Pop14]; $(ii)$\, a technique for constructing what we call {\it approximate fixed points} for Monge-Amp\`ere equations when we allow the r.h.s. to vary (cf. proof of Proposition \ref{Prop:power-diff}) whose rough idea originates in and was suggested by discussions the author had several years ago in a completely different context with different equations and for very different purposes with J.-P. Demailly to whom we are very grateful.

\section{Special case of Conjecture \ref{Conj:volume_lbound} when $\{\alpha - \beta\}$ is nef}\label{section:nef_difference}

 We start by noticing the following elementary inequality.

\begin{Lem}\label{Lem:semi-pos} Let $\alpha>0$ and $\beta\geq 0$ be $C^{\infty}$ $(1,\,1)$-forms on a complex manifold $X$ with $\mbox{dim}_{\C}X=n$ such that $\alpha-\beta \geq 0$. Then : 

\begin{equation}\label{eqn:form-ineq}(\alpha-\beta)^n \geq \alpha^n - n\,\alpha^{n-1}\wedge\beta   \hspace{3ex} \mbox{at every point in}\hspace{1ex} X.\end{equation}

\noindent If $d\alpha = d\beta = 0$ and if $X$ is compact, then taking integrals we get:

\begin{equation}\label{eqn:class-ineq}\mbox{Vol}\,\,(\{\alpha - \beta\}) = \int_X(\alpha-\beta)^n \geq \int_X\alpha^n - n\,\int_X\alpha^{n-1}\wedge\beta = \{\alpha\}^n - n\,\{\alpha\}^{n-1}.\,\{\beta\}.\end{equation}

\end{Lem}

\noindent {\it Proof.} Let $x_0\in X$ be an arbitrary point and let $z_1, \dots , z_n$ be local holomorphic coordinates centred at $x_0$ such that at $x_0$ we have:

$$\alpha = \sum\limits_{j=1}^nidz_j\wedge d\bar{z}_j  \hspace{2ex} \mbox{and}  \hspace{2ex} \beta = \sum\limits_{j=1}^n\beta_j\,idz_j\wedge d\bar{z}_j.$$

\noindent Then $\alpha - \beta = \sum\limits_{j=1}^n(1 - \beta_j)\,idz_j\wedge d\bar{z}_j$ at $x_0$, while $\beta_j\geq 0$ and $1-\beta_j\geq 0$ at $x_0$ for all $j$. Thus inequality (\ref{eqn:form-ineq}) at $x_0$ translates to

\begin{equation}\label{eqn:elementary-ineq}(1 - \beta_1)\cdots(1 - \beta_n) \geq 1 - (\beta_1 + \cdots + \beta_n) \hspace{3ex} \mbox{for all}\hspace{1ex}\beta_1, \dots , \beta_n\in[0,\,1].\end{equation}

\noindent This elementary inequality is easily proved by induction on $n\geq 1$. Indeed, (\ref{eqn:elementary-ineq}) is an identity for $n=1$, while if (\ref{eqn:elementary-ineq}) has been proved for $n$, then we have:

\begin{eqnarray}\label{eqn:induction}\nonumber (1 - \beta_1)\cdots(1 - \beta_n)\,(1 - \beta_{n+1}) & \stackrel{(i)}{\geq} & \bigg(1 - (\beta_1 + \cdots + \beta_n)\bigg)\,(1 - \beta_{n+1})\\
\nonumber  & = & 1 - (\beta_1 + \cdots + \beta_n + \beta_{n+1}) + \beta_{n+1}\,(\beta_1 + \cdots + \beta_n)\\
\nonumber & \geq & 1 - (\beta_1 + \cdots + \beta_n + \beta_{n+1})\end{eqnarray}

\noindent since $\beta_j\geq 0$ for all $j$. (We used $1-\beta_{n+1}\geq 0$ to get $(i)$ from the induction hypothesis.) Thus (\ref{eqn:elementary-ineq}) is proved and (\ref{eqn:form-ineq}) follows from it.

 Now, if $\alpha$ and $\beta$ are $d$-closed, they define Bott-Chern cohomology classes. Since $\alpha-\beta$ is a semi-positive $C^{\infty}$ $(1,\,1)$-form, its Bott-Chern class is nef (and even a bit more), hence its volume equals $\int_X(\alpha-\beta)^n$ by [Bou02, Theorem 4.1] if $X$ is compact. (Note that $X$ is compact K\"ahler since $\alpha$ is a K\"ahler metric under the present assumptions.) The remaining part of (\ref{eqn:class-ineq}) follows at once from (\ref{eqn:form-ineq}) by integration.   \hfill $\Box$

\vspace{2ex}

 An immediate consequence of Lemma \ref{Lem:semi-pos} is the desired volume lower bound (\ref{eqn:volume_lbound}) in the special case when the class $\{\alpha-\beta\}$ is assumed to be nef. Note, however, that $\{\alpha-\beta\}$ need not be nef in general even a posteriori in the setting of Conjecture \ref{Conj:volume_lbound}.  

\begin{Prop}\label{Prop:nef-case} Let $X$ be a compact K\"ahler manifold with $\mbox{dim}_{\C}X=n$ and let $\{\alpha\}, \{\beta\}\in H^{1,\,1}_{BC}(X,\,\R)$ be {\bf nef} Bott-Chern cohomology classes such that the class $\{\alpha - \beta\}$ is {\bf nef}. Then 

\begin{equation}\label{eqn:volume_lbound_re}\mbox{Vol}\,\,(\{\alpha - \beta\}) \geq \{\alpha\}^n - n\,\{\alpha\}^{n-1}.\,\{\beta\}.\end{equation}

\end{Prop}

\noindent {\it Proof.} It suffices to prove inequality (\ref{eqn:volume_lbound_re}) in the case when the classes $\{\alpha\}, \{\beta\}$ and $\{\alpha-\beta\}$ are all K\"ahler. (Otherwise, we can add $2\varepsilon\{\omega\}$ to $\{\alpha\}$ and $\varepsilon\{\omega\}$ to $\{\beta\}$ for a fixed K\"ahler class $\{\omega\}$ and let $\varepsilon\downarrow 0$ in the end. The volume function is known to be continuous by [Bou02, Corollary $4.11$].) If we define the form $\alpha$ as the sum of any K\"ahler metric in the class $\{\alpha-\beta\}$ with any K\"ahler metric $\beta$ in the class $\{\beta\}$, the forms $\alpha$, $\beta$ and $\alpha-\beta$ obtained in this way satisfy the hypotheses of Lemma \ref{Lem:semi-pos}, hence also the elementary inequality (\ref{eqn:form-ineq}) and its consequence (\ref{eqn:class-ineq}).  \hfill $\Box$

\vspace{2ex}

 Recall that the class $\{\alpha - \beta\}$ is big under the assumptions of Conjecture \ref{Conj:volume_lbound} by the main result in [Pop14]. However, big positivity is quite different in nature to nef positivity. The general (i.e. possibly non-nef) case is discussed in the next sections.

\section{Applications of Monge-Amp\`ere equations}\label{section:M-A}

 In this section, we rewrite in a more effective way and observe certain consequences of the arguments in [Pop14, $\S.3$].

\begin{Lem}\label{Lem:prod-traces-int} Let $X$ be any compact complex manifold with $\mbox{dim}_{\C}X=n$. With any $C^{\infty}$ $(1,\,1)$-forms $\alpha, \beta>0$ and any Gauduchon metric $\gamma$, we associate the $C^{\infty}$ $(1,\,1)$-form $\widetilde{\alpha} = \alpha + i\partial\bar\partial u >0$ defined as the unique normalised solution (whose existence is guaranteed by the Tosatti-Weinkove theorem in [TW10]) of the Monge-Amp\`ere equation:

\begin{equation}\label{eqn:MA}(\alpha + i\partial\bar\partial u)^n = c\,\beta\wedge\gamma^{n-1} \hspace{2ex} \mbox{such that} \hspace{2ex} \sup\limits_X u=0,\end{equation}

\noindent where $c>0$ is the unique constant for which the above equation admits a solution $u\,:\,X\rightarrow\R$. (Of course, a posteriori, $c = (\int_X(\alpha + i\partial\bar\partial u)^n)/(\int_X\beta\wedge\gamma^{n-1})$, while if $d\alpha=0$ then $c = \int_X\alpha^n = \{\alpha\}^n>0$.) 

 Then the following inequality holds:

\begin{equation}\label{eqn:prod-traces-int}\bigg(\int\limits_X\widetilde{\alpha}\wedge\gamma^{n-1}\bigg)\cdot\bigg(\int\limits_X\widetilde{\alpha}^{n-1}\wedge\beta\bigg) \geq \frac{1}{n}\,\bigg(\int\limits_X\widetilde{\alpha}^n\bigg)\,\bigg(\int\limits_X\beta\wedge\gamma^{n-1}\bigg).\end{equation}

\end{Lem}

\noindent {\it Proof.} Let us define $\det_{\gamma}\widetilde{\alpha}$ by requiring $\widetilde{\alpha}^n = (\det_{\gamma}\widetilde{\alpha})\,\gamma^n$ on $X$. Since $\beta\wedge\gamma^{n-1} = (1/n)\,(\Lambda_{\gamma}\beta)\,\gamma^n$, the Monge-Amp\'ere equation (\ref{eqn:MA}) translates to

\begin{equation}\label{eqn:MA-translated}\mbox{det}_{\gamma}\widetilde{\alpha} = \frac{c}{n}\,\Lambda_{\gamma}\beta.\end{equation}

\noindent Hence, we get the following identities and inequalities:\\

\noindent $\displaystyle\bigg(\int_X\widetilde{\alpha}\wedge\gamma^{n-1}\bigg)\,\bigg(\int\limits_X\widetilde{\alpha}^{n-1}\wedge\beta\bigg) = \bigg(\int_X\frac{1}{n}\,(\Lambda_{\gamma}\widetilde{\alpha})\,\gamma^n\bigg)\,\bigg(\int_X\frac{1}{n}\,(\Lambda_{\widetilde{\alpha}}\beta)\,(\mbox{det}_{\gamma}\widetilde{\alpha})\,\gamma^n\bigg)$

\noindent \begin{eqnarray}\nonumber \stackrel{(a)}{\geq} \frac{1}{n^2}\,\bigg(\int\limits_X[(\Lambda_{\gamma}\widetilde{\alpha})\,(\Lambda_{\widetilde{\alpha}}\beta)]^{\frac{1}{2}}\,(\mbox{det}_{\gamma}\widetilde{\alpha})^{\frac{1}{2}}\,\gamma^n\bigg)^2 & \stackrel{(b)}{\geq} & \frac{1}{n^2}\,\bigg(\int\limits_X(\Lambda_{\gamma}\beta)^{\frac{1}{2}}\,(\mbox{det}_{\gamma}\widetilde{\alpha})^{\frac{1}{2}}\,\gamma^n\bigg)^2\\
\nonumber \stackrel{(c)}{=} \frac{1}{n^2}\,\bigg(\sqrt{\frac{c}{n}}\,\int\limits_X(\Lambda_{\gamma}\beta)\,\gamma^n\bigg)^2 & = & \frac{1}{n^2}\,\bigg(\sqrt{\frac{c}{n}}\,n\,\int\limits_X\beta\wedge\gamma^{n-1}\bigg)^2 = \frac{c}{n}\,\bigg(\int\limits_X\beta\wedge\gamma^{n-1}\bigg)^2,\end{eqnarray}

\noindent which prove (\ref{eqn:prod-traces-int}) since $c = (\int_X\widetilde{\alpha}^n)/(\int_X\beta\wedge\gamma^{n-1})>0$, where $(a)$ is the Cauchy-Schwarz inequality, $(b)$ follows from the inequality $(\Lambda_{\gamma}\widetilde{\alpha})\,(\Lambda_{\widetilde{\alpha}}\beta)\geq\Lambda_{\gamma}\beta$ proved in [Pop14, Lemma 3.1], while $(c)$ follows from (\ref{eqn:MA-translated}).  \hfill $\Box$

\begin{Cor}\label{Cor:prod-traces-direct}Let $X$ be a compact K\"ahler manifold with $\mbox{dim}_{\C}X=n$. Then, for every K\"ahler metrics $\alpha, \beta$ and every Gauduchon metric $\gamma$ on $X$, the following inequality holds:

\begin{equation}\label{eqn:prod-traces-direct}\bigg(\int\limits_X\alpha\wedge\gamma^{n-1}\bigg)\cdot\bigg(\int\limits_X\alpha^{n-1}\wedge\beta\bigg) \geq \frac{1}{n}\,\bigg(\int\limits_X\alpha^n\bigg)\,\bigg(\int\limits_X\beta\wedge\gamma^{n-1}\bigg).\end{equation}

\end{Cor}

\noindent {\it Proof.} It is clear that (\ref{eqn:prod-traces-direct}) follows immediately from (\ref{eqn:prod-traces-int}) since the assumption $d\alpha = d\beta = 0$ ensures that $\int_X\widetilde{\alpha}\wedge\gamma^{n-1}  = \int_X\alpha\wedge\gamma^{n-1}, \int_X\widetilde{\alpha}^{n-1}\wedge\beta  = \int_X\alpha^{n-1}\wedge\beta$ and $\int_X\widetilde{\alpha}^n = \int_X\alpha^n$.  \hfill $\Box$

\begin{Rem}\label{Rem:different-normalisation} Under the hypotheses of Corollary \ref{Cor:prod-traces-direct}, for any $\gamma$ satisfying the inequality $(\Lambda_{\gamma}\widetilde{\alpha})\,(\Lambda_{\widetilde{\alpha}}\beta)\geq n\Lambda_{\gamma}\beta$ (an improved version of [Pop14, Lemma 3.1] which need not hold in general, but holds for some special choices of $\gamma$ -- cf. proof of Lemma \ref{Lem:starstar-ineq}), the lower bound on the r.h.s. of (\ref{eqn:prod-traces-direct}) improves to $(\int_X\alpha^n)\,(\int_X\beta\wedge\gamma^{n-1}).$ If this improved lower bound held for all Gauduchon metrics $\gamma$, Conjecture \ref{Conj:volume_lbound} would follow immediately (see Theorem \ref{The:volume_n2_lbound} below).

\end{Rem}

 We first notice a consequence of Corollary \ref{Cor:prod-traces-direct} for nef classes.

\begin{Cor}\label{Cor:non-orthogonality} If $\{\alpha\}$ and $\{\beta\}\in H^{1,\,1}_{BC}(X,\,\R)$ are {\bf nef} classes on a compact K\"ahler manifold $X$ with $\mbox{dim}_{\C}X=n$ such that $\{\alpha\}^n - n\,\{\alpha\}^{n-1}.\,\{\beta\}>0$, then $\{\alpha\}^n>0$ and, unless $\{\beta\}=0$, the following non-orthogonality property holds: $\{\alpha\}^{n-1}.\,\{\beta\}>0$.

\end{Cor}

\noindent {\it Proof.} The nef hypothesis on $\{\alpha\}$ and $\{\beta\}$ ensures that $\{\alpha\}^{n-1}.\,\{\beta\}\geq 0$, hence $\{\alpha\}^n>0$ since $\{\alpha\}^n > n\,\{\alpha\}^{n-1}.\,\{\beta\}$ by assumption. For the rest of the proof, we reason by contradiction: suppose that $\{\alpha\}^{n-1}.\,\{\beta\}= 0$ and that $\{\beta\}\neq 0$. By the nef hypothesis on $\{\alpha\}$ and $\{\beta\}$, for every $\varepsilon>0$, there exist $C^{\infty}$ forms $\alpha_{\varepsilon}\in\{\alpha\}, \beta_{\varepsilon}\in\{\beta\}$ such that $\alpha_{\varepsilon} + \varepsilon\,\omega>0$ and $\beta_{\varepsilon} + \varepsilon\,\omega>0$ for some arbitrary fixed K\"ahler metric $\omega$ on $X$. Applying (\ref{eqn:prod-traces-direct}) to the K\"ahler metrics $\alpha_{\varepsilon} + \varepsilon\,\omega$ and $\beta_{\varepsilon} + \varepsilon\,\omega$ in place of $\alpha$ and $\beta$ and letting $\varepsilon\downarrow 0$, we get $\int_X\beta\wedge\gamma^{n-1}=0$ for every Gauduchon metric $\gamma$ on $X$. (Note that $\int_X\alpha_{\varepsilon}\wedge\gamma^{n-1} = \int_X\alpha\wedge\gamma^{n-1}$, $\int_X\beta_{\varepsilon}\wedge\gamma^{n-1} = \int_X\beta\wedge\gamma^{n-1}$ and $\int_X\alpha_{\varepsilon}^{n-1}\wedge\beta_{\varepsilon} = \{\alpha\}^{n-1}.\{\beta\} = 0$.) If we fix a $d$-closed positive current $T\geq 0$ in the class $\{\beta\}$ (such a current exists since the nef class $\{\beta\}$ is, in particular, pseudo-effective), this means that $\int_XT\wedge\gamma^{n-1}=0$ for every Gauduchon metric $\gamma$ on $X$. Consequently, $T=0$, hence $\{\beta\} = \{T\} = 0$, a contradiction.  \hfill $\Box$

\vspace{2ex}

 An immediate consequence of Corollary \ref{Cor:prod-traces-direct} is the following result in which the volume lower bound (\ref{eqn:volume_n2_lbound}) falls short of the expected inequality (\ref{eqn:volume_lbound}). However, (\ref{eqn:volume_n2_lbound}) solves the qualitative part of [BDPP13, Conjecture $10.1,$\, $(ii)$] already solved in [Pop14], while (\ref{eqn:T_n2_lbound}) gives moreover an effective estimate of the largest $t>0$ for which the class $\{\alpha - t\beta\}$ remains pseudo-effective. This estimate will prompt the discussion of the psef and nef thresholds in the next section.

\begin{The}\label{The:volume_n2_lbound}Let $X$ be a compact K\"ahler manifold with $\mbox{dim}_{\C}X=n$ and let $\alpha, \beta>0$ be K\"ahler metrics such that $\{\alpha\}^n - n\,\{\alpha\}^{n-1}.\,\{\beta\}>0$. 

 Then, for every $t\in[0, \,\,+\infty)$, there exists a real $(1,\,1)$-current $T_t\in\{\alpha - t\beta\}$ such that

\begin{equation}\label{eqn:T_n2_lbound}T_t\geq\bigg(1 - nt\,\frac{\{\alpha\}^{n-1}.\,\{\beta\}}{\{\alpha\}^n}\bigg)\,\alpha  \hspace{2ex}\mbox{on}\hspace{1ex}X.\end{equation}

\noindent In particular, $T_t$ is a K\"ahler current for all $0\leq t<\frac{\{\alpha\}^n}{n\,\{\alpha\}^{n-1}.\,\{\beta\}}$, so taking $t=1$ (which is allowed by the assumption $\{\alpha\}^n - n\,\{\alpha\}^{n-1}.\,\{\beta\}>0$) we get that the class $\{\alpha - \beta\}$ contains a K\"ahler current. Moreover, its volume satisfies:

\begin{equation}\label{eqn:volume_n2_lbound}\mbox{Vol}\,\,(\{\alpha - \beta\}) \geq (\{\alpha\}^n - n\,\{\alpha\}^{n-1}.\,\{\beta\})\,\bigg(\frac{\{\alpha\}^n - n\,\{\alpha\}^{n-1}.\,\{\beta\}}{\{\alpha\}^n}\bigg)^{n-1} \geq \{\alpha\}^n - n^2\,\{\alpha\}^{n-1}.\,\{\beta\}.\end{equation}

\end{The}

\noindent {\it Proof.} Thanks to Lamari's positivity criterion (Lemma \ref{Lem:Lamari}), the existence of a current $T_t\in\{\alpha - t\beta\}$ satisfying (\ref{eqn:T_n2_lbound}) is equivalent to

$$\int\limits_X\bigg(\alpha - t\,\beta - \alpha + nt\,\frac{\{\alpha\}^{n-1}.\,\{\beta\}}{\{\alpha\}^n}\,\alpha\bigg)\wedge\gamma^{n-1} \geq 0$$

\noindent for every Gauduchon metric $\gamma$ on $X$. This, in turn, is equivalent to

$$nt\,\frac{\{\alpha\}^{n-1}.\,\{\beta\}}{\{\alpha\}^n}\,\int\limits_X\alpha\wedge\gamma^{n-1}\geq t\,\int\limits_X\beta\wedge\gamma^{n-1} \hspace{2ex}\mbox{for every Gauduchon metric}\hspace{1ex}\gamma.$$

\noindent The last inequality is nothing but (\ref{eqn:prod-traces-direct}) which was proved in Corollary \ref{Cor:prod-traces-direct}. This completes the proof of the existence of a current $T_t\in\{\alpha-t\beta\}$ satisfying (\ref{eqn:T_n2_lbound}). 

 Now, (\ref{eqn:T_n2_lbound}) implies that the absolutely continuous part $T_{ac}$ of $T:=T_1\in\{\alpha - \beta\}$ has the same lower bound as $T$. Moreover, if $\{\alpha\}^n - n\,\{\alpha\}^{n-1}.\,\{\beta\}>0$, then 

\begin{equation}\nonumber\mbox{Vol}\,\,(\{\alpha - \beta\}) \geq \int\limits_XT_{ac}^n\geq\bigg(1 - n\,\frac{\{\alpha\}^{n-1}.\,\{\beta\}}{\{\alpha\}^n}\bigg)^n\,\int\limits_X\alpha^n \stackrel{(i)}{\geq} \bigg(1 - n^2\,\frac{\{\alpha\}^{n-1}.\,\{\beta\}}{\{\alpha\}^n}\bigg)\,\{\alpha\}^n,\end{equation}

\noindent which proves the claim (\ref{eqn:volume_n2_lbound}). To obtain $(i)$, we have used the elementary inequality $(1-\lambda)^n\geq 1-n\lambda$ which holds for every $\lambda\in[0,\,1]$.   \hfill $\Box$

\vspace{2ex}

 The above proof shows that a current $T_t\in\{\alpha - t\beta\}$ satisfying (\ref{eqn:T_n2_lbound}) exists even if we do not assume $\{\alpha\}^n - n\,\{\alpha\}^{n-1}.\,\{\beta\}>0$, although this information will be of use only under this assumption. Note that the non-orthogonality property $\{\alpha\}^{n-1}.\,\{\beta\}>0$ ensured by Corollary \ref{Cor:non-orthogonality} constitutes the obstruction to the volume lower estimate (\ref{eqn:volume_n2_lbound}) being optimal (i.e. coinciding with the expected estimate (\ref{eqn:volume_lbound})). We now point out an alternative way of inferring the same suboptimal volume lower bound (\ref{eqn:volume_n2_lbound}) from the proof of Theorem \ref{The:volume_n2_lbound}.

\vspace{2ex}

\noindent {\it Alternative wording of the proof of the volume lower estimate (\ref{eqn:volume_n2_lbound}).} By Lamari's positivity criterion (Lemma \ref{Lem:Lamari}), the existence of a current $T$ in the class $\{\alpha - \beta\}$ such that $T\geq\delta\alpha$ for some constant $\delta>0$ (which must be such that $\delta<1$) is equivalent to 

$$\int\limits_X\bigg((1-\delta)\,\alpha - \beta\bigg)\wedge\gamma^{n-1}\geq 0,  \hspace{2ex} \mbox{i.e. to}\hspace{2ex} \int\limits_X\bigg((\alpha-\frac{1}{1-\delta}\,\beta\bigg)\wedge\gamma^{n-1}\geq 0,$$

\noindent for all Gauduchon metrics $\gamma$ on $X$. Applying again Lamari's positivity criterion, this is still equivalent to the class $\{\alpha\} - (1/(1-\delta))\,\{\beta\}$ being pseudo-effective. Inequality (\ref{eqn:T_n2_lbound}) shows that the largest $\delta$ we can choose with this property is larger than or equal to

\begin{equation}\label{eqn:best-delta}\delta_0 = 1 - n\,\frac{\{\alpha\}^{n-1}.\{\beta\}}{\{\alpha\}^n}, \hspace{2ex} \mbox{which gives} \hspace{2ex} 1-\delta_0 = n\,\frac{\{\alpha\}^{n-1}.\{\beta\}}{\{\alpha\}^n}.\end{equation}

\noindent On the other hand, we can write:

\begin{equation}\label{eqn:vol_sum}\mbox{Vol}\,\,(\{\alpha - \beta\}) = \mbox{Vol}\,\,\bigg((1-t)\,\{\alpha\}  + t(\{\alpha\} - \frac{1}{t}\,\{\beta\})\bigg) \stackrel{(a)}{\geq} \mbox{Vol}\,\,((1-t)\,\{\alpha\}) = (1-t)^n\,\{\alpha\}^n,\end{equation}

\noindent where inequality $(a)$ holds for every $t\in [0,\,1]$ for which the class $\{\alpha\} - (1/t)\,\{\beta\}$ is {\it pseudo-effective}. By (\ref{eqn:best-delta}), $t: = 1-\delta_0 = n\,\{\alpha\}^{n-1}.\{\beta\}/\{\alpha\}^n$ satisfies this property. With this choice of $t$, inequality (\ref{eqn:vol_sum}) translates to the first inequality in (\ref{eqn:volume_n2_lbound}).  \hfill $\Box$

\vspace{2ex}

\begin{Cor}\label{Cor:big_nef_classes} Let $\{\alpha\}, \{\beta\}\in H^{1,\,1}_{BC}(X,\,\R)$ be {\bf nef} classes on a compact K\"ahler manifold $X$ with $\mbox{dim}_{\C}X=n$ such that $\{\alpha\}^n-n\,\{\alpha\}^{n-1}.\{\beta\}>0$. If $\{\beta\}=0$, then $\{\alpha\}$ is big, while if $\{\beta\}\neq 0$, then $\{\alpha\} -t\,\{\beta\}$ is big for all $0\leq t<\frac{\{\alpha\}^n}{n\,\{\alpha\}^{n-1}.\,\{\beta\}}$. Moreover, the volume lower bound (\ref{eqn:volume_n2_lbound}) holds.

\end{Cor}

 The case when $\{\beta\}=0$ is the key Theorem 2.12 in [DP04]. So, in particular, our method produces a much quicker proof of this fundamental result of [DP04]. The case when $\{\beta\}\neq 0$ is new, although the case $t=1$ and the method of proof are those of [Pop14]. Notice that the quantity $\{\alpha\}^n/n\,\{\alpha\}^{n-1}.\,\{\beta\}>0$ is well defined when $\{\beta\}\neq 0$ by Corollary \ref{Cor:non-orthogonality}.

\vspace{2ex}

\noindent {\it Proof.} We fix an arbitrary K\"ahler metric $\omega$ on $X$ and a constant $t\geq 0$ that will be specified shortly. The nefness assumption on $\{\alpha\}, \{\beta\}$ means that for every $\varepsilon>0$, smooth forms $\alpha\in\{\alpha\}$ and $\beta\in\{\beta\}$ depending on $\varepsilon$ can be found such that $\alpha_{\varepsilon}:=\alpha + \varepsilon\,\omega$ and $\beta_{\varepsilon}:=\beta + \frac{\varepsilon}{t}\,\omega$ are K\"ahler metrics. Notice that the class $\{\alpha_{\varepsilon} -t\beta_{\varepsilon}\} = \{\alpha -t\beta\}$ is independent of $\varepsilon$. On the other hand, the quantities $\{\alpha_{\varepsilon}\}^n = \{\alpha\}^n + \sum\limits_{k=1}^n\varepsilon^k\,{n \choose k}\,\{\alpha\}^{n-k}.\{\omega\}^k$ and $\{\alpha_{\varepsilon}\}^{n-1}.\,\{\beta_{\varepsilon}\} = (\{\alpha\}^{n-1} + \sum\limits_{l=1}^{n-1}\varepsilon^l\,{n-1 \choose l}\,\{\alpha\}^{n-1-l}.\{\omega\}^l).(\{\beta\} + \frac{\varepsilon}{t}\,\{\omega\})$ converge to $\{\alpha\}^n$ and respectively $\{\alpha\}^{n-1}.\{\beta\}$ when $\varepsilon\rightarrow 0$. Thus, $\{\alpha_{\varepsilon}\}^n-n\,\{\alpha_{\varepsilon}\}^{n-1}.\{\beta_{\varepsilon}\}>0$ if $\varepsilon>0$ is small enough. Applying Theorem \ref{The:volume_n2_lbound} to the K\"ahler metrics $\alpha_{\varepsilon}$ and $\beta_{\varepsilon}$, we infer that the class $\{\alpha_{\varepsilon} -t\beta_{\varepsilon}\} = \{\alpha -t\beta\}$ is big whenever $0\leq t<\{\alpha_{\varepsilon}\}^n/n\,\{\alpha_{\varepsilon}\}^{n-1}.\,\{\beta_{\varepsilon}\}$. In particular, if $\{\beta\}= 0$, this means that the class $\{\alpha\}$ is big (since we can fix $\varepsilon>0$ and choose $t=0$). Meanwhile, if $\{\beta\}\neq 0$ and if we choose $t<\{\alpha\}^n/n\,\{\alpha\}^{n-1}.\{\beta\}$, then $t<\{\alpha_{\varepsilon}\}^n/n\,\{\alpha_{\varepsilon}\}^{n-1}.\{\beta_{\varepsilon}\}$ for all $\varepsilon>0$ small enough and we conclude that $\{\alpha -t\beta\}$ is big. The volume lower bound (\ref{eqn:volume_n2_lbound}) holds for $\{\alpha_{\varepsilon}\}$ and $\{\beta_{\varepsilon}\}$ for $t=1$ and all sufficiently small $\varepsilon>0$, so letting $\varepsilon\rightarrow 0$ and using the continuity of the volume, we get it for $\{\alpha\}$ and $\{\beta\}$.   \hfill $\Box$

\section{Trace and volume of $(1,\,1)$-cohomology classes}\label{section:dir-trace}

 The implicit discussion of the relative positivity thresholds of a cohomology class with respect to another in Theorem \ref{The:volume_n2_lbound} and in Corollary \ref{Cor:big_nef_classes} prompts a further investigation of their relationships with the volume that we undertake to study in this section.

\subsection{The psef threshold}\label{subsection:psef}

 Let $X$ be a compact complex manifold in Fujiki's {\it class} ${\cal C}$, $n:=\mbox{dim}_{\C}X$. 

\begin{Def}\label{Def:dir-trace}  For every {\bf big} Bott-Chern class $\{\beta\} = [\beta]_{BC}\in H^{1,\,1}_{BC}(X,\,\R)$, we define the {\bf $\beta$-directed trace} (or the {\bf psef threshold in the $\beta$-direction}) to be the function:

\begin{equation}\label{eqn:dir-trace}P^{(\beta)}\,\,:\,\,H^{1,\,1}_{BC}(X,\,\R)\rightarrow\R, \hspace{2ex} P^{(\beta)}(\alpha):=\inf\int\limits_X\alpha\wedge\gamma^{n-1},\end{equation}

\noindent for all Bott-Chern classes $\{\alpha\} = [\alpha]_{BC}\in H^{1,\,1}_{BC}(X,\,\R)$, where the infimum is taken over all the Gauduchon metrics $\gamma$ on $X$ normalised such that

\begin{equation}\label{eqn:dir-trace-normalisation}[\beta]_{BC}.[\gamma^{n-1}]_A=\int\limits_X\beta\wedge\gamma^{n-1}=1.\end{equation}

\end{Def}

 All the integrals involved in the above definition are clearly independent of the representatives $\alpha$, $\beta$ of the Bott-Chern classes $[\alpha]_{BC} ,[\beta]_{BC}$ and of the representative $\gamma^{n-1}$ of the Aeppli-Gauduchon class $[\gamma^{n-1}]_A\in H^{n-1,\,n-1}_A(X,\,\R)$. Thus the infimum is taken over the subset $S_{\beta}$ of the Gauduchon cone ${\cal G}_X$ consisting of classes $[\gamma^{n-1}]_A$ normalised by $[\beta]_{BC}.[\gamma^{n-1}]_A=1$. The bigness assumption on $[\beta]_{BC}$ has been imposed to ensure that $[\beta]_{BC}.[\gamma^{n-1}]_A>0$, hence that $[\gamma^{n-1}]_A$ can be normalised w.r.t. $[\beta]_{BC}$ as in (\ref{eqn:dir-trace-normalisation}), for {\it every} class $[\gamma^{n-1}]_A\in{\cal G}_X$. 

 This definition is motivated in part by the next observation which is an immediate consequence of Lamari's positivity criterion: the ${\beta}$-directed trace $P^{(\beta)}$ coincides with the slope function introduced for big classes $\{\alpha\}$ in [BFJ09, Definition 3.7] and thus gauges the positivity of real $(1,\,1)$-classes $\{\alpha\}$ w.r.t. a reference big class $\{\beta\}$. The quantity on the r.h.s. of (\ref{eqn:trace-psef}) below (i.e. the slope) may well be called the {\bf psef threshold of $\{\alpha\}$ in the $\{\beta\}$-direction} (a term already used in the literature).

\begin{Prop}\label{Prop:trace-slope} Suppose that $\{\beta\}\in H^{1,\,1}_{BC}(X,\,\R)$ is a fixed {\bf big} class. Then

\begin{eqnarray}\label{eqn:trace-psef}P^{(\beta)}(\alpha) & = & \sup\,\{t\in\R\,\,\slash\,\,\mbox{the class}\,\,\{\alpha\} - t\,\{\beta\}\,\,\mbox{is psef}\}\\
\nonumber  & = & \sup\,\{t\in\R\,\,\slash\,\,\exists\, T\in\{\alpha\}\,\,\mbox{current,}\,\,\exists\, \widetilde{\beta}\in\{\beta\}\,\,C^{\infty}\mbox{-form s.t.}\,\,T\geq t\,\widetilde{\beta}\}\\
\nonumber & = & \sup\,\{t\in\R\,\,\slash\,\,\forall\widetilde{\beta}\in\{\beta\}\,\,C^{\infty}\mbox{-form,}\,\,\exists\, T\in\{\alpha\}\,\,\mbox{current s.t.}\,\,T\geq t\,\widetilde{\beta}\},\end{eqnarray}

\noindent for every class $\{\alpha\}\in H^{1,\,1}_{BC}(X,\,\R)$. In particular, the set $\{t\in\R\,\slash\,\mbox{the class}\,\,\{\alpha\} - t\,\{\beta\}\,\mbox{is psef}\}$ equals the interval $(-\infty,\,P^{(\beta)}(\alpha)]$.

\end{Prop}

\noindent {\it Proof.} Let $A_{\alpha}^{\beta}:= \{t\in\R\,\,\slash\,\,\mbox{the class}\,\,\{\alpha\} - t\,\{\beta\}\,\,\mbox{is psef}\}$ and let $t_{\alpha}^{\beta}:=\sup\,A_{\alpha}^{\beta}$. By Lamari's positivity criterion, the class $\{\alpha\} - t\,\{\beta\}$ is psef iff 

$$\int\limits_X\alpha\wedge\gamma^{n-1}\geq t\int\limits_X\beta\wedge\gamma^{n-1} \hspace{2ex} \mbox{for all}\,\,[\gamma^{n-1}]_A\in{\cal G}_X \iff \int\limits_X\alpha\wedge\gamma^{n-1}\geq t$$

\noindent for all $[\gamma^{n-1}]_A\in{\cal G}_X$ normalised such that $\int\limits_X\beta\wedge\gamma^{n-1}=1$. This proves the inequality $P^{(\beta)}(\alpha)\geq t_{\alpha}^{\beta}$. To prove that equality holds, we reason by contradiction. Suppose that $P^{(\beta)}(\alpha)>t_{\alpha}^{\beta}$. Pick any $t_1$ such that $P^{(\beta)}(\alpha)>t_1>t_{\alpha}^{\beta}$. Then $\int_X\alpha\wedge\gamma^{n-1}>t_1$ for all Gauduchon metrics $\gamma$ such that $[\beta]_{BC}.[\gamma^{n-1}]_A=1$. This is equivalent to $\int_X\alpha\wedge\gamma^{n-1}>t_1\,\int_X\beta\wedge\gamma^{n-1}$ for all Gauduchon metrics $\gamma$, which thanks to Lamari's positivity criterion implies:

\vspace{1ex}

\noindent $\exists\,T\in\{\alpha\} - t_1\,\{\beta\} \hspace{2ex}\mbox{such that}\hspace{2ex}T\geq 0, \hspace{2ex}\mbox{i.e. the class}\hspace{1ex}\{\alpha\} - t_1\,\{\beta\}\hspace{1ex}\mbox{is psef}.$

\vspace{1ex}

\noindent Thus $t_1\in A_{\alpha}^{\beta}$, contradicting the choice $t_1>t_{\alpha}^{\beta} = \sup\,A_{\alpha}^{\beta}$.  \hfill $\Box$

\vspace{2ex}

 An immediate consequence is the next statement showing that the ${\beta}$-directed trace (= the psef threshold) gauges the positivity of real Bott-Chern $(1,\,1)$-classes much as the volume does.

\begin{Cor}\label{Cor:beta-dir-positivity} $(i)$\, Suppose that $\{\beta\}\in H^{1,\,1}_{BC}(X,\,\R)$ is a fixed {\bf big} class. For any class $\{\alpha\}\in H^{1,\,1}_{BC}(X,\,\R)$, the following equivalences hold:\\

 $(i)$\, $\{\alpha\}$ is psef $\iff$ $P^{(\beta)}(\alpha)\geq 0.$

\vspace{1ex}

 $(ii)$\,  $\{\alpha\}$ is big $\iff$ $P^{(\beta)}(\alpha)> 0.$

\end{Cor}

\noindent {\it Proof.} $(i)$ follows at once from (\ref{eqn:trace-psef}) and so does $(ii)$ after we (trivially) notice that the class $\{\alpha\}$ is big iff there exists $\varepsilon >0$ such that $\{\alpha\} - \varepsilon\{\beta\}$ is psef. Indeed, this is a consequence of the fixed class $\{\beta\}$ being supposed big.     \hfill $\Box$

\vspace{2ex}

 Next, we observe some easy but useful properties of the $\beta$-directed trace.

\begin{Prop}\label{Prop:dir-trace-properties}Suppose that $\{\beta\}\in H^{1,\,1}_{BC}(X,\,\R)$ is a fixed {\bf big} class.

\noindent $(i)$\, For all classes $\{\alpha_1\}, \{\alpha_2\}\in H^{1,\,1}_{BC}(X,\,\R)$, we have

\begin{equation}\label{eqn:super-add}P^{(\beta)}(\alpha_1 + \alpha_2) \geq P^{(\beta)}(\alpha_1) + P^{(\beta)}(\alpha_2).\end{equation}

\noindent In particular, $P^{(\beta)}(\alpha_1)\geq P^{(\beta)}(\alpha_2)$ whenever $\{\alpha_1\}\geq_{psef}\{\alpha_2\}$ (in the sense that $\{\alpha_1-\alpha_2\}$ is psef).

\noindent $(ii)$\, For any class $\{\alpha\}\in H^{1,\,1}_{BC}(X,\,\R)$ and any $t\in\R$, we have

\begin{equation}\label{eqn:homogeneity}P^{(\beta)}(t\,\alpha) = t\, P^{(\beta)}(\alpha) \hspace{2ex}\mbox{and, if $t>0$,}\hspace{2ex} P^{(t\,\beta)}(\alpha) = \frac{1}{t}\, P^{(\beta)}(\alpha).\end{equation}

\noindent $(iii)$\, For every {\bf big} class $\{\alpha\}\in H^{1,\,1}_{BC}(X,\,\R)$, we have

\begin{equation}\label{eqn:Palpha-alpha}P^{(\alpha)}(\alpha) = 1.\end{equation}

\end{Prop}

\noindent {\it Proof.} Let $\{\alpha_1\}, \{\alpha_2\}\in H^{1,\,1}_{BC}(X,\,\R)$. Since $\int_X(\alpha_1 + \alpha_2)\wedge\gamma^{n-1} = \int_X\alpha_1\wedge\gamma^{n-1} + \int_X\alpha_2\wedge\gamma^{n-1}$ for every $[\gamma^{n-1}]_A\in H^{n-1,\,n-1}_A(X,\,\R)$, we get 

\vspace{1ex}

\hspace{6ex} $\inf\int_X(\alpha_1 + \alpha_2)\wedge\gamma^{n-1} \geq \inf\int_X\alpha_1\wedge\gamma^{n-1} + \inf\int_X\alpha_1\wedge\gamma^{n-1},$

\vspace{1ex}

\noindent where the infima are taken over all $[\gamma^{n-1}]\in S_{\beta}$. This proves $(i)$. 

\noindent $(ii)$ follows immediately from $\int_Xt\alpha\wedge\gamma^{n-1}=t\,\int_X\alpha\wedge\gamma^{n-1}$ and from the fact that $[\gamma^{n-1}]_A$ is $(t\beta)$-normalised if and only if $t[\gamma^{n-1}]_A$ is $\beta$-normalised.

\noindent $(iii)$ follows from $\int_X\alpha\wedge\gamma^{n-1}=1$ for all $[\gamma^{n-1}]_A$ such that $[\alpha]_{BC}.[\gamma^{n-1}]_A=1$. \hfill $\Box$

\vspace{2ex}

 The next observation deals with the variation of $P^{(\beta)}$ when $\{\beta\}$ varies. As usual, an inequality $\{\alpha\}\geq_{psef}\{\beta\}$ between real $(1,\,1)$-classes will mean that the difference class $\{\alpha-\beta\}$ is psef.

\begin{Prop}\label{Prop:Pbeta-variation} Let $\{\beta_1\}, \{\beta_2\}\in H^{1,\,1}_{BC}(X,\,\R)$ be {\bf big} classes.

\noindent $(i)$\, If $\{\beta_1\}\geq_{psef}C\,\{\beta_2\}$ for some constant $C>0$, then

\begin{equation}\label{eqn:Pbeta-variation1}P^{(\beta_1)}\leq\frac{1}{C}\,P^{(\beta_2)}   \hspace{2ex} \mbox{on the psef cone} \hspace{1ex} {\cal E}_X\subset H^{1,\,1}_{BC}(X,\,\R).\end{equation}

\noindent $(ii)$\, The following inequality holds:

\begin{equation}\label{eqn:Pbeta-variation2}P^{(\beta_2)}(\beta_1)\,P^{(\beta_1)}\leq
P^{(\beta_2)}    \hspace{2ex} \mbox{on the psef cone} \hspace{1ex} {\cal E}_X\subset H^{1,\,1}_{BC}(X,\,\R).\end{equation}

\end{Prop}

\noindent {\it Proof.} If $\{\beta_1-C\,\beta_2\}$ is psef, then $\int_X(\beta_1-C\,\beta_2)\wedge\gamma^{n-1}\geq 0$, i.e. $[\beta_1]_{BC}.[\gamma^{n-1}]_A\geq C\,[\beta_2]_{BC}.[\gamma^{n-1}]_A$, for all classes $[\gamma^{n-1}]_A\in{\cal G}_X$. It follows that, for every psef class $\{\alpha\}\in H^{1,\,1}_{BC}(X,\,\R)$, we have:

$$\int\limits_X\alpha\wedge\frac{\gamma^{n-1}}{\int\limits_X\beta_1\wedge\gamma^{n-1}} \leq \frac{1}{C}\, \int\limits_X\alpha\wedge\frac{\gamma^{n-1}}{\int\limits_X\beta_2\wedge\gamma^{n-1}} \hspace{2ex}\mbox{for all}\hspace{1ex} [\gamma^{n-1}]_A\in{\cal G}_X.$$

\noindent Taking infima over all $[\gamma^{n-1}]_A\in{\cal G}_X$, we get (\ref{eqn:Pbeta-variation1}). On the other hand, it follows from (\ref{eqn:trace-psef}) that

$$\{\beta_1\}\geq_{psef}P^{(\beta_2)}(\beta_1)\,\{\beta_2\},$$

\noindent which in turn implies (\ref{eqn:Pbeta-variation2}) thanks to (\ref{eqn:Pbeta-variation1}) applied with $C=P^{(\beta_2)}(\beta_1)$. \hfill $\Box$  

\subsection{The nef threshold}\label{subsection:nef}

 We now observe that the discussion of the psef threshold in $\S.$\ref{subsection:psef} can be run analogously in the nef context using the following important result of [DP04, Corollary 0.4].

\begin{The}(Demailly-Paun 2004)\label{The:DP04} Let $X$ be a compact K\"ahler manifold, $\mbox{dim}_{\C}X=n$. Then the dual of the nef cone $\overline{{\cal K}_X}\subset H^{1,\,1}(X,\,\R)$ under the Serre duality is the closed convex cone ${\cal N}_X\subset H^{n-1,\,n-1}(X,\,\R)$ generated by classes of currents of the shape $[Y]\wedge\omega^{n-p-1}$, where $Y$ runs over the irreducible analytic subsets of $X$ of any codimension $p=0, 1, \dots , n-1$ and $\{\omega\}$ runs over the K\"ahler classes of $X$.  

\end{The}

 Let $\{\alpha\}, \{\beta\}\in H^{1,\,1}(X,\,\R)$ be arbitrary classes on a compact K\"ahler $n$-fold $X$. By Theorem \ref{The:DP04}, for any $s\in\R$, the class $\{\alpha - s\beta\}$ is nef iff

$$\int\limits_Y\alpha\wedge\omega^{n-p-1}\geq s\,\int\limits_Y\beta\wedge\omega^{n-p-1},  \hspace{2ex}p=0,1, \dots , n-1,\, \mbox{codim}_X\,Y=p, \{\omega\}\in{\cal K}_X.$$

\noindent (As usual, ${\cal K}_X$ denotes the K\"ahler cone of $X$.) This immediately implies the following statement.

\begin{Prop}\label{Prop:nef-threshold_def} Let $\{\beta\}\in H^{1,\,1}_{BC}(X,\,\R)$ be any {\bf K\"ahler} class on a compact K\"ahler $n$-fold $X$. The {\bf nef threshold} of any $\{\alpha\}\in H^{1,\,1}_{BC}(X,\,\R)$ {\bf in the $\{\beta\}$-direction}, defined by the first identity below, also satisfies the second identity:

 \begin{equation}\label{eqn:trace-nef} N^{(\beta)}(\alpha):=\inf\int\limits_Y\alpha\wedge\omega^{n-p-1} = \sup\{s\in\R\,\,/\,\,\mbox{the class}\,\,\{\alpha\}-s\{\beta\}\,\,\mbox{is nef}\,\},\end{equation}

\noindent where the infimum is taken over all $p=0,1, \dots , n-1$, over all the irreducible analytic subsets $Y\subset X$ such that $\mbox{codim}\,Y=p$ and over all K\"ahler classes $\{\omega\}$ normalised such that $\int_Y\beta\wedge\omega^{n-p-1}=1$. In particular, the set $\{s\in\R\,\,/\,\,\mbox{the class}\,\,\{\alpha\}-s\{\beta\}\,\,\mbox{is nef}\,\}$ equals the interval $(-\infty,\,N^{(\beta)}(\alpha)]$.

 Thus, we obtain a function $N^{(\beta)}:H^{1,\,1}_{BC}(X,\,\R)\rightarrow\R$. It is clear that

\begin{equation}\label{equation:psef-nef_ineq}N^{(\beta)}(\alpha)\leq P^{(\beta)}(\alpha) \hspace{3ex} \mbox{for all}\hspace{1ex}\{\alpha\}\in H^{1,\,1}_{BC}(X,\,\R)\end{equation}

\noindent thanks to the supremum characterisations of the two thresholds and to the well-known implication ``$\mbox{nef}\implies\mbox{psef}$''.

\end{Prop}

 It is precisely in order to ensure that $\int_Y\beta\wedge\omega^{n-p-1}>0$, hence that $\{\omega\}$ can be normalised as stated, for any K\"ahler class $\{\omega\}$ and any $Y\subset X$ that we assumed $\{\beta\}$ to be K\"ahler. The two-fold characterisations of the nef and the psef thresholds yield at once the following consequence.

\begin{Obs}\label{Obs:P=N} Suppose that no analytic subset $Y\subset X$ exists except in codimensions $0$ and $n$. Then $N^{(\beta)}(\alpha)= P^{(\beta)}(\alpha)$ for all K\"ahler classes $\{\alpha\}, \{\beta\}.$ 

\end{Obs}

\noindent {\it Proof.} If $Y=X$ is the only analytic subset of $X$ of codimension $p<n$, then $N^{(\beta)}(\alpha)=\inf\int_X\alpha\wedge\omega^{n-1}$ where the infimum is taken over all the K\"ahler classes $\{\omega\}$, i.e. over all the Aeppli-Gauduchon classes $[\omega^{n-1}]_A$ representable by the $(n-1)^{st}$ power of a K\"ahler metric, normalised such that $\int_X\beta\wedge\omega^{n-1}=1$. Since these classes form a subset of all the Aeppli-Gauduchon classes $[\gamma^{n-1}]_A$ normalised by $\int_X\beta\wedge\gamma^{n-1}=1$, we get $N^{(\beta)}(\alpha)\geq P^{(\beta)}(\alpha)$. However, the reverse inequality always holds, hence equality holds.  \hfill $\Box$

\vspace{2ex}

 An immediate consequence of Proposition \ref{Prop:nef-threshold_def} is the following analogue of Corollary \ref{Cor:beta-dir-positivity} for the nef/K\"ahler context.

\begin{Cor}\label{Cor:beta-dir-positivity_nef} $(i)$\, Suppose that $\{\beta\}\in H^{1,\,1}_{BC}(X,\,\R)$ is a fixed {\bf K\"ahler} class. For any class $\{\alpha\}\in H^{1,\,1}_{BC}(X,\,\R)$, the following equivalences hold:\\

 $(i)$\, $\{\alpha\}$ is nef $\iff$ $N^{(\beta)}(\alpha)\geq 0.$

\vspace{1ex}

 $(ii)$\,  $\{\alpha\}$ is K\"ahler $\iff$ $N^{(\beta)}(\alpha)> 0.$

\vspace{1ex}

\noindent In particular, if no analytic subset $Y\subset X$ exists except in codimensions $0$ and $n$, then the following (actually known, see [Dem92]) equivalences hold:

\vspace{1ex}

    $(a)\hspace{2ex} \{\alpha\}\,\mbox{is nef}\, \iff \{\alpha\} \mbox{is psef}, \hspace{2ex} (b)\hspace{2ex}  \{\alpha\}\, \mbox{is K\"ahler}\, \iff \{\alpha\}\, \mbox{is big}.$

\end{Cor}

\noindent {\it Proof.} $(i)$ follows at once from (\ref{eqn:trace-nef}) and so does $(ii)$ after we (trivially) notice that the class $\{\alpha\}$ is K\"ahler iff there exists $\varepsilon >0$ such that $\{\alpha\} - \varepsilon\{\beta\}$ is nef. Indeed, this is a consequence of the fixed class $\{\beta\}$ being supposed K\"ahler and of the K\"ahler cone being the interior of the nef cone. \hfill $\Box$

\vspace{2ex}

 We immediately get analogues of Propositios \ref{Prop:dir-trace-properties} and \ref{Prop:Pbeta-variation} for $N^{(\beta)}(\alpha)$ in place of $P^{(\beta)}(\alpha)$ and for the order relation $\geq_{nef}$ in place of $\geq_{psef}$, where $\{\alpha\}\geq_{nef}\{\beta\}$ means that the class $\{\alpha-\beta\}$ is nef.

\subsection{Relations of the psef/nef threshold to the volume}\label{subsection:trace-volume}

 We now relate the ${\beta}$-directed trace of a K\"ahler class $\{\alpha\}$ to the volume of $\{\alpha-\beta\}$.

\begin{Prop}\label{Prop:P-vol}$(i)$\, For every {\bf K\"ahler} classes $\{\alpha\}, \{\beta\}$ on a compact K\"ahler $n$-fold $X$, we have:

\begin{equation}\label{eqn:P-vol1}\frac{\{\alpha\}^n}{n\,\{\alpha\}^{n-1}.\{\beta\}} \stackrel{(a)}{\leq} P^{(\beta)}(\alpha) \stackrel{(b)}{\leq} \frac{\{\alpha\}^n}{\{\alpha\}^{n-1}.\{\beta\}}.\end{equation}

\noindent In fact, it suffices to suppose that $\{\beta\}$ is big in the inequality $(b)$. In particular, if $\{\alpha\}^n - n\,\{\alpha\}^{n-1}.\{\beta\}>0$, then $P^{(\beta)}(\alpha)>1$ (hence we find again that $\{\alpha-\beta\}$ is big in this case).

\vspace{1ex}

 \noindent $(ii)$\, For every {\bf K\"ahler} classes $\{\alpha\}, \{\beta\}$ such that $\{\alpha\}^n - n\,\{\alpha\}^{n-1}.\{\beta\}>0$, we have:

\begin{equation}\label{eqn:P-vol2}\mbox{Vol}\,(\{\alpha-\beta\})\geq\bigg(1 - \frac{1}{P^{(\beta)}(\alpha)}\bigg)^n\,\{\alpha\}^n.\end{equation}

\noindent Note that the combination of (\ref{eqn:P-vol2}) and part $(a)$ of (\ref{eqn:P-vol1}) is the volume lower bound (\ref{eqn:volume_n2_lbound}).

\end{Prop}

\noindent {\it Proof.} $(i)$\, Inequality $(b)$ is trivial: it suffices to choose $[\gamma^{n-1}]_A = t\,[\alpha^{n-1}]_A$ for the constant $t>0$ satisfying the $\beta$-normalisation condition $[\beta]_{BC}.t\,[\alpha^{n-1}]_A=1$, i.e. $t=1/\{\alpha\}^{n-1}.\{\beta\}$, and to use the definition of $P^{(\beta)}(\alpha)$ as an infimum.

 Inequality $(a)$ follows from Corollary \ref{Cor:prod-traces-direct} by taking the infimum over all the Gauduchon metrics $\gamma$ normalised by $\int_X\beta\wedge\gamma^{n-1}=1$ in (\ref{eqn:prod-traces-direct}). 

 $(ii)$\, We saw in the second proof of the lower estimate (\ref{eqn:volume_n2_lbound}) that (\ref{eqn:vol_sum}) holds for every $t\in [0,\,1]$ such that $\{\alpha\} - (1/t)\,\{\beta\}$ is psef. Now, (\ref{eqn:trace-psef}) shows that the infimum of all these $t$ is $1/P^{(\beta)}(\alpha)$. Thus (\ref{eqn:vol_sum}) holds for $t=1/P^{(\beta)}(\alpha)$, yielding (\ref{eqn:P-vol2}).   \hfill $\Box$

\vspace{2ex}

 A similar link between the volume and the nef threshold is given in the next result by considering Monge-Amp\`ere equations on analytic subsets $Y\subset X$. 

\begin{Prop}\label{Prop:N-vol} For every {\bf K\"ahler} classes $\{\alpha\}, \{\beta\}$ on a compact K\"ahler $n$-fold $X$, we have:

\begin{equation}\label{eqn:N-vol1}\inf_{\stackrel{p=0,1,\dots , n-1,}{Y\subset X,\,\mbox{codim}\,Y=p}}\frac{\mbox{Vol}_Y(\alpha)}{(n-p)\,\{\alpha\}^{n-p-1}.\{\beta\}.\{[Y]\}} \stackrel{(a)}{\leq} N^{(\beta)}(\alpha) \stackrel{(b)}{\leq}\inf_{\stackrel{p=0,1,\dots , n-1,}{Y\subset X,\,\mbox{codim}\,Y=p}}\frac{\mbox{Vol}_Y(\alpha)}{\{\alpha\}^{n-p-1}.\{\beta\}.\{[Y]\}},\end{equation}

\noindent where the infima are taken over the analytic subsets $Y\subset X$. We have set $\mbox{Vol}_Y(\alpha):=\int_Y\alpha^{n-p} = \int_X\alpha^{n-p}\wedge[Y]$ and ${\{\alpha\}^{n-p-1}.\{\beta\}.\{[Y]\}}:=\int_Y\alpha^{n-p-1}\wedge\beta =\int_X\alpha^{n-p-1}\wedge\beta\wedge[Y]$ (both quantities depending only on $p$ and the classes $\{\alpha\}, \{\beta\}, \{[Y]\}$).

\end{Prop}

\noindent {\it Proof.} Pick any K\"ahler metrics $\alpha\in\{\alpha\}$ and $\beta\in\{\beta\}$. Let $Y\subset X$ be any analytic subset of arbitrary codimension $p\in\{0,1, \dots , n-1\}$ and let $\omega$ be any K\"ahler metric on $X$ normalised such that $\int_Y\beta\wedge\omega^{n-p-1}=1$. We can solve the following Monge-Amp\`ere equation: 

\begin{equation}\label{eqn:MA-Y}\widetilde{\alpha}_Y^{n-p} = \mbox{Vol}_Y(\alpha)\,\beta\wedge\omega^{n-p-1}  \hspace{3ex} \mbox{on}\hspace{1ex}Y\end{equation}

\noindent in the sense that there exists a $d$-closed (weakly) positive $(1,\,1)$-current $\widetilde{\alpha}_Y$ on $Y$ (cf. Definition $1.2$ in [Dem85]) lying in the restricted class $\{\alpha\}_{|Y}$ such that $\widetilde{\alpha}_Y$ is $C^{\infty}$ on the regular part $Y_{reg}$ of $Y$. We defer to the end of the proof the explanation of how this follows from results in the literature. We adopt the standard point of view (see [Dem85, $\S.1$]) according to which $C^{\infty}$ forms on a singular variety $Y$ are defined locally as restrictions to $Y_{reg}$ of $C^{\infty}$ forms on an open subset of some $\C^N$ into which $Y$ locally embeds. In what follows, the exterior powers and products involving $\widetilde{\alpha}_Y$ are to be understood on $Y_{reg}$ even when we write $Y$.

 If we define $\mbox{det}_{\omega}\widetilde{\alpha}_Y$ by requiring $\widetilde{\alpha}_Y^{n-p} = (\mbox{det}_{\omega}\widetilde{\alpha}_Y)\,\omega^{n-p}$ on $Y$, then (\ref{eqn:MA-Y}) translates to the identity:

\begin{equation}\label{eqn:MA-Y_detform}\mbox{det}_{\omega}\widetilde{\alpha}_Y = \frac{\mbox{Vol}_Y(\alpha)}{n-p}\,\Lambda_{\omega}\beta_{|Y} \hspace{3ex} \mbox{on}\hspace{1ex}Y.\end{equation}

Thus, the argument in the proof of Lemma \ref{Lem:prod-traces-int} can be rerun on $Y$ as follows:

\vspace{2ex}

\noindent $\displaystyle\bigg(\int_Y\widetilde{\alpha}_Y\wedge\omega^{n-p-1}\bigg)\,\bigg(\int\limits_Y\widetilde{\alpha}_Y^{n-p-1}\wedge\beta\bigg) = \frac{1}{(n-p)^2}\bigg(\int_Y(\Lambda_{\omega}\widetilde{\alpha}_Y)\,\omega^{n-p}\bigg)\,\bigg(\int_Y(\Lambda_{\widetilde{\alpha}_Y}\beta_{|Y})\,(\mbox{det}_{\omega}\widetilde{\alpha}_Y)\,\omega^{n-p}\bigg)$

\noindent \begin{eqnarray}\nonumber & \stackrel{(a)}{\geq} & \frac{1}{(n-p)^2}\,\bigg(\int\limits_Y[(\Lambda_{\omega}\widetilde{\alpha}_Y)\,(\Lambda_{\widetilde{\alpha}_Y}\beta_{|Y})]^{\frac{1}{2}}\,(\mbox{det}_{\omega}\widetilde{\alpha}_Y)^{\frac{1}{2}}\,\omega^{n-p}\bigg)^2\\\nonumber  & \stackrel{(b)}{\geq} & \frac{1}{(n-p)^2}\,\bigg(\int\limits_Y(\Lambda_{\omega}\beta_{|Y})^{\frac{1}{2}}\,\bigg(\frac{\mbox{Vol}_Y(\alpha)}{n-p}\bigg)^{\frac{1}{2}}\,(\Lambda_{\omega}\beta_{|Y})^{\frac{1}{2}}\,\omega^{n-p}\bigg)^2\\
\nonumber & = & \frac{\mbox{Vol}_Y(\alpha)}{n-p}\,\bigg(\frac{1}{(n-p)}\,\int\limits_Y(\Lambda_{\omega}\beta_{|Y})\,\omega^{n-p}\bigg)^2 = \frac{\mbox{Vol}_Y(\alpha)}{n-p}\,\bigg(\int\limits_Y\beta\wedge\omega^{n-p-1}\bigg)^2 \stackrel{(c)}{=} \frac{\mbox{Vol}_Y(\alpha)}{n-p},\end{eqnarray}

\noindent where $(a)$ is an application of the Cauchy-Schwarz inequality, $(b)$ follows from the pointwise inequality $(\Lambda_{\omega}\widetilde{\alpha}_Y)\,(\Lambda_{\widetilde{\alpha}_Y}\beta)\geq\Lambda_{\omega}\beta$ (cf. [Pop14, Lemma 3.1]) and from (\ref{eqn:MA-Y_detform}), while $(c)$ follows from the normalisation $\int_X\beta\wedge\omega^{n-p-1}=1$. 

 Thus, since $\int_Y\widetilde{\alpha}_Y\wedge\omega^{n-p-1} = \int_Y\alpha\wedge\omega^{n-p-1}$ and $\int\limits_Y\widetilde{\alpha}_Y^{n-p-1}\wedge\beta = \int\limits_Y\alpha^{n-p-1}\wedge\beta$, we get:

$$\int_Y\alpha\wedge\omega^{n-p-1}\geq\frac{\mbox{Vol}_Y(\alpha)}{(n-p)\,\{\alpha\}^{n-p-1}.\{\beta\}.\{[Y]\}}$$

\noindent for every analytic subset $Y\subset X$ and every K\"ahler metric $\omega$ normalised by $\int_Y\beta\wedge\omega^{n-p-1}=1$. This proves inequality $(a)$ in (\ref{eqn:N-vol1}). 

 The proof of inequality $(b)$ in (\ref{eqn:N-vol1}) follows immediately by choosing the K\"ahler metric $\omega$ to be proportional to $\alpha$, i.e. $\omega = t\alpha$ for the constant $t = t_Y>0$ determined by the normalisation condition $\int_Y\beta\wedge\omega^{n-p-1}=1$ once $Y\subset X$ has been chosen. Indeed, for every $p=0, 1, \dots , n-1$ and every analytic subset $Y\subset X$, we immediately get: 

$$\inf\limits_{\omega}\int_Y\alpha\wedge\omega^{n-p-1}\leq\int_Y\alpha\wedge(t\alpha)^{n-p-1} = \frac{\int\limits_Y\alpha^{n-p}}{\int\limits_Y\beta\wedge\alpha^{n-p-1}}$$ 

\noindent which implies part $(b)$ of (\ref{eqn:N-vol1}) after taking the infimum over $p$ and $Y$.  

\vspace{2ex}

 It remains to explain how the solution of equation (\ref{eqn:MA-Y}) is obtained. If $Y$ is smooth, Yau's classical theorem in [Yau78] ensures the existence and uniqueness of a K\"ahler metric $\widetilde{\alpha}_Y$ in $\{\alpha\}_{|Y}$ which solves (\ref{eqn:MA-Y}). If $Y$ is singular, we choose a desingularisation $\widetilde{Y}$ of $Y$ that is a finite sequence of blow-ups with smooth centres in $X$:

$$\mu : \widetilde{Y} \longrightarrow Y, \hspace{2ex} \mbox{which is the restriction of} \hspace{2ex} \mu : \widetilde{X} \longrightarrow X.$$

\noindent Thus, $\mu : \widetilde{X}\setminus\mu^{-1}(Z) \longrightarrow X\setminus Z$ is a biholomorphism above the complement of the analytic subset $Z:=Y_{sing}$ and $\widetilde{X}$ is a compact K\"ahler manifold, hence so is the submanifold $\widetilde{Y}$. Moreover, $\mu^{\star}(\beta\wedge\omega^{n-p-1})$ is a $C^{\infty}$ semi-positive $(n-p,\,n-p)$-form on $\widetilde{X}$ that is strictly positive on $\widetilde{X}\setminus\mu^{-1}(Z)$. Clearly, $\mu^{\star}\{\alpha\} = \{\mu^{\star}\alpha\}$ is a semi-positive (hence also nef) big class on $\widetilde{X}$ and

$$\mbox{Vol}_{\widetilde{Y}}(\mu^{\star}\{\alpha\}) = \int\limits_{\widetilde{X}}(\mu^{\star}\alpha)^{n-p}\wedge [\widetilde{Y}] = \int\limits_X\alpha^{n-p}\wedge [Y] = \mbox{Vol}_{Y}(\alpha)>0.$$ 

\noindent We consider the following Monge-Amp\`ere equation on the (smooth) compact K\"ahler manifold $\widetilde{Y}$:

\begin{equation}\label{eqn:MA-Ytilde}\widetilde{\alpha}_{\widetilde{Y}}^{n-p} = \mbox{Vol}_{\widetilde{Y}}(\mu^{\star}\{\alpha\})\,\mu^{\star}(\beta\wedge\omega^{n-p-1})  \hspace{3ex} \mbox{on}\hspace{1ex}\widetilde{Y}.\end{equation}

\noindent If the class $\mu^{\star}\{\alpha\}$ were K\"ahler, Yau's Theorem $3$ in [Yau78] on solutions of the Monge-Am\`ere equation with a degenerate (i.e. semi-positive) smooth r.h.s. would yield a unique $d$-closed $(1,\,1)$-current $\widetilde{\alpha}_{\widetilde{Y}}\in\mu^{\star}\{\alpha\}_{|\widetilde{Y}}$ solving equation (\ref{eqn:MA-Ytilde}) such that $\widetilde{\alpha}_{\widetilde{Y}}\geq 0$ on $\widetilde{Y}$, $\widetilde{\alpha}_{\widetilde{Y}}$ is $C^{\infty}$ on $\widetilde{Y}\setminus\mu^{-1}(Z)$ and $\widetilde{\alpha}_{\widetilde{Y}}$ has locally bounded coefficients on $\widetilde{Y}$. In our more general case where the class $\mu^{\star}\{\alpha\}$ is only semi-positive and big, Theorems A, B, C in [BEGZ10] yield a unique $d$-closed $(1,\,1)$-current $\widetilde{\alpha}_{\widetilde{Y}}\in\mu^{\star}\{\alpha\}_{|\widetilde{Y}}$ such that $\widetilde{\alpha}_{\widetilde{Y}}\geq 0$ on $\widetilde{Y}$ and 

$$\langle\widetilde{\alpha}_{\widetilde{Y}}^{n-p}\rangle = \mbox{Vol}_{\widetilde{Y}}(\mu^{\star}\{\alpha\})\,\mu^{\star}(\beta\wedge\omega^{n-p-1})  \hspace{3ex} \mbox{on}\hspace{1ex}\widetilde{Y},$$

\noindent where $\langle\hspace{2ex}\rangle$ stands for the non-pluripolar product introduced in [BEGZ10]. Moreover, $\widetilde{\alpha}_{\widetilde{Y}}$ is $C^{\infty}$ on the ample locus of the class $\mu^{\star}\{\alpha\}_{|\widetilde{Y}}$ (cf. Theorem C in [BEGZ10]), which in our case coincides with $\widetilde{Y}\setminus\mu^{-1}(Z)$, and $\widetilde{\alpha}_{\widetilde{Y}}$ has minimal singularities (cf. Theorem B in [BEGZ10]) among the positive currents in the class $\mu^{\star}\{\alpha\}_{|\widetilde{Y}}$. Since this class contains $C^{\infty}$ semi-positive forms (e.g. $(\mu^{\star}\alpha)_{|\widetilde{Y}}$), its currents with minimal singularties have locally bounded potentials. Thus, $\widetilde{\alpha}_{\widetilde{Y}}$ has locally bounded (and even continuous) potentials, so $\langle\widetilde{\alpha}_{\widetilde{Y}}^{n-p}\rangle$ equals the exterior power $\widetilde{\alpha}_{\widetilde{Y}}^{n-p}$ in the sense of Bedford and Taylor [BT82]. In particular, $[\widetilde{\alpha}_{\widetilde{Y}}^k]_{BC} = [\mu^{\star}\alpha^k]_{BC}$ for all $k$, so $\int_{\widetilde{Y}}\widetilde{\alpha}_{\widetilde{Y}}^{n-p-1}\wedge\mu^{\star}\beta = \int_{\widetilde{Y}}\mu^{\star}\alpha^{n-p-1}\wedge\mu^{\star}\beta$. It remains to set 

$$\widetilde{\alpha}_Y:=\mu_{\star}\widetilde{\alpha}_{\widetilde{Y}}.$$ 

\noindent We thus get a $d$-closed positive $(1,\,1)$-current $\widetilde{\alpha}_Y\in\{\alpha\}_{|Y}$ whose restriction to $Y_{reg}=Y\setminus Z$ is $C^{\infty}$ and which solves the Monge-Amp\`ere equation (\ref{eqn:MA-Y}).

\hfill $\Box$

\vspace{2ex}

\vspace{2ex}

 We can now relate both the psef and the nef thresholds $P^{(\beta)}(\alpha), N^{(\beta)}(\alpha)$ to the volume of $\{\alpha-\beta\}$. The next result confirms Conjecture \ref{Conj:volume_lbound} in the case when these thresholds are sufficiently {\it close} to each other.

\begin{Prop}\label{Prop:psef-nef-vol} Let $X$ be a compact K\"ahler manifold, $\mbox{dim}_{\C}X=n$, and let $\{\alpha\}, \{\beta\}\in H^{1,\,1}_{BC}(X,\,\R)$ be K\"ahler classes such that

\begin{equation}\label{eqn:initial-pos-cond-re}\{\alpha\}^n - n\,\{\alpha\}^{n-1}.\,\{\beta\} >0.\end{equation}

\noindent If either of the following two conditions is satisfied:

\begin{equation}\label{eqn:sufficient-cond}(i)\,\,N^{(\beta)}(\alpha)\geq 1  \hspace{6ex}\mbox{or}\hspace{6ex} (ii)\,\,N^{(\beta)}(\alpha)\geq\frac{\frac{\{\alpha\}^n}{\{\alpha\}^{n-1}.\{\beta\}} - P^{(\beta)}(\alpha)}{n-1},\end{equation}

\noindent then

\begin{equation}\label{eqn:volume_lbound_next}\mbox{Vol}\,\,(\{\alpha - \beta\}) \geq \{\alpha\}^n - n\,\{\alpha\}^{n-1}.\,\{\beta\}.\end{equation}

\end{Prop}

 Note that $P^{(\beta)}(\alpha)\geq\frac{\frac{\{\alpha\}^n}{\{\alpha\}^{n-1}.\{\beta\}} - P^{(\beta)}(\alpha)}{n-1}$ thanks to inequality $(a)$ in (\ref{eqn:P-vol1}). Since $P^{(\beta)}(\alpha)\geq N^{(\beta)}(\alpha)$ (cf. (\ref{equation:psef-nef_ineq})), this shows that condition $(ii)$ requires $N^{(\beta)}(\alpha)$ to be ``close'' to $P^{(\beta)}(\alpha)$. In particular, $(ii)$ holds if $N^{(\beta)}(\alpha)$ and $P^{(\beta)}(\alpha)$ coincide.

\vspace{2ex}

\noindent {\it Proof of Proposition \ref{Prop:psef-nef-vol}.} If $N^{(\beta)}(\alpha)\geq 1$, then the class $\{\alpha - \beta\}$ is nef (cf. Proposition \ref{Prop:nef-threshold_def}), so (\ref{eqn:volume_lbound_next}) follows from Proposition \ref{Prop:nef-case} in this case. 

 Let us now suppose that $N^{(\beta)}(\alpha)< 1$ and that condition $(ii)$ is satisfied. We set $s_0:=N^{(\beta)}(\alpha)$ and $t_0:=P^{(\beta)}(\alpha)$, so $s_0<1<t_0$ (where the last inequality follows from $\{\alpha-\beta\}$ being big --- the main result in [Pop14]). We have

\begin{equation}\label{eqn:convex-combination}\{\alpha-\beta\} = \frac{t_0-1}{t_0-s_0}\,\{\alpha-s_0\beta\} + \frac{1-s_0}{t_0-s_0}\,\{\alpha-t_0\beta\}.\end{equation}

\noindent Since the class $(1-s_0)/(t_0-s_0)\cdot\{\alpha-t_0\beta\}$ is psef, we get the first inequality below:

\begin{equation}\label{eqn:s0-vol_lbound}\mbox{Vol}(\{\alpha-\beta\}) \geq \bigg(\frac{t_0-1}{t_0-s_0}\bigg)^n\,\mbox{Vol}(\{\alpha-s_0\,\beta\}) \geq \bigg(1 - \frac{1-s_0}{t_0-s_0}\bigg)^n\,\bigg(\{\alpha\}^n-ns_0\{\alpha\}^{n-1}.\{\beta\}\bigg),\end{equation}

\noindent where the second inequality follows from Proposition \ref{Prop:nef-case} since the class $\{\alpha-s_0\,\beta\}$ is nef. Let

\begin{equation}\label{eqn:f_def}f:[0,\,1]\rightarrow [0,\,+\infty),  \hspace{3ex} f(s):=\bigg(1 - \frac{1-s}{t_0-s}\bigg)^n\,\bigg(\{\alpha\}^n-ns\{\alpha\}^{n-1}.\{\beta\}\bigg).\end{equation}

\noindent Thus $f(1)=\{\alpha\}^n-n\{\alpha\}^{n-1}.\{\beta\}$ and (\ref{eqn:s0-vol_lbound}) translates to $\mbox{Vol}(\{\alpha-\beta\})\geq f(s_0).$

 We will now show that $f$ is non-increasing on the interval $[\frac{R-t_0}{n-1},\,1]$, where we set: 

\begin{equation}\label{eqn:R-def}R:=\frac{\{\alpha\}^n}{\{\alpha\}^{n-1}.\{\beta\}} \hspace{3ex} \mbox{or equivalently} \hspace{3ex} R=\sup\{r>0\,\,/\,\,\{\alpha\}^n-r\{\alpha\}^{n-1}.\{\beta\}>0\}.\end{equation}

\noindent Assumption (\ref{eqn:initial-pos-cond-re}) means that $R>n$. Deriving $f$, we get:

\begin{eqnarray}\nonumber f'(s) & = & -n\,\frac{(t_0-1)^n}{(t_0-s)^n}\,\{\alpha\}^{n-1}.\{\beta\} + n\,\frac{(t_0-1)^{n-1}}{(t_0-s)^{n-1}}\,\frac{t_0-1}{(t_0-s)^2}\,\bigg(\{\alpha\}^n-ns\{\alpha\}^{n-1}.\{\beta\}\bigg)\\ 
\nonumber & = & n\,\frac{(t_0-1)^n}{(t_0-s)^{n+1}}\,\bigg(\{\alpha\}^n-((n-1)s+t_0)\,\{\alpha\}^{n-1}.\{\beta\}\bigg),  \hspace{6ex} s\in[0,\,1].\end{eqnarray}

\noindent Since $t_0-1>0$ and $t_0-s>0$, the definition of $R$ implies that $f'(s)\leq 0$ for all $s$ such that $(n-1)s+t_0\geq R$, i.e. for all $s\geq\frac{R-t_0}{n-1}$.

 Recall that we are working under the assumption $s_0\in[\frac{R-t_0}{n-1},\,1)$, so from $f$ being non-increasing on $[\frac{R-t_0}{n-1},\,1]$ we infer that $f(s_0)\geq f(1) = \{\alpha\}^n-n\{\alpha\}^{n-1}.\{\beta\}$. Since $\mbox{Vol}(\{\alpha-\beta\})\geq f(s_0)$ by (\ref{eqn:s0-vol_lbound}), we get (\ref{eqn:volume_lbound_next}).  \hfill $\Box$

\subsection{Nef/psef thresholds and volume revisited}\label{subsection:nef/psef/vol_again}   

 We now prove Theorem \ref{The:volume_lbound_nefT}. In so doing, we use a different method for obtaining a lower bound for the volume of $\{\alpha - \beta\}$ that takes into account the ``angles'' between $\{\alpha - s_0\,\beta\}$ and $\{\alpha - t\,\beta\}$ when $t$ varies in a subinterval of $[1,\,t_0)$.  
 
 We start with a useful observation in linear algebra generalising inequality (\ref{eqn:form-ineq}).

\begin{Lem}\label{Lem:alpha-beta-nef-angles} Let $\alpha>0$ and $\beta\geq 0$ be $C^{\infty}$ $(1,\,1)$-forms on an arbitrary complex manifold $X$ with $\mbox{dim}_{\C}X=n$ such that $\alpha - \beta\geq 0$. Then, for every $k\in\{0, 1, \dots , n\}$, the following inequality holds:

\begin{equation}\label{eqn:alpha-beta-nef-angles}(\alpha - \beta)^{n-k}\wedge\alpha^k\geq\alpha^n - (n-k)\,\alpha^{n-1}\wedge\beta.\end{equation}

\end{Lem}

\noindent {\it Proof.} Let $x_0\in X$ be any point and $z_1, \dots , z_n$ local holomorphic coordinates about $x_0$ such that 

$$\alpha =\sum\limits_{j=1}^n i\,dz_j\wedge d\bar{z}_j  \hspace{2ex}\mbox{and}\hspace{2ex} \beta = \sum\limits_{j=1}^n \beta_j\,i\,dz_j\wedge d\bar{z}_j, \hspace{2ex}\mbox{hence}\hspace{2ex} \alpha - \beta = \sum\limits_{j=1}^n (1-\beta_j)\,i\,dz_j\wedge d\bar{z}_j  \hspace{2ex}\mbox{at}\hspace{1ex} x_0.$$

\noindent Thus $\beta_j\in[0,\,1]$ for all $j=1, \dots , n$ by our assumptions and inequality (\ref{eqn:alpha-beta-nef-angles}) at $x_0$ translates to

$$\frac{k!\,(n-k)!}{n!}\,\sum\limits_{1\leq j_1<\dots <j_{n-k}\leq n}(1-\beta_{j_1})\dots (1-\beta_{j_{n-k}}) \geq 1 -\frac{n-k}{n}\,\sum\limits_{l=1}^n\beta_l,$$

\noindent which, in turn, translates to the following inequality after we set $\gamma_j:=1-\beta_j\in[0,\,1]$:

\begin{equation}\label{eqn:gamma_angles}\frac{k!\,(n-k)!}{n!}\,\bigg(\sum\limits_{1\leq j_1<\dots <j_k\leq n}\frac{1}{\gamma_{j_1}\dots\gamma_{j_k}}\bigg)\,\gamma_1\dots\gamma_n \geq \frac{n-k}{n}\,\sum\limits_{l=1}^n\gamma_l + k+1-n.\end{equation}

\noindent Note that the l.h.s. of (\ref{eqn:gamma_angles}) is meaningful even if some $\gamma_j$ vanishes because it reappears in $\gamma_1\dots\gamma_n$. We will prove inequality (\ref{eqn:gamma_angles}) by induction on $n\geq 1$ (where $k\in\{1,\dots , n\}$ is fixed arbitrarily). 
         
 If $n=1$, (\ref{eqn:gamma_angles}) reads $1\geq 1$. Although it is not required by the induction procedure, we now prove (\ref{eqn:gamma_angles}) for $n=3$ and $k=1$ since this case will be used further down, i.e. we prove

\begin{equation}\label{eqn:n=3k=1}\frac{1}{3}\,(\gamma_1\gamma_2 + \gamma_2\gamma_3  + \gamma_3\gamma_1) \geq \frac{2}{3}\,(\gamma_1 + \gamma_2 + \gamma_3) - 1  \hspace{3ex}\mbox{for all}\hspace{1ex} \gamma_1, \gamma_2, \gamma_3\in[0,\,1].\end{equation}

\noindent It is clear that (\ref{eqn:n=3k=1}) is equivalent to $(\gamma_1 -1)\,(\gamma_2 -1) + (\gamma_2 -1)\,(\gamma_3 -1) + (\gamma_3 -1)\,(\gamma_1 -1)\geq 0$ which clearly holds since $\gamma_j -1 \leq 0$ for all $j$.

 Now we perform the induction step. Suppose that we have proved (\ref{eqn:gamma_angles}) for all $1\leq m\leq n$. Proving (\ref{eqn:gamma_angles}) for $n+1$ amounts to proving the following inequality:

\begin{equation}\label{eqn:gamma_angles1}A_{k,\,n+1}:=\frac{k!\,(n+1-k)!}{(n+1)!}\,\sum\limits_{1\leq j_1<\dots <j_k\leq n+1}\frac{\gamma_1\dots\gamma_{n+1}}{\gamma_{j_1}\dots\gamma_{j_k}} \geq \frac{n+1-k}{n+1}\,\sum\limits_{l=1}^{n+1}\gamma_l + k-n.\end{equation}

\noindent The left-hand term $A_{k,\,n+1}$ of (\ref{eqn:gamma_angles1}) can be re-written as

$$\frac{k!\,(n+1-k)!}{(n+1)!}\,\frac{1}{n+1-k}\,\bigg(\gamma_1\,\sum\limits_{\stackrel{\neq 1}{1\leq r_1<\dots <r_{n-k}\leq n+1}}\gamma_{r_1}\dots\gamma_{r_{n-k}}  + \dots + \gamma_{n+1}\,\sum\limits_{\stackrel{\neq n+1}{1\leq r_1<\dots <r_{n-k}\leq n+1}}\gamma_{r_1}\dots\gamma_{r_{n-k}}\bigg),$$ 

\noindent where the meaning of the notation is that the sum whose coefficient is $\gamma_s$ runs over all the ordered indices $r_1<\dots <r_{n-k}$ selected from the set $\{1,\dots , n+1\}\setminus\{s\}$. Now, using inequality (\ref{eqn:gamma_angles}) for $n$ (the induction hypothesis), for every $s\in\{1,\dots , n+1\}$ we get: 

\vspace{1ex}

$\displaystyle\sum\limits_{\stackrel{\neq s}{1\leq r_1<\dots <r_{n-k}\leq n+1}}\gamma_{r_1}\dots\gamma_{r_{n-k}} \geq \frac{n!}{k!\,(n-k)!}\,\bigg(\frac{n-k}{n}\,\sum\limits_{l\in\{1, \dots , n+1\}\setminus\{s\}}\gamma_l + k+1-n\bigg).$

\vspace{1ex}

\noindent Plugging these inequalities into the last (re-written) expression for $A_{k,\,n+1}$, we get:

$$A_{k,\,n+1}\geq\frac{n-k}{n(n+1)}\,\bigg(\sum\limits_{l\in\{1, \dots , n+1\}\setminus\{1\}}\gamma_1\,\gamma_l + \dots + \sum\limits_{l\in\{1, \dots , n+1\}\setminus\{n+1\}}\gamma_{n+1}\,\gamma_l\bigg) + \frac{k+1-n}{n+1}\,(\gamma_1 + \dots + \gamma_{n+1}),$$

\noindent hence \begin{eqnarray}\nonumber A_{k,\,n+1} & \geq & \frac{2(n-k)}{n(n+1)}\,\sum\limits_{1\leq j<k\leq n+1}\gamma_j\,\gamma_k + \frac{k+1-n}{n+1}\,\sum\limits_{l=1}^{n+1}\gamma_l\\
\nonumber & = & \frac{2(n-k)}{n(n+1)}\,\frac{\sum\limits_{1\leq j<k<l\leq n+1}(\gamma_j\,\gamma_k + \gamma_k\,\gamma_l + \gamma_l\,\gamma_j)}{n-1} + \frac{k+1-n}{n+1}\,\sum\limits_{l=1}^{n+1}\gamma_l\\
\nonumber & \stackrel{(a)}{\geq} & \frac{2(n-k)}{(n-1)n(n+1)}\,\bigg(2\sum\limits_{1\leq j<k<l\leq n+1}(\gamma_j + \gamma_k + \gamma_l) - 3\,{n+1 \choose 3}\bigg) + \frac{k+1-n}{n+1}\,\sum\limits_{l=1}^{n+1}\gamma_l\\ 
\nonumber & = & \frac{2(n-k)}{(n-1)n(n+1)}\,\bigg(2\,{n \choose 2}\sum\limits_{l=1}^{n+1}\gamma_l - 3\,{n+1 \choose 3}\bigg) + \frac{k+1-n}{n+1}\,\sum\limits_{l=1}^{n+1}\gamma_l\\
\nonumber & = & \frac{1}{n+1}\, \bigg(\frac{4(n-k)}{n(n-1)}\,\frac{n(n-1)}{2} + k+1-n\bigg)\,\sum\limits_{l=1}^{n+1}\gamma_l - \frac{2(n-k)}{(n-1)n(n+1)}\,3\frac{(n-1)n(n+1)}{2\cdot 3}\\
\nonumber & = & \frac{n-k+1}{n+1}\,\sum\limits_{l=1}^{n+1}\gamma_l - (n-k),\end{eqnarray}

\noindent where inequality $(a)$ above follows from (\ref{eqn:n=3k=1}) applied to each sum $\gamma_j\gamma_k + \gamma_k\gamma_l + \gamma_l\gamma_j$. Thus we have got precisely the inequality (\ref{eqn:gamma_angles1}) that we set out to prove. The proof of Lemma \ref{Lem:alpha-beta-nef-angles} is complete.  \hfill $\Box$

\vspace{3ex}

Now suppose we are in the setting of Conjecture \ref{Conj:volume_lbound}. We keep the notation of $\S.$\ref{subsection:trace-volume}. Recall that $s_0:=N^{(\beta)}(\alpha)$ and $t_0:=P^{(\beta)}(\alpha)$. We assume that $s_0<1$ (since Conjecture \ref{Conj:volume_lbound} has been proved in the case when $s_0\geq 1$).

 We express the class $\{\alpha -\beta\}$ as a convex combination of the nef class $\{\alpha -s_0\,\beta\}$ and the big class $\{\alpha -t\,\beta\}$ for every $t\in[1,\,t_0)$ (cf. Theorem \ref{The:volume_n2_lbound}) in the following more flexible version of (\ref{eqn:convex-combination}):

\begin{equation}\label{eqn:convex-combination-t}\{\alpha-\beta\} = \frac{t-1}{t-s_0}\,\{\alpha-s_0\beta\} + \frac{1-s_0}{t-s_0}\,\{\alpha-t\beta\}, \hspace{3ex} t\in[1,\,t_0).\end{equation}

\noindent We know from Theorem \ref{The:volume_n2_lbound} that for every $t<\frac{R}{n}$ (cf. notation (\ref{eqn:R-def})) there exists a K\"ahler current $T_t$ in the class $\{\alpha-t\beta\}$ such that $T_t\geq (1-\frac{n}{R}\,t)\,\alpha$.
Thus we get the following K\"ahler current in the class $\{\alpha-\beta\}$:

\begin{equation}\label{eqn:St-def} S_t:= \frac{t-1}{t-s_0}\,(\alpha-s_0\beta) + \frac{1-s_0}{t-s_0}\,T_t \geq \frac{t-1}{t-s_0}\,(\alpha-s_0\beta) + \frac{1-s_0}{t-s_0}\,\bigg(1-\frac{n}{R}\,t\bigg)\,\alpha,   \hspace{3ex} t\in\bigg[1,\,\frac{R}{n}\bigg],\end{equation}

\noindent since the class $\{\alpha -s_0\,\beta\}$ being nef allows us to assume without loss of generality that $\alpha-s_0\beta\geq 0$ (after possibly adding $\varepsilon\omega$ and letting $\varepsilon\downarrow 0$ in the end). Since the r.h.s. of (\ref{eqn:St-def}) is smooth, it also provides a lower bound for the absolutely continuous part of $S_t$, so we get the following lower bound for the volume for all $t\in[1,\,\frac{R}{n}]$:

\begin{equation}\label{eqn:vol-angles_lbound}\mbox{Vol}\,(\{\alpha - \beta\}) \geq \int\limits_XS_{t,\,ac}^n \geq \frac{1}{(t-s_0)^n}\,\sum\limits_{k=0}^n{n \choose k} (t-1)^{n-k}(1-s_0)^k\bigg(1-\frac{nt}{R}\bigg)^k\int\limits_X(\alpha-s_0\beta)^{n-k}\wedge\alpha^k.\end{equation}

\noindent Since the class $\{\alpha-s_0\beta\}$ is nef, using Lemma \ref{Lem:alpha-beta-nef-angles}, we get the following

\begin{Lem}\label{Lem:vol_lbound_s0t} Let $X$ be a compact K\"ahler manifold, $\mbox{dim}_{\C}X=n$, and let $\{\alpha\}, \{\beta\}\in H^{1,\,1}_{BC}(X,\,\R)$ be K\"ahler classes such that $\{\alpha\}^n - n\,\{\alpha\}^{n-1}.\,\{\beta\} >0.$ Suppose that $s_0:=N^{(\beta)}(\alpha)< 1$. Then the following estimate holds:

\begin{equation}\label{eqn:volume_lbound_A}\mbox{Vol}\,\,(\{\alpha - \beta\}) \geq \bigg(\frac{A\,t - s_0}{t-s_0}\bigg)^n\,\bigg(\{\alpha\}^n - \frac{s_0(t-1)}{A\,t-s_0}\, n\,\{\alpha\}^{n-1}.\,\{\beta\}\bigg), \hspace{2ex} \mbox{for all}\hspace{1ex} t\in\bigg[1,\,\frac{R}{n}\bigg],\end{equation}

\noindent where we denote $R:=\{\alpha\}^n/\{\alpha\}^{n-1}.\,\{\beta\}>n$ and $A:=1-\frac{n}{R}\,(1-s_0)\in(s_0,\,1)$.

\end{Lem}

\noindent {\it Proof.} From $A-s_0=(1-s_0)\,(1-\frac{n}{R})\in(0,\,1)$ (because $s_0\in(0,\,1)$ and $1-\frac{n}{R}\in(0,\,1)$), we infer that $A>s_0$. That $A<1$ is obvious.

 Without loss of generality, we may assume that $\alpha-s_0\beta\geq 0$, so (\ref{eqn:alpha-beta-nef-angles}) applies to $\alpha$ and $s_0\beta$ and from (\ref{eqn:vol-angles_lbound}) we get:

\begin{eqnarray}\nonumber\mbox{Vol}\,(\{\alpha - \beta\}) & \geq & \frac{1}{(t-s_0)^n}\,\sum\limits_{k=0}^n{n \choose k} (t-1)^{n-k}(1-s_0)^k\bigg(1-\frac{nt}{R}\bigg)^k\bigg(\{\alpha\}^n - (n-k)\,s_0\,\{\alpha\}^{n-1}.\,\{\beta\}\bigg)\\
 \nonumber & = & \frac{1}{(t-s_0)^n}\,\bigg[t-1 + (1-s_0)\,\bigg(1-\frac{nt}{R}\bigg)\bigg]^n\,\{\alpha\}^n\\
\nonumber &   & \hspace{20ex} - \frac{t-1}{(t-s_0)^n}\,\bigg[t-1 + (1-s_0)\,\bigg(1-\frac{nt}{R}\bigg)\bigg]^{n-1}\,ns_0\,\{\alpha\}^{n-1}.\,\{\beta\},\end{eqnarray}

\noindent which proves (\ref{eqn:volume_lbound_A}) since   $t-1 + (1-s_0)\,(1-\frac{nt}{R}) = At-s_0$.  \hfill $\Box$

\vspace{3ex}

 Thus, it becomes necessary to study the variation of the folowing function:

\begin{equation}\label{eqn:g-def}g\,\,:\,\,\bigg[1,\,\frac{R}{n}\bigg]\rightarrow\R, \hspace{3ex} g(t):= \bigg(\frac{A\,t - s_0}{t-s_0}\bigg)^n\,\bigg(\{\alpha\}^n - \frac{s_0(t-1)}{A\,t-s_0}\, n\,\{\alpha\}^{n-1}.\,\{\beta\}\bigg),\end{equation}

\noindent since $\mbox{Vol}\,\,(\{\alpha - \beta\}) \geq g(t)$ for all $t\in[1,\,\frac{R}{n}]$. From $A-s_0 = (1-s_0)\,(1-\frac{n}{R})\in(0,\,1)$, we get:

\begin{equation}\label{eqn:g_ends}g(1)=\bigg(1-\frac{n}{R}\bigg)^n\,\{\alpha\}^n, \hspace{2ex} \mbox{while} \hspace{2ex} g\bigg(\frac{R}{n}\bigg) = \bigg(\frac{R-n}{R-ns_0}\bigg)^n\,(\{\alpha\}^n - ns_0\,\{\alpha\}^{n-1}.\,\{\beta\}).\end{equation}

\noindent We see that $g(1)$ is precisely the lower bound obtained for the volume of $\{\alpha - \beta\}$ in (\ref{eqn:volume_n2_lbound}), so this lower bound will be improved if $g(t)>g(1)$ for some $t\in(1,\,R/n]$.

\vspace{2ex}

\noindent {\bf Variation of $g$.} Since $[(t-1)/(At-s_0)]'= (A-s_0)/(At-s_0)^2$ and $[(At-s_0)/(t-s_0)]'=(1-A)\,s_0/(t-s_0)^2$, for the derivative of $g(t)$ we get: $g'(t) = $

$$n(1-A)\,s_0\,\frac{(At-s_0)^{n-1}}{(t-s_0)^{n+1}}\,\bigg(\{\alpha\}^n - \frac{(ns_0 - ns_0A + A - s_0)\,t - ns_0\,(1-A) - s_0\,(A-s_0)}{(1-A)\,(At-s_0)}\,\{\alpha\}^{n-1}.\,\{\beta\}\bigg).$$

\noindent Now, $At-s_0 > 0$ for all $t\in[1,\,R/n]$ since $At-s_0 \geq A-s_0 = (1-s_0)(1-\frac{n}{R})>0$. Since $t\geq 1>s_0$, from the definition (\ref{eqn:R-def}) of $R$, we get the equivalences: 

\begin{eqnarray}\label{eqn:equiv-AR} g'(t)\geq 0 & \iff & [ns_0\,(1 - A) + A - s_0]\,t - ns_0\,(1-A) - s_0\,(A-s_0) \leq (1-A)\,(At-s_0)\,R\\
\nonumber      & \iff & -[RA^2 -(ns_0-1+R)\,A + (n-1)\,s_0]\,t + s_0\,[A-s_0 + (n-R)\,(1-A)]\geq 0.\end{eqnarray}

\noindent $\bullet$ {\it Sign of $RA^2 -(ns_0-1+R)\,A + (n-1)\,s_0$.} The discriminant of this $2^{nd}$ degree polynomial in $A$ is

\begin{equation}\label{eqn:Delta-A}\Delta_R = R^2 - 2\bigg((n-2)\,s_0 +1\bigg)\,R + (ns_0-1)^2.\end{equation}

\noindent The discriminant of $\Delta_R$ (viewed as a polynomial in $R$) is

\begin{equation}\label{eqn:Delta-R}\Delta'=16(n-1)\,s_0(1-s_0)>0 \hspace{2ex}\mbox{since}\hspace{2ex} s_0\in(0,\,1).\end{equation}

\noindent Thus, the $\Delta_R$ vanishes at $R_1=(n-2)s_0 + 1 -2\,\sqrt{(n-1)s_0(1-s_0)}$ and $R_2=(n-2)s_0 + 1 + 2\,\sqrt{(n-1)s_0(1-s_0)}$.

\begin{Lem}\label{Lem:R12} With our usual notation $R:= \{\alpha\}^n/\{\alpha\}^{n-1}.\,\{\beta\}$, we have: $R_1<R_2\leq n<R$.  \end{Lem}

\noindent {\it Proof.} Only the inequality $R_2\leq n$ needs a proof. It is equivalent to 

\begin{eqnarray}\nonumber (n-2)s_0 + 2\,\sqrt{(n-1)s_0(1-s_0)} \leq n-1 & \iff & 2\,\sqrt{(n-1)s_0(1-s_0)}\leq (n-1)(1-s_0) + s_0\\
\nonumber & \iff & (\sqrt{(n-1)\,(1-s_0)} - \sqrt{s_0})^2 \geq 0,\end{eqnarray}

\noindent which clearly holds.   \hfill $\Box$

\vspace{3ex} 

 The upshot is that $\Delta_R>0$, so $RA^2 -(ns_0-1+R)\,A + (n-1)\,s_0$ vanishes at $A_1 = (ns_0-1+R - \sqrt{\Delta_R})/2R$ and $A_2 = (ns_0-1+R + \sqrt{\Delta_R})/2R$.

\begin{Lem}\label{Lem:A12} With our notation $A:= 1 - \frac{n}{R}\,(1-s_0)\in[0,\,1)$, we have: $A_1<A<A_2$.\end{Lem}

\noindent {\it Proof.} The inequality $A_1<A$ is equivalent to

\begin{equation}\label{eqn:A1A}\frac{ns_0-1+R-\sqrt{\Delta_R}}{2R} < \frac{R-n+ns_0}{R} \iff n(2-s_0)-1-R < \sqrt{\Delta_R}.\end{equation}

\noindent If $n(2-s_0)-1-R\leq 0$, (\ref{eqn:A1A}) is obvious. If $n(2-s_0)-1-R>0$, inequality (\ref{eqn:A1A}) is equivalent to

\begin{eqnarray}\nonumber R^2 + [n(2-s_0)-1]^2 -2[n(2-s_0)-1]R & < & R^2 -2[(n-2)s_0+1]R + (ns_0-1)^2 \iff\\
\nonumber [n(2-s_0)-ns_0]\,[n(2-s_0)+ns_0 +2] & < & 2[n(2-s_0)-(n-2)s_0-2]\,R\\
\nonumber n(n-1)(1-s_0) < (1-s_0)(n-1)R \iff n & < & R,\end{eqnarray}

\noindent where the last inequality holds thanks to our assumption $\{\alpha\}^n-n\{\alpha\}^{n-1}.\,\{\beta\}>0$.

 The inequality $A<A_2$ is equivalent to

\begin{equation}\label{eqn:A2A}\frac{R-n+ns_0}{R} < \frac{ns_0-1+R+\sqrt{\Delta_R}}{2R} \iff R+1-(2-s_0)n < \sqrt{\Delta_R}.\end{equation}

\noindent If $R+1-(2-s_0)n\leq 0$, (\ref{eqn:A2A}) is obvious. If $R+1-(2-s_0)n>0$, inequality (\ref{eqn:A2A}) is equivalent to

\begin{eqnarray}\nonumber R^2 + [1-(2-s_0)n]^2 + 2\,[1-(2-s_0)n]\,R & < & R^2 -2\,[(n-2)s_0+1]\,R + (ns_0-1)^2 \iff\\
\nonumber 2\,[2-(2-s_0)n + (n-2)s_0]\,R & < & [ns_0-(2-s_0)n]\,[ns_0-2+(2-s_0)n] \iff\\ 
\nonumber (n-1)\,(s_0-1)\,R < n(n-1)(s_0-1) \iff R & > & n \hspace{2ex} \mbox{since}\hspace{1ex} s_0-1<0,\end{eqnarray}

\noindent where the last inequality holds thanks to our assumption $\{\alpha\}^n-n\{\alpha\}^{n-1}.\,\{\beta\}>0$.   \hfill $\Box$

\vspace{3ex}

 The obvious corollary of Lemma \ref{Lem:A12} is the following inequality:

\begin{equation}\label{eqn:sign_RA2}RA^2 -(ns_0-1+R)\,A + (n-1)\,s_0 <0.\end{equation}

\noindent $\bullet$ {\it Monotonicity of $g\,\,:\,\,[1,\,\frac{R}{n}]\rightarrow\R$.} Picking up where we left off in (\ref{eqn:equiv-AR}), we get the equivalence:

\begin{eqnarray}\label{eqn:g'sign} g'(t)\geq 0 & \iff & t\geq s_0\,\frac{A-s_0 + (n-R)\,(1-A)}{RA^2 -(ns_0-1+R)\,A + (n-1)\,s_0}.\end{eqnarray}

\begin{Lem}\label{Lem:1s_0-frac} The following inequalities hold:

\begin{equation}\label{eqn:1s_0-frac} (a)\,\,\, A-s_0 + (n-R)\,(1-A) < 0  \hspace{2ex} \mbox{and} \hspace{2ex} (b)\,\,\, 1 > s_0\,\frac{A-s_0 + (n-R)\,(1-A)}{RA^2 -(ns_0-1+R)\,A + (n-1)\,s_0}.\end{equation}

\end{Lem}

\noindent {\it Proof.} $(a)$\, We have: $A-s_0 + (n-R)\,(1-A) = (1-s_0)\,(1-\frac{n}{R}) + \frac{n}{R}\,(1-s_0)\,(n-R) = \frac{(1-s_0)(n-R)(n-1)}{R}$ and the last expression is negative since $n-R<0$ while $1-s_0>0$ and $n-1>0$.

 $(b)$\, Thanks to (\ref{eqn:sign_RA2}), inequality $(b)$ in (\ref{eqn:1s_0-frac}) is equivalent to

\begin{eqnarray}\label{eqn:b-t-g'}\nonumber & & RA^2 -(ns_0-1+R)\,A + (n-1)\,s_0 < s_0\,[A-s_0 + (n-R)\,(1-A)] \iff \\
      & & RA^2 - (Rs_0+R+s_0-1)\,A + s_0\,(s_0+R-1) < 0.\end{eqnarray}

\noindent The discriminant of the l.h.s. in (\ref{eqn:b-t-g'}), viewed as a $2^{nd}$ degree polynomial in $A$, is $\Delta''=(R-1)^2(1-s_0)^2$, so the l.h.s. of (\ref{eqn:b-t-g'}) vanishes at

$$A_3=\frac{R(s_0+1)+s_0-1-(R-1)(1-s_0)}{2R}=s_0 \hspace{2ex} \mbox{and} \hspace{2ex}  A_4 = 1 -\frac{1-s_0}{R}, \hspace{3ex}\mbox{where clearly}\hspace{1ex} A_3<A_4.\hspace{2ex}.$$

\noindent Thus, inequality (\ref{eqn:b-t-g'}) is equivalent to $s_0<A<1 -\frac{1-s_0}{R}$. We have seen in Lemma \ref{Lem:vol_lbound_s0t} that $A>s_0$. On the other hand, proving $A<1 -\frac{1-s_0}{R}$ amounts to proving 

$$1-\frac{n}{R}\,(1-s_0) < 1 -\frac{1-s_0}{R} \iff 1<n \hspace{2ex}(\mbox{since}\hspace{1ex} 1-s_0>0 \hspace{1ex}\mbox{and}\hspace{1ex}R>0).$$

\noindent The last inequality being obvious, the proof of  $(b)$ in (\ref{eqn:1s_0-frac}) is complete.  \hfill $\Box$

\begin{Conc}\label{Conc:g-increasing} Inequality (\ref{eqn:g'sign}) holds strictly for every $t\geq 1$ thanks to part $(b)$ of (\ref{eqn:1s_0-frac}). So, in particular, $g'(t)>0$ for all $t\in[1,\,\frac{R}{n}]$, i.e. the function $g\,\,:\,\,[1,\,\frac{R}{n}]\rightarrow\R$ is increasing.

 Since $\mbox{Vol}\,\,(\{\alpha - \beta\}) \geq g(t)$ for all $t\in[1,\,\frac{R}{n}]$ (cf. Lemma \ref{Lem:vol_lbound_s0t}), the best lower bound for $\mbox{Vol}\,\,(\{\alpha - \beta\})$ that we get through this method in the case when $s_0:=N^{(\beta)}(\alpha)<1$ is

\begin{equation}\label{eqn:volume_lbound_conc}\mbox{Vol}\,\,(\{\alpha - \beta\}) \geq
g\bigg(\frac{R}{n}\bigg) = (\{\alpha\}^n - n\,\{\alpha\}^{n-1}.\,\{\beta\})\,\bigg(\frac{\{\alpha\}^n - n\,\{\alpha\}^{n-1}.\,\{\beta\}}{\{\alpha\}^n - ns_0\,\{\alpha\}^{n-1}.\,\{\beta\}}\bigg)^{n-1}.\end{equation}

\noindent This proves Theorem \ref{The:volume_lbound_nefT}. Note that this lower bound for the volume improves on the lower bound $g(1)$ (cf. (\ref{eqn:g_ends})) obtained in (\ref{eqn:volume_n2_lbound}).

\end{Conc}

\section{Intersection numbers}\label{section:intersection-numbers}

 In this section, we prove Theorem \ref{The:s-intersection_lbound}. We start by deriving analogues in bidegree $(p,\,p)$ with $p\geq 2$ of the inequalities established in $\S.$\ref{section:M-A}. We will use the standard notion of positivity for $(q,\,q)$-forms whose definition is recalled at the beginning of the Appendix before Lemma \ref{Lem:star-ineq}.

\begin{Prop}\label{Prop:power-diff} Let $X$ be a compact K\"ahler manifold with $\mbox{dim}_{\C}X=n$ and let $\alpha, \beta$ be K\"ahler metrics on $X$. Then, for every $t\in[0,\,+\infty)$, every $p\in\{1, \dots , n\}$ and every $C^{\infty}$ positive $(n-p,\,n-p)$-form $\Omega^{n-p,\,n-p}\geq 0$ on $X$ such that $\partial\bar\partial\Omega^{n-p,\,n-p}=0$, we have:

\begin{equation}\label{eqn:power-diff1} \int\limits_X(\alpha^p-tp\,\alpha^{p-1}\wedge\beta)\wedge\Omega^{n-p,\,n-p} \geq \bigg(1 - t\,\frac{n}{R}\bigg)\,\int\limits_X\alpha^p\wedge\Omega^{n-p,\,n-p},\end{equation}

\noindent where, as usual, we let $R:=\frac{\{\alpha\}^n}{\{\alpha\}^{n-1}.\{\beta\}}$. We also have:

\begin{equation}\label{eqn:power-diff2} \int\limits_X(\alpha^p-t^p\beta^p)\wedge\Omega^{n-p,\,n-p} \geq \bigg(1 - t^p\,\frac{{n \choose p}}{R_p}\bigg)\,\int\limits_X\alpha^p\wedge\Omega^{n-p,\,n-p},\end{equation}

\noindent where we let $R_p:=\frac{\{\alpha\}^n}{\{\alpha\}^{n-p}.\{\beta\}^p}$.

\end{Prop}

\noindent {\it Proof.} We may and will assume without loss of generality that $\Omega^{n-p,\,n-p}$ is strictly positive. Inequality (\ref{eqn:power-diff1}) is equivalent to 

$$t\,\frac{n}{R}\,\int\limits_X\alpha^p\wedge\Omega^{n-p,\,n-p} \geq tp\,\int\limits_X\alpha^{p-1}\wedge\beta\wedge\Omega^{n-p,\,n-p},$$

\noindent which, in turn, after the simplification of $t\geq 0$ and the unravelling of $R$, is equivalent to

\begin{equation}\label{eqn:power-diff1-1}\frac{n}{p}\,\bigg(\int\limits_X\alpha^p\wedge\Omega^{n-p,\,n-p}\bigg)\cdot \bigg(\int\limits_X\alpha^{n-1}\wedge\beta\bigg) \geq \{\alpha\}^n\,\int\limits_X\alpha^{p-1}\wedge\beta\wedge\Omega^{n-p,\,n-p}.\end{equation}

\noindent This inequality can be proved using the method in the proof of Lemma \ref{Lem:prod-traces-int}, the pointwise inequality (\ref{eqn:star-ineq-bis}) proved in the Appendix and an {\it approximate fixed point technique} that we now describe. Here are the details. 

\vspace{2ex}

\noindent {\it Approximate fixed point technique} 

\vspace{2ex}

 We consider the following Monge-Amp\`ere equation whose unique $C^{\infty}$ solution in the K\"ahler class $\{\alpha\}$ is denoted by $\widetilde\alpha: = \alpha + i\partial\bar\partial\varphi>0$:

\begin{equation}\label{eqn:M-A_power-diff}\widetilde{\alpha}^n = \frac{\{\alpha\}^n}{\{\alpha\}^{p-1}.\{\beta\}.[\Omega^{n-p,\,n-p}]_A}\,\alpha^{p-1}\wedge\beta\wedge\Omega^{n-p,\,n-p}.\end{equation}

\noindent By $\{\alpha\}^{p-1}.\{\beta\}.[\Omega^{n-p,\,n-p}]_A$ we mean the positive real number $\int_X\alpha^{p-1}\wedge\beta\wedge\Omega^{n-p,\,n-p}$ which clearly depends only on the Bott-Chern classes $\{\alpha\}, \{\beta\}\in H^{1,\,1}(X,\,\R)$ and on the Aeppli class $[\Omega^{n-p,\,n-p}]_A\in H^{n-p,\,n-p}_A(X,\,\R)$. 

 We will vary the form $\alpha$ on the r.h.s. of (\ref{eqn:M-A_power-diff}) in its K\"ahler class $\{\alpha\}$. Let ${\cal E}_{\alpha}:=\{T\in\{\alpha\}\,\,/\,\, T\geq 0\}$ be the set of $d$-closed positive $(1,\,1)$-currents in the K\"ahler class $\{\alpha\}$. Thus ${\cal E}_{\alpha}$ is a {\it compact} convex subset of the locally convex space ${\cal D}^{'1,\,1}(X,\,\R)$ endowed with the weak topology of currents. (The compactness is a consequence of the existence of Gauduchon metrics and holds for any psef class $\{\alpha\}$ even if $X$ is not K\"ahler.) Fix an arbitrary K\"ahler metric $\omega$ in $\{\alpha\}$. For every $\varepsilon>0$, we associate with equation (\ref{eqn:M-A_power-diff}) the map:

\begin{equation}\label{eqn:R-function-def}R_{\varepsilon}:{\cal E}_{\alpha}\rightarrow{\cal E}_{\alpha}, \hspace{2ex} R_{\varepsilon}(T)=\alpha_{T,\,\varepsilon},\end {equation}  

\noindent defined in three steps as follows. Let $T\in{\cal E}_{\alpha}$ be arbitrary. 

\vspace{1ex}

$(i)$\, By the Blocki-Kolodziej version [BK07] for {\it K\"ahler} classes of Demailly's regularisation-of-currents theorem [Dem92, Theorem 1.1], there exist $C^{\infty}$ $d$-closed $(1,\,1)$-forms $\omega_{\varepsilon}\in\{\alpha\} = \{T\}$ for $\varepsilon>0$ such that $\omega_{\varepsilon}\geq -\varepsilon\omega$ and $\omega_{\varepsilon}\rightarrow T$ in the weak topology of currents as $\varepsilon\rightarrow 0$. (The K\"ahler assumption on the class $\{\alpha\}$ crucially ensures that the possible negative part of $\omega_{\varepsilon}$ does not exceed $\varepsilon\omega$, see [BK07].)

 Note that for every sequence of currents $T_j\in{\cal E}_{\alpha}$ converging weakly to a current $T\in{\cal E}_{\alpha}$ and for every fixed $\varepsilon>0$, the sequence of $C^{\infty}$ forms $(\omega_{j,\,\varepsilon})_j$ (obtained by applying to each $T_j$ the Blocki-Kolodziej regularisation procedure just described producing a family $\omega_{j,\,\varepsilon}\rightarrow T_j$ as $\varepsilon\rightarrow 0$) converges in the $C^{\infty}$ topology to the $C^{\infty}$ form $\omega_{\varepsilon}$ (obtained by applying to $T$ the Blocki-Kolodziej regularisation procedure producing a family $\omega_{\varepsilon}\rightarrow T$ as $\varepsilon\rightarrow 0$). In other words, for every fixed $\varepsilon>0$, the map ${\cal E}_{\alpha}\ni T\mapsto\omega_{\varepsilon}\in C^{\infty}_{1,\,1}(X,\,\C)$ is {\it continuous} if ${\cal E}_{\alpha}$ has been equipped with the weak topology of currents and the space of smooth $(1,\,1)$-forms has been given the $C^{\infty}$ topology. 

 To see this, it suffices to work locally with currents $T_j = i\partial\bar\partial\psi_j\geq 0$ and $T=i\partial\bar\partial\psi\geq 0$ for which the psh potentials have the property that $\psi_j\longrightarrow\psi$ in the $L^1$ topology as $j\rightarrow +\infty$, and to show that for every fixed $\varepsilon>0$ we have $i\partial\bar\partial\psi_j\star\rho_{\varepsilon}\longrightarrow i\partial\bar\partial\psi\star\rho_{\varepsilon}$ in the $C^0$ topology as $j\rightarrow +\infty$. (The convergence in the $C^{\infty}$ topology follows from this by taking derivatives.) Indeed, currents are regularised in [BK07] by convolution of their local potentials with regularising kernels $\rho_{\varepsilon}$. Since $i\partial\bar\partial\psi_j\star\rho_{\varepsilon} = \psi_j\star i\partial\bar\partial\rho_{\varepsilon}$ and $i\partial\bar\partial\psi\star\rho_{\varepsilon} = \psi\star i\partial\bar\partial\rho_{\varepsilon}$, we have to ensure, for every fixed $\varepsilon>0$, that

$$\int_{U'}(\psi_j-\psi)(y)\,u_{\varepsilon}(x-y) \underset{j\rightarrow +\infty}{\longrightarrow} 0  \hspace{2ex} \mbox{locally uniformly w.r.t.} \hspace{1ex} x\in U'\Subset U$$

\noindent for every $C^{\infty}$ function $u_{\varepsilon}$ (which is an arbitrary coefficient of $i\partial\bar\partial\rho_{\varepsilon}$ in this case) defined on the open subset $U\subset X$ on which we work. This is clear from the $L^1$ convergence $\psi_j\longrightarrow\psi$ on $U$.

\vspace{1ex}

$(ii)$\, Set $u_{T,\,\varepsilon}:=(1-\varepsilon)\,\omega_{\varepsilon} + \varepsilon\omega$. Thus $u_{T,\,\varepsilon}$ is a K\"ahler metric in the class $\{\alpha\}$ since it is $C^{\infty}$ and $u_{T,\,\varepsilon} \geq -(1-\varepsilon)\,\varepsilon\,\omega + \varepsilon\,\omega = \varepsilon^2\,\omega>0$. Moreover, $u_{T,\,\varepsilon}\rightarrow T$ in the weak topology of currents as $\varepsilon\rightarrow 0$.

\vspace{1ex}

$(iii)$\, Solve equation (\ref{eqn:M-A_power-diff}) with right-hand term defined by $u_{T,\,\varepsilon}$ instead of $\alpha$:

\begin{equation}\label{eqn:MAalpha-beta-p_u}\alpha_{T,\,\varepsilon}^n = \frac{\{\alpha\}^n}{\{\alpha\}^{p-1}.\{\beta\}.[\Omega^{n-p,\,n-p}]_A}\,u_{T,\,\varepsilon}^{p-1}\wedge\beta\wedge\Omega^{n-p,\,n-p}.\end{equation}

\noindent This means that we denote by $\alpha_{T,\,\varepsilon}$ the unique K\"ahler metric in the K\"ahler class $\{\alpha\}$ solving equation (\ref{eqn:MAalpha-beta-p_u}) whose existence is ensured by Yau's theorem [Yau78]. We put $R_{\varepsilon}(T):=\alpha_{T,\,\varepsilon}$. Thus, in particular, the image of $R_{\varepsilon}$ consists of (smooth) K\"ahler metrics in $\{\alpha\}$.

\vspace{1ex}

 Now, $R_{\varepsilon}$ is a continuous self-map of the compact convex subset ${\cal E}_{\alpha}$ of the locally convex space ${\cal D}^{'1,\,1}(X,\,\R)$, so by the Schauder fixed point theorem, there exists a current $T_{\varepsilon}\in{\cal E}_{\alpha}$ such that $T_{\varepsilon} = R_{\varepsilon}(T_{\varepsilon}) = \alpha_{T_{\varepsilon},\,\varepsilon}$. Since $\alpha_{T_{\varepsilon},\,\varepsilon}:=\widetilde\alpha_{\varepsilon}$ is $C^{\infty}$, by construction, the fixed-point current $T_{\varepsilon}$ must be a $C^{\infty}$ form, so $T_{\varepsilon} = \widetilde\alpha_{\varepsilon}$ and $\omega_{\varepsilon}\geq\widetilde\alpha_{\varepsilon} - \delta_{\varepsilon}\omega$ for some $\delta_{\varepsilon}\downarrow 0$ when $\varepsilon\rightarrow 0$. (The last statement follows from the fact that $\omega_{\varepsilon}$ converges in the $C^{\infty}$ topology to $T$ if $T$ is $C^{\infty}$ -- see the explanations under $(3)$ below.) Hence $u_{T_{\varepsilon},\,\varepsilon} = (1-\varepsilon)\,\omega_{\varepsilon} + \varepsilon\omega \geq (1-\varepsilon)\,\widetilde\alpha_{\varepsilon} + [\varepsilon - (1-\varepsilon)\,\delta_{\varepsilon}]\,\omega$. We put $\eta_{\varepsilon}:= \varepsilon - (1-\varepsilon)\,\delta_{\varepsilon}$, so $\eta_{\varepsilon}\rightarrow 0$ when $\varepsilon\rightarrow 0$.

To conclude, for every $\varepsilon>0$, we have got a K\"ahler metric $\widetilde\alpha_{\varepsilon}$ in the K\"ahler class $\{\alpha\}$ such that

\begin{eqnarray}\label{eqn:M-A_power-diff_epsilon}\nonumber\widetilde{\alpha}^n_{\varepsilon} & = & \frac{\{\alpha\}^n}{\{\alpha\}^{p-1}.\{\beta\}.[\Omega^{n-p,\,n-p}]_A}\,[(1-\varepsilon)\,\omega_{\varepsilon} + \varepsilon\omega]^{p-1}\wedge\beta\wedge\Omega^{n-p,\,n-p} \\
  & \geq & (1-\varepsilon)^{p-1}\,\frac{\{\alpha\}^n}{\{\alpha\}^{p-1}.\{\beta\}.[\Omega^{n-p,\,n-p}]_A}\,\widetilde{\alpha}_{\varepsilon}^{p-1}\wedge\beta\wedge\Omega^{n-p,\,n-p} - O(|\eta_{\varepsilon}|),\end{eqnarray}

\noindent where $\omega$ is an arbitrary, fixed K\"ahler metric in the class $\{\alpha\}$ and $O(|\eta_{\varepsilon}|)$ is a quantity that converges to zero as $\varepsilon\rightarrow 0$. The K\"ahler metric $\widetilde\alpha_{\varepsilon}$ can be viewed as an {\it approximate fixed point} in the class $\{\alpha\}$ of equation (\ref{eqn:M-A_power-diff}). 

\vspace{2ex}

\noindent {\bf Explanations.} Here are a few additional comments on the choice of a {\it continuous} regularising operator $R_{\varepsilon}$ for every $\varepsilon>0$. We are indebted to J.-P. Demailly and to A. Zeriahi for many of the ensuing remarks that were left out of the first version of this paper.

\vspace{1ex}

 $(1)$\, The existence of a {\it continuous} regularising operator is an easy consequence of the regularisation theorem (whatever version of it may be used, be it Demailly's regularisation of currents [Dem92, Theorem 1.1] or the Blocki-Kolodziej one [BK07] or any other one) applied to {\it finitely many} currents. The argument for this statement, which has been very kindly and effectively explained to the author by J.-P. Demailly, makes use of the compactness and convexity of ${\cal E}_{\alpha}$ in ${\cal D}^{'1,\,1}(X,\,\R)$. Indeed, the cone of positive currents has a compact and metrisable, hence countable, base. For this reason, there are several different topologies that induce the same topology on this cone ($=$ the weak topology of currents). If we fix a smooth representative $\alpha$ of the class $\{\alpha\}$ (and $\alpha$ can be chosen to be a K\"ahler metric in this case, but this is irrelevant here), any current $T\in{\cal E}_{\alpha}$ can be written as $T=\alpha + i\partial\bar\partial\varphi\geq 0$ for a global quasi-psh (hence $L^2$ and indeed $L^p$ for every $p\in[1,\,+\infty)$) potential $\varphi$ that is unique up to a constant. We can equip the space of potentials $\{\varphi\,\,/\,\,i\partial\bar\partial\varphi \geq -\alpha\}$ with the topology induced by the $L^2$ Hilbert space topology, which is separable, hence has a countable orthonormal base. This topology induces on ${\cal E}_{\alpha}$ the weak topology of currents.

 Now, by compactness of ${\cal E}_{\alpha}$, for every $\varepsilon$, there is a {\it finite} covering of ${\cal E}_{\alpha}$ by open balls of radius $\varepsilon$. Let $T_1, \dots , T_{N_{\varepsilon}}\in{\cal E}_{\alpha}$ be the centres of these balls and let ${\cal E}_{\alpha,\,\varepsilon}$ be the convex polyhedron generated by $T_1, \dots , T_{N_{\varepsilon}}$. We can take $\varepsilon=1/m$ for $m\in\N^{\star}$ and by convexity of ${\cal E}_{\alpha}$ we get

$${\cal E}_{\alpha} = \bigcup\limits_{m=1}^{+\infty}{\cal E}_{\alpha,\,\frac{1}{m}}.$$

\noindent Thus, it suffices to regularise the finitely many currents $T_1, \dots , T_{N_{\varepsilon}}$ and to extend the regularisation to all the currents $T\in{\cal E}_{\alpha,\,\varepsilon}$ by mere convex combinations. This clearly produces a {\it continuous} regularising operator.

\vspace{1ex}

 $(2)$\, The main result of [BK07], namely that in a K\"ahler class positive currents can be regularised with only an $O(\varepsilon)$ loss of positivity (so, ultimately, with no loss at all, as explained above -- hence the K\"ahler metrics in a given K\"ahler class are dense in the positive currents of that class) can also be obtained as an easy consequence of Demailly's regularisation theorem [Dem92]. The argument for this statement, which was very kindly explained to the author by A. Zeriahi, proceeds by first regularising by a mere cut-off operation. Indeed, let $T=\alpha + i\partial\bar\partial\varphi\geq 0$ be an arbitrary positive current in the K\"ahler class $\{\alpha\}$, where $\alpha > 0$ is a K\"ahler metric in this class. For every $\varepsilon >0$, put $T_{\varepsilon}:= \alpha + i\partial\bar\partial\max(\varphi,\,-\frac{1}{\varepsilon})\geq 0$. The current $T_{\varepsilon}$ is still positive since the maximum of any two $\alpha$-psh functions is still $\alpha$-psh when $\alpha$ is a K\"ahler metric ([GZ05, Proposition $2.3$, $(4)$]). We have $\max(\varphi,\,-\frac{1}{\varepsilon})\downarrow\varphi$ pointwise and $T_{\varepsilon}\rightarrow T$ weakly as $\varepsilon\rightarrow 0$. Moreover, the currents $T_{\varepsilon}$ have {\it bounded} potentials, so we can apply Demailly's regularisation theorem [Dem92] to each of them to write $T_{\varepsilon}$ as the weak limit of a sequence of $C^{\infty}$ forms $T_{\varepsilon,\,\delta}\in\{\alpha\}$ as $\delta\rightarrow 0$. Since all the Lelong numbers of $T_{\varepsilon}$ vanish (because the potential is bounded), Demailly's theorem [Dem92] ensures that only a loss of positivity of $O(\delta)$ is introduced by the regularisation process. Taking the diagonal sequence with $\varepsilon = \delta$, we get an approximation of the original current $T$ by $C^{\infty}$ forms in its class with only an $O(\varepsilon)$ loss of positivity.

 The interest in the Blocki-Kolodziej regularisation procedure [BK07] lies in its giving a much simpler proof of the existence of a good regularisation of currents (which is by no means unique) for the special case of a K\"ahler class than Demailly's proof of the general case.

\vspace{1ex}

 $(3)$\, The Blocki-Kolodziej regularisation [BK07] proceeds by convolution of the local potentials of the current $T$ with regularising kernels $\rho_{\varepsilon}$. This method produces a {\it continuous} regularising operator $R_{\varepsilon}$ for every $\varepsilon$. Moreover, if $T$ is a $C^{\infty}$ form in the class $\{\alpha\}$, the $C^{\infty}$ forms $T_{\varepsilon}$ obtained by regularising $T$ converge to $T$ in the $C^{\infty}$ topology as $\varepsilon\rightarrow 0$. This is because locally, if $\psi$ is a $C^{\infty}$ function defined on an open subset $\Omega\subset\C^n$ containing the origin, then all the derivatives of the convolutions $\rho_{\varepsilon}\star\psi$ converge uniformly on all the compact subsets $K\subset\Omega$ to the corresponding derivatives of $\psi$ and the (standard) patching procedure used in [BK07] does not destroy this property. On the other hand, Yau's theorem [Yau78] gives uniform estimates in all the $C^l$ norms of the solution of the Monge-Amp\`ere equation in terms of the r.h.s. term of this equation. Putting these facts together, we get that the regularising operator $R_{\varepsilon}$ obtained by regularisation followed by an application of Yau's theorem is indeed {\it continuous} in the weak topology of currents and, moreover, $R_{\varepsilon}(T)$ converges in the $C^{\infty}$ topology to $T$ whenever $T\in{\cal E}_{\alpha}$ is $C^{\infty}$.

\vspace{2ex}

\noindent {\it Use of the approximate fixed point} 

\vspace{2ex}

 Let us fix any smooth volume form $dV>0$ on $X$. The l.h.s. term in (\ref{eqn:power-diff1-1}) reads:

\vspace{2ex}

$\displaystyle\frac{n}{p}\,\bigg(\int\limits_X\alpha^p\wedge\Omega^{n-p,\,n-p}\bigg)\cdot \bigg(\int\limits_X\alpha^{n-1}\wedge\beta\bigg) = \frac{n}{p}\,\bigg(\int\limits_X\frac{\widetilde{\alpha}_{\varepsilon}^p\wedge\Omega^{n-p,\,n-p}}{dV}\,dV\bigg)\cdot \bigg(\int\limits_X\frac{\widetilde{\alpha}^{n-1}_{\varepsilon}\wedge\beta}{\widetilde{\alpha}_{\varepsilon}^n}\,\frac{\widetilde{\alpha}^n_{\varepsilon}}{dV}\,dV\bigg)$
\begin{eqnarray} \nonumber & \stackrel{(a)}{\geq} & \bigg[\int\limits_X\bigg(\frac{n}{p}\,\frac{\widetilde{\alpha}_{\varepsilon}^p\wedge\Omega^{n-p,\,n-p}}{dV}\,\frac{\widetilde{\alpha}^{n-1}_{\varepsilon}\wedge\beta}{\widetilde{\alpha}_{\varepsilon}^n}\bigg)^{\frac{1}{2}}\, \bigg(\frac{\widetilde{\alpha}^n_{\varepsilon}}{dV}\bigg)^{\frac{1}{2}}\,dV\bigg]^2 \stackrel{(b)}{\geq} \bigg[\int\limits_X\bigg(\frac{\widetilde{\alpha}_{\varepsilon}^{p-1}\wedge\beta\wedge\Omega^{n-p,\,n-p}}{dV}\bigg)^{\frac{1}{2}}\,\bigg(\frac{\widetilde{\alpha}^n_{\varepsilon}}{dV}\bigg)^{\frac{1}{2}}\,dV\bigg]^2 \\
\nonumber & \stackrel{(c)}{\geq} & (1-\varepsilon)^{p-1}\,\frac{\{\alpha\}^n}{\{\alpha\}^{p-1}.\{\beta\}.[\Omega^{n-p,\,n-p}]_A}\,\bigg[\int\limits_X\bigg(\frac{\widetilde{\alpha}_{\varepsilon}^{p-1}\wedge\beta\wedge\Omega^{n-p,\,n-p}}{dV}\bigg)^{\frac{1}{2}}\,\bigg(\frac{\widetilde{\alpha}_{\varepsilon}^{p-1}\wedge\beta\wedge\Omega^{n-p,\,n-p}}{dV}\bigg)^{\frac{1}{2}}\,dV\bigg]^2\\ 
\nonumber & & - O(|\eta_{\varepsilon}|) \\
\nonumber & \stackrel{(d)}{=} & (1-\varepsilon)^{p-1}\,\frac{\{\alpha\}^n}{\{\alpha\}^{p-1}.\{\beta\}.[\Omega^{n-p,\,n-p}]_A}\,\bigg[\int\limits_X\alpha^{p-1}\wedge\beta\wedge\Omega^{n-p,\,n-p}\bigg]^2 - O(|\eta_{\varepsilon}|)\\
\nonumber & = & (1-\varepsilon)^{p-1}\,\{\alpha\}^n\,\int\limits_X\alpha^{p-1}\wedge\beta\wedge\Omega^{n-p,\,n-p} - O(|\eta_{\varepsilon}|),\end{eqnarray}

\noindent for every $\varepsilon>0$. Letting $\varepsilon\rightarrow 0$, we get the desired inequality (\ref{eqn:power-diff1-1}) since $\eta_{\varepsilon}\rightarrow 0$. Inequality $(a)$ was an application of the Cauchy-Schwarz inequality, $(b)$ was an application of the pointwise inequality (\ref{eqn:star-ineq-bis}) of Lemma \ref{Lem:star-ineq} in the Appendix, $(c)$ followed from (\ref{eqn:M-A_power-diff_epsilon}), while identity $(d)$ followed from $\widetilde{\alpha}_{\varepsilon}$ belonging to the class $\{\alpha\}$. The proof of (\ref{eqn:power-diff1}), which is equivalent to (\ref{eqn:power-diff1-1}), is complete.

 The proof of (\ref{eqn:power-diff2}) runs along the same lines. Indeed, (\ref{eqn:power-diff2}) is equivalent to

$$t^p\,\frac{{n \choose p}}{R_p}\,\int\limits_X\alpha^p\wedge\Omega^{n-p,\,n-p}\geq t^p\,\int\limits_X\beta^p\wedge\Omega^{n-p,\,n-p},$$

\noindent which, in turn, after the simplification of $t^p\geq 0$ and the unravelling of $R_p$, is equivalent to

\begin{equation}\label{eqn:power-diff2-1}{n \choose p}\,\bigg(\int\limits_X\alpha^p\wedge\Omega^{n-p,\,n-p}\bigg)\cdot\bigg(\int\limits_X\alpha^{n-p}\wedge\beta^p\bigg)\geq \{\alpha\}^n\,\int\limits_X\beta^p\wedge\Omega^{n-p,\,n-p}.\end{equation}

\noindent The proof of (\ref{eqn:power-diff2-1}) is almost identical to that of (\ref{eqn:power-diff1-1}) spelt out above except for the replacement of equation (\ref{eqn:M-A_power-diff}) with the following Monge-Amp\`ere equation:

\begin{equation}\label{eqn:M-A_power-diff2}\widetilde{\alpha}^n = \frac{\{\alpha\}^n}{\{\beta\}^p.[\Omega^{n-p,\,n-p}]_A}\,\beta^p\wedge\Omega^{n-p,\,n-p}\end{equation}

\noindent and for the replacement of the pointwise inequality (\ref{eqn:star-ineq-bis}) with (\ref{eqn:star-ineq}). Note that, since $\alpha$ does not feature on the r.h.s. of equation (\ref{eqn:M-A_power-diff2}), the {\it approximate fixed point technique} is no longer necessary in this case. It suffices to work with the unique K\"ahler-metric solution $\widetilde{\alpha}$ of (\ref{eqn:M-A_power-diff2}). \hfill $\Box$

\begin{Rem}\label{Rem:exact-fixed-point}{\rm If an exact (rather than an approximate) fixed point for equation (\ref{eqn:M-A_power-diff}) had been sought, we would have needed to consider the following equation in which the K\"ahler-metric solution $\widetilde{\alpha}\in\{\alpha\}$ features on both sides:

\vspace{1ex}

\hspace{20ex}$\displaystyle\widetilde{\alpha}^n = \frac{\{\alpha\}^n}{\{\alpha\}^{p-1}.\{\beta\}.[\Omega^{n-p,\,n-p}]_A}\,\widetilde{\alpha}^{p-1}\wedge\beta\wedge\Omega^{n-p,\,n-p}.$

\vspace{1ex}

\noindent Equations of this type, going back to Donaldson's $J$-flow and to work by Chen, admit a solution under a certain assumption on the class $\{\alpha\}$. See [FLM11] and the references therein for details. Our {\it approximate fixed point technique} does not require any particular assumption on $\{\alpha\}$.}

\end{Rem}
 
\vspace{2ex}

 We can now prove the main result of this section which subsumes Theorem \ref{The:s-intersection_lbound}.

\begin{The}\label{The:alpha-beta_nth} Let $X$ be a compact K\"ahler manifold with $\mbox{dim}_{\C}X=n$ and let $\{\alpha\}, \{\beta\}$ be K\"ahler classes such that $\{\alpha\}^n - n\,\{\alpha\}^{n-1}.\{\beta\}>0$. Then, for every $k\in\{1, 2, \dots , n\}$ and every smooth positive $(n-k,\,n-k)$-form $\Omega^{n-k,\,n-k}\geq 0$ such that $\partial\bar\partial\Omega^{n-k,\,n-k}=0$, the following inequalities hold:

\begin{eqnarray}\label{eqn:alpha-beta-gamma_kth}\nonumber \{\alpha^k-\beta^k\}.[\Omega^{n-k,\,n-k}]_A \stackrel{(I_k)}{\geq} \{\alpha-\beta\}^k.[\Omega^{n-k,\,n-k}]_A & \stackrel{(II_k)}{\geq} & \{\alpha^k-k\,\alpha^{k-1}\wedge\beta\}.[\Omega^{n-k,\,n-k}]_A \\
  & \stackrel{(III_k)}{\geq} & \bigg(1-\frac{n}{R}\bigg)\,\{\alpha\}^k.[\Omega^{n-k,\,n-k}]_A \geq 0,\end{eqnarray}

\noindent where, as usual, $R:=\frac{\{\alpha\}^n}{\{\alpha\}^{n-1}.\{\beta\}}.$ (Thus $R>n$ by assumption.) In particular, $(II_n)$ and $(III_n)$ read:

\begin{eqnarray}\label{eqn:alpha-beta_kth}\{\alpha-\beta\}^n \geq \{\alpha\}^n-n\,\{\alpha\}^{n-1}.\{\beta\} = \bigg(1-\frac{n}{R}\bigg)\,\{\alpha\}^n >0.\end{eqnarray}

\end{The}

\noindent {\it Proof.} We may and will assume without loss of generality that $\Omega^{n-k,\,n-k}$ is strictly positive. 

Inequality $(III_k)$ is nothing but (\ref{eqn:power-diff1}) for $t=1$ and $p=k$. 

 We will now prove $(II_k)$ by induction on $k\in\{1, \dots, n\}$. Let us fix K\"ahler metrics $\alpha, \beta$ in the classes $\{\alpha\}$, resp. $\{\beta\}$. For $k=1$, $(II_1)$ is obviously an identity. Now, proving $(II_k)$ for an arbitrary $k$ amounts to proving that the quantity

\begin{equation}\label{eqn:S_k-def}S_k:=\int\limits_X\bigg((\alpha-\beta)^k-\alpha^k +k\,\alpha^{k-1}\wedge\beta\bigg)\wedge\Omega^{n-k,\,n-k}\end{equation}

\noindent is non-negative. To this end, we first prove the identity:

\begin{equation}\label{eqn:S_k-identity}S_k = \sum\limits_{l=1}^{k-1}l\,\int\limits_X(\alpha-\beta)^{k-l-1}\wedge\alpha^{l-1}\wedge\beta^2\wedge\Omega^{n-k,\,n-k},  \hspace{3ex} k=1, \dots, n.\end{equation}

\noindent This follows immediately by writing the next pointwise identities:

\begin{eqnarray}\nonumber & & (\alpha-\beta)^k-\alpha^k + k\,\alpha^{k-1}\wedge\beta   = -\beta\wedge\sum\limits_{l=1}^k(\alpha-\beta)^{k-l}\wedge\alpha^{l-1} + k\,\alpha^{k-1}\wedge\beta\\
\nonumber & = & \sum\limits_{l=1}^{k-1}\alpha^{l-1}\wedge\beta\wedge\bigg(\alpha^{k-l} - (\alpha-\beta)^{k-l}\bigg) = \sum\limits_{l=1}^{k-1}\alpha^{l-1}\wedge\beta^2\wedge\sum\limits_{r=0}^{k-l-1}\alpha^{k-l-1-r}\wedge(\alpha-\beta)^r\\
\nonumber & = & \sum\limits_{l=1}^{k-1}\sum\limits_{r=0}^{k-l-1} \alpha^{k-r-2}\wedge\beta^2\wedge(\alpha-\beta)^r\\
\nonumber & = & \sum\limits_{r=0}^{k-2} \alpha^{k-r-2}\wedge\beta^2\wedge(\alpha-\beta)^r + \dots + \sum\limits_{r=0}^{k-l-1} \alpha^{k-r-2}\wedge\beta^2\wedge(\alpha-\beta)^r + \dots + \alpha^{k-2}\wedge\beta^2\\
\nonumber & = & \beta^2\wedge(\alpha-\beta)^{k-2} + 2\,\alpha\wedge\beta^2\wedge(\alpha-\beta)^{k-3} + \dots + l\,\alpha^{l-1}\wedge\beta^2\wedge(\alpha-\beta)^{k-l-1} + \dots + (k-1)\,\alpha^{k-2}\wedge\beta^2\\
\nonumber & = & \sum\limits_{l=1}^{k-1}l\,(\alpha-\beta)^{k-l-1}\wedge\alpha^{l-1}\wedge\beta^2.\end{eqnarray}

\noindent This clearly proves (\ref{eqn:S_k-identity}). 

 Now we can run the induction on $k\in\{1, \dots, n\}$ to prove $(II_k)$. Suppose that $(II_1), \dots , (II_{k-1})$ have been proved. Combining them with $(III_k)$ that was proved in (\ref{eqn:power-diff1}) for all $k\in\{1, \dots, n\}$, we deduce that the classes $\{\alpha-\beta\}^{k-r}$ are positive in the following sense:

$$\{\alpha-\beta\}^{k-r}.[\Omega^{n-k+r,\,n-k+r}]_A \geq 0$$

\noindent for all $r\in\{1, \dots , k\}$ and for all $C^{\infty}$ strictly positive $(n-k+r,\, n-k+r)$-forms $\Omega^{n-k+r,\,n-k+r}> 0$ such that $\partial\bar\partial\Omega^{n-k+r,\,n-k+r} = 0$.

 Choosing forms of the shape $\Omega^{n-k+r,\,n-k+r}:= \alpha^{r-2}\wedge\beta^2\wedge\Omega^{n-k,\,n-k}$ with $\Omega^{n-k,\,n-k}>0$ of bidegree $(n-k,\,n-k)$ satisfying $\partial\bar\partial\Omega^{n-k,\,n-k} = 0$, we get:

$$\{\alpha-\beta\}^{k-r}.\{\alpha\}^{r-2}.\{\beta\}^2.[\Omega^{n-k,\,n-k}]_A \geq 0, \hspace{3ex} r\in\{2, \dots , k\}.$$

\noindent Setting $r:=l+1$, this translates to

$$\int\limits_X(\alpha-\beta)^{k-l-1}\wedge\alpha^{l-1}\wedge\beta^2\wedge\Omega^{n-k,\,n-k} \geq 0, \hspace{3ex} l\in\{1, \dots , k-1\},$$

\noindent which means precisely that all the terms in the sum expressing $S_k$ in (\ref{eqn:S_k-identity}) are non-negative. Hence, $S_k\geq 0$, which proves $(II_k)$ (see \ref{eqn:S_k-def}). 

 Let us now prove $(I_k)$ as a consequence of $(II_k)$ and $(III_k)$. For every $k\in\{1, \dots , n\}$, the following pointwise identities are obvious:

\begin{eqnarray}\label{eqn:I_k-II_k}\nonumber \alpha^k-\beta^k-(\alpha-\beta)^k & = & \beta\wedge\sum\limits_{l=0}^{k-1}\alpha^{k-l-1}\wedge(\alpha-\beta)^l - \beta^k = \beta\wedge\bigg(\alpha^{k-1}-\beta^{k-1} + \sum\limits_{l=1}^{k-1}\alpha^{k-l-1}\wedge(\alpha-\beta)^l\bigg)\\
 \nonumber & = &  \beta\wedge\bigg((\alpha-\beta)\wedge\sum\limits_{r=0}^{k-2}\alpha^{k-r-2}\wedge\beta^r + \sum\limits_{l=1}^{k-1}\alpha^{k-l-1}\wedge(\alpha-\beta)^l   \bigg).\end{eqnarray}

\noindent Hence, for every smooth $(n-k,\,n-k)$-form $\Omega^{n-k,\,n-k}\geq 0$ such that $\partial\bar\partial\Omega^{n-k,\,n-k}=0$, we have:

\begin{eqnarray}\nonumber (\{\alpha^k-\beta^k\} - \{\alpha-\beta\}^k).[\Omega^{n-k,\,n-k}]_A  & = & \sum\limits_{r=0}^{k-2}\int\limits_X(\alpha-\beta)\wedge\alpha^{k-r-2}\wedge\beta^{r+1}\wedge\Omega^{n-k,\,n-k}\\
\nonumber & + & \sum\limits_{l=1}^{k-1}\int\limits_X(\alpha-\beta)^l\wedge\alpha^{k-l-1}\wedge\beta\wedge\Omega^{n-k,\,n-k}\\
\nonumber & = & \sum\limits_{r=0}^{k-2}\{\alpha-\beta\}.[\Omega_r^{n-1,\,n-1}]_A + \sum\limits_{l=1}^{k-1}\{\alpha-\beta\}^l.[\Gamma_l^{n-l,\,n-l}]_A\\
\nonumber & \geq 0 &,\end{eqnarray}

\noindent where we have put $\Omega_r^{n-1,\,n-1}:=\alpha^{k-r-2}\wedge\beta^{r+1}\wedge\Omega^{n-k,\,n-k}$ and $\Gamma_l^{n-l,\,n-l}:=\alpha^{k-l-1}\wedge\beta\wedge\Omega^{n-k,\,n-k}$. It is clear that $\Omega_r^{n-1,\,n-1}$ and $\Gamma_l^{n-l,\,n-l}$ are positive $\partial\bar\partial$-closed forms of bidegree $(n-1,\,n-1)$, resp. $(n-l,\,n-l)$, so the last inequality follows from the combination of $(II_k)$ and $(III_k)$. Thus $(I_k)$ is proved.  \hfill $\Box$

\vspace{2ex}

 We immediately get the following consequence of Theorem \ref{The:alpha-beta_nth} which is the analogue of Theorem \ref{The:volume_n2_lbound} in bidegree $(k,\,k)$ for an arbitrary $k$.

\begin{Cor}\label{Cor:pos-current-power-classes} Let $X$ be a compact K\"ahler manifold with $\mbox{dim}_{\C}X=n$ and let $\alpha, \beta>0$ be K\"ahler metrics such that $\{\alpha\}^n - n\,\{\alpha\}^{n-1}.\{\beta\}>0$. Then, for every $k\in\{1, 2, \dots , n\}$, there exist closed positive $(k,\,k)$-currents $U_k\in\{\alpha^k - \beta^k\}$ and $S_k\in\{(\alpha- \beta)^k\}$ such that

\begin{equation}\label{eqn:pos-current-power-classes}U_k\geq \bigg(1-\frac{n}{R}\bigg)\,\alpha^k \hspace{3ex} \mbox{and} \hspace{3ex} S_k\geq \bigg(1-\frac{n}{R}\bigg)\,\alpha^k\end{equation}

\noindent on $X$, where, as usual, we let $R:=\frac{\{\alpha\}^n}{\{\alpha\}^{n-1}.\{\beta\}}$. (So $R>n$ by assumption.)

\end{Cor}

\noindent {\it Proof.} This follows immediately from Theorem \ref{The:alpha-beta_nth} by using the analogue of Lamari's positivity criterion [Lam99, Lemme 3.3] in bidegree $(k,\,k)$ for every $k$.   \hfill $\Box$

\section{A conjecture in the non-K\"ahler context}\label{section:nKapplication}

Let $X$ be a compact complex manifold with $\mbox{dim}_{\C}X=n$. It is standard that if $X$ is of {\it class} ${\cal C}$, then $X$ is both {\it balanced} (i.e. it admits a balanced metric: a Hermitian metric $\omega$ such that $d\omega^{n-1}=0$) by [AB93, Corollary 4.5] and a $\partial\bar\partial$-manifold (i.e. the $\partial\bar\partial$-lemma holds on $X$). On the other hand, there are a great deal of examples of balanced manifolds that are not $\partial\bar\partial$-manifolds (e.g. the Iwasawa manifold), but it is still an open problem to find out whether or not every $\partial\bar\partial$-manifold admits a balanced metric. To the author's knowledge, all the examples of $\partial\bar\partial$-manifolds known so far are also balanced. We now briefly indicate how a generalised version of Demailly's {\it Transcendental Morse Inequalities Conjecture} for a difference of two nef classes might answer a stronger version of this question. The main idea is borrowed from Toma's work [Tom10] in the projective setting and was also exploited in [CRS14] in the K\"ahler setting. 

 It is standard that the canonical linear map induced in cohomology by the identity:

\begin{equation}\label{eqn:mapBC-A}I_{n-1}: H^{n-1,\,n-1}_{BC}(X,\,\C)\rightarrow H^{n-1,\,n-1}_A(X,\,\C), \hspace{3ex} [\Omega]_{BC}\mapsto [\Omega]_A,\end{equation}

\noindent is well defined on every $X$, but it is neither injective, nor surjective in general. Moreover, the balanced cone of $X$ consisting of Bott-Chern cohomology classes of bidegree $(n-1,\,n-1)$ representable by balanced metrics $\omega^{n-1}$:

$${\cal B}_X=\{[\omega^{n-1}]_{BC}\,/\,\omega>0,\hspace{1ex}C^{\infty}\,\,(1,\,1)\mbox{-form such that}\hspace{1ex}d\omega^{n-1}=0 \hspace{1ex} \mbox{on} \hspace{1ex} X\}\subset H^{n-1,\,n-1}_{BC}(X,\,\R),$$

\noindent maps under $I_{n-1}$ to a subset of the Gauduchon cone of $X$ (introduced in [Pop13]) consisting of Aeppli cohomology classes of bidegree $(n-1,\,n-1)$ representable by Gauduchon metrics $\omega^{n-1}$:

$${\cal G}_X=\{[\omega^{n-1}]_A\,/\,\omega>0,\hspace{1ex}C^{\infty}\,\,(1,\,1)\mbox{-form such that}\hspace{1ex}\partial\bar\partial\omega^{n-1}=0 \hspace{1ex} \mbox{on} \hspace{1ex} X\}\subset H^{n-1,\,n-1}_A(X,\,\R).$$ 

\noindent Clearly, the inclusion $I_{n-1}({\cal B}_X)\subset{\cal G}_X$ is strict in general. So is the inclusion $I_{n-1}(\overline{{\cal B}_X}) \subset {\overline{\cal G}_X}$ involving the closures of these two open convex cones.

Now, if $X$ is a $\partial\bar\partial$-manifold, $I_{n-1}$ is an isomorphism of the vector spaces $H^{n-1,\,n-1}_{BC}(X,\,\C)$ and $H^{n-1,\,n-1}_A(X,\,\C)$, as is well known. It is tempting to make the following

\begin{Conj}\label{Conj:cone-equality} If $X$ is a compact $\partial\bar\partial$-manifold of dimension $n$, then $I_{n-1}(\overline{{\cal B}_X}) = {\overline{\cal G}_X}.$

\end{Conj}

 If proved to hold, this conjecture would imply that every $\partial\bar\partial$-manifold is actually balanced since the Gauduchon cone is never empty (due to the existence of Gauduchon metrics by [Gau77]), so the balanced cone would also have to be non-empty in this case. Moreover, a positive answer to this conjecture would have far-reaching implications for a possible future non-K\"ahler mirror symmetry theory since it would remove the ambiguity of choice between the balanced and the Gauduchon cones on $\partial\bar\partial$-manifolds. These two cones would be canonically equivalent on $\partial\bar\partial$-manifolds in this event.

 One piece of evidence supporting Conjecture \ref{Conj:cone-equality} is that it holds on every {\it class} ${\cal C}$ manifold $X$ if the whole of Demailly's {\it Transcendental Morse Inequalities Conjecture} for a difference of two nef classes is confirmed when $X$ is K\"ahler. This is the gist of the observations made in [Tom10] and in [CRS14] alluded to above. Indeed, if $X$ is of {\it class} ${\cal C}$, we may assume without loss of generality that $X$ is actually compact K\"ahler. As proved in [BDPP13], a complete positive answer to Conjecture \ref{Conj:volume_lbound} would imply that the pseudo-effective cone ${\cal E}_X\subset H^{1,\,1}(X,\R)$ of classes of $d$-closed positive $(1,\,1)$-currents $T$ is the dual of the cone ${\cal M}_X\subset H^{n-1,\,n-1}(X,\R)$ of {\it movable classes} (i.e. the closure of the cone generated by classes of currents of the shape $\mu_{\star}(\widetilde{\omega}_1\wedge\dots\wedge\widetilde{\omega}_{n-1})$, where $\mu:\widetilde{X}\rightarrow X$ is any modification of compact K\"ahler manifolds and the $\widetilde{\omega}_j$ are any K\"ahler metrics on $\widetilde{X}$ -- see [BDPP13, Definition 1.3]). Since on $\partial\bar\partial$-manifolds (hence, in particular, on compact K\"ahler ones) the Bott-Chern, Dolbeault and Aeppli cohomologies are canonically equivalent, it is irrelevant in which of these cohomologies the groups $H^{1,\,1}(X,\R)$ and $H^{n-1,\,n-1}(X,\R)$ are considered. 

 The closure $\overline{\cal G}_X\subset H^{n-1,\,n-1}(X,\R)$ of the Gauduchon cone is dual to the pseudo-effective cone ${\cal E}_X\subset H^{1,\,1}(X,\R)$ by Lamari's positivity criterion (Lemma \ref{Lem:Lamari}), while the same kind of argument (i.e. duality and Hahn-Banach) going back to Sullivan shows that the closure $\overline{\cal B}_X\subset H^{n-1,\,n-1}(X,\R)$ of the balanced cone is dual to the cone

$${\cal S}_X=\{[T]_A\,/\,T\geq 0,\,\,T\,\,\mbox{is a}\,(1,\,1)-\mbox{current such that}\hspace{1ex}\partial\bar\partial T=0 \hspace{1ex} \mbox{on} \hspace{1ex} X\}\subset H^{1,\,1}_A(X,\,\R).$$

\noindent Note that ${\cal S}_X$ is closed if $X$ admits a balanced metric $\omega^{n-1}$ (against which the masses of positive $\partial\bar\partial$-closed $(1,\,1)$-currents $T$ can be considered), hence so is it when $X$ is K\"ahler. Thus, by duality, the identity $I_{n-1}(\overline{{\cal B}_X}) = {\overline{\cal G}_X}$ is equivalent to $I_1({\cal E}_X) = {\cal S}_X,$ where $I_1$ is the canonical linear map induced in cohomology by the identity:

\begin{equation}\label{eqn:mapBC-A}I_1: H^{1,\,1}_{BC}(X,\,\C)\rightarrow H^{1,\,1}_A(X,\,\C), \hspace{3ex} [\gamma]_{BC}\mapsto [\gamma]_A.\end{equation}

\noindent In general, $I_1$ is neither injective, nor surjective, but it is an isomorphism when $X$ is a $\partial\bar\partial$-manifold.

 With these facts understood, the identity $I_1({\cal E}_X) = {\cal S}_X$ can be proved when $X$ is K\"ahler (provided that Conjecture \ref{Conj:volume_lbound} can be solved in the affirmative) as explained in [CRS14, Proposition 2.5] by an argument generalising to transcendental classes an earlier argument from [Tom10] that we now recall for the reader's convenience. 

 The inclusion $I_1({\cal E}_X) \subset {\cal S}_X$ is obvious. To prove the reverse inclusion, let $[T]_A\in{\cal S}_X$, i.e. $T\geq 0$ is a $(1,\,1)$-current such that $\partial\bar\partial T=0$. Since $I_1$ is an isomorphism, there exists a unique class $[\gamma]_{BC}\in H^{1,\,1}_{BC}(X,\,\R)$ such that $I_1([\gamma]_{BC}) = [T]_A$. This means that $[\gamma]_A = [T]_A$. We will show that $[\gamma]_{BC}\in{\cal E}_X$. If the [BDPP13] conjecture (predicated on Conjecture \ref{Conj:volume_lbound}) predicting that ${\cal E}_X$ is dual to ${\cal M}_X$ is confirmed, showing that $[\gamma]_{BC}\in{\cal E}_X$ amounts to showing that

\begin{equation}\label{eqn:pos-evaluation}[\gamma]_{BC}\,.\,[\mu_{\star}(\widetilde{\omega}_1\wedge\dots\wedge\widetilde{\omega}_{n-1})]_A \geq 0\end{equation}

\noindent for all modifications $\mu:\widetilde{X}\rightarrow X$ and all K\"ahler metrics $\widetilde{\omega}_j$ on $\widetilde{X}$. On the other hand, Alessandrini and Bassanelli proved in [AB96, Theorem 5.6] the existence and uniqueness of the inverse image under proper modifications $\mu:\widetilde{X}\rightarrow X$ of arbitrary complex manifolds of any positive $\partial\bar\partial$-closed $(1,\,1)$-current $T\geq 0$ in such a way that the Aeppli cohomology class $[T]_A$ is preserved:

$$\exists\, !\,\, (1,\,1)-\mbox{current}\,\,\mu^{\star}T\geq 0\,\,\mbox{on}\,\,\widetilde{X}\hspace{1ex}\mbox{such that}\hspace{1ex}\partial\bar\partial\,(\mu^{\star}T) = 0, \hspace{1ex} [\mu^{\star}T]_A = \mu^{\star}([T]_A) \hspace{1ex}\mbox{and}\hspace{1ex} \mu_{\star}(\mu^{\star}T) = T.$$ 

\noindent (Note that the inverse image $\mu^{\star}([T]_A)$ of any Aeppli class is trivially well defined by taking smooth representatives of the class and pulling them back. Indeed, $\partial\bar\partial$-closedness is preserved, while pullbacks of Aeppli-cohomologous smooth forms are trivially seen to be Aeppli-cohomologous.) Using this key ingredient from [AB96], we get:

\begin{eqnarray}\nonumber [\gamma]_{BC}\,.\,[\mu_{\star}(\widetilde{\omega}_1\wedge\dots\wedge\widetilde{\omega}_{n-1})]_A & = & \int\limits_X\gamma\wedge\mu_{\star}(\widetilde{\omega}_1\wedge\dots\wedge\widetilde{\omega}_{n-1}) = \int\limits_{\widetilde{X}}(\mu^{\star}\gamma)\wedge(\widetilde{\omega}_1\wedge\dots\wedge\widetilde{\omega}_{n-1})\\
\nonumber & = & [\mu^{\star}\gamma]_A\,.\,[\widetilde{\omega}_1\wedge\dots\wedge\widetilde{\omega}_{n-1}]_{BC} = \int\limits_{\widetilde{X}}(\mu^{\star}T)\wedge(\widetilde{\omega}_1\wedge\dots\wedge\widetilde{\omega}_{n-1})\geq 0,\end{eqnarray}

\noindent which proves (\ref{eqn:pos-evaluation}). Note that $\gamma$ and $\mu^{\star}\gamma$ have no sign, so the key point has been the replacement in the integral over $\widetilde{X}$ of $\mu^{\star}\gamma$ by $\mu^{\star}T\geq 0$ which was made possible by $\widetilde{\omega}_1\wedge\dots\wedge\widetilde{\omega}_{n-1}$ being $d$-closed (so we could switch the roles of the Bott-Chern and the Aeppli cohomologies) and by the identity $[\mu^{\star}\gamma]_A = [\mu^{\star}T]_A$ following from $[\gamma]_A = [T]_A$ (see above) and from $[\mu^{\star}T]_A = \mu^{\star}([T]_A)$.

\vspace{2ex}

 The techniques employed in this paper do not seem to be using the full force of the K\"ahler assumption on $X$ and many of the arguments are valid in a more general context. This is part of the justification for proposing Conjecture \ref{Conj:cone-equality}.

\section{Appendix: Hovanskii-Teissier-type inequalities}\label{Appendix}

 In this section, we prove the pointwise inequalities for Hermitian metrics that were used in earlier sections. They generalise the inequality in [Pop14, Lemma 3.1]. 

 For the sake of enhanced flexibility, we shall deal with positive $(q,\,q)$-forms that are not necessarily the $q^{th}$ power of a positive $(1,\,1)$-form. Given any $q\in\{0,\dots , n\}$ and any $C^{\infty}$ real $(q,\,q)$-form $\Omega^{q,\,q}$ on $X$, we make use of the standard notion of (weak) positivity (see e.g. [Dem97, III.1.1]): $\Omega^{q,\,q}$ is said to be {\it positive} (resp. {\it strictly positive}) if for any $(1,\,0)$-forms $\alpha_1, \dots \alpha_{n-q}$, the $(n,\,n)$-form $\Omega^{q,\,q}\wedge i\alpha_1\wedge\overline{\alpha}_1\wedge\dots\wedge i\alpha_{n-q}\wedge\overline{\alpha}_{n-q}$ is non-negative (resp. positive). We write $\Omega^{q,\,q}\geq 0$ (resp. $\Omega^{q,\,q}>0$) in this case. If, in local holomorphic coordinates $z_1, \dots , z_n$, we write

\begin{equation}\label{eqn:Omega_local-form}\frac{\Omega^{q,\,q}}{q!} = \sum\limits_{|L|=|R|=q}\Omega_{L\bar{R}}\,idz_L\wedge d\bar{z}_R,\end{equation}

\noindent then it is clear by considering $\Omega^{q,\,q}\wedge idz_{s_1}\wedge d\bar{z}_{s_1}\wedge\dots\wedge idz_{s_{n-q}}\wedge d\bar{z}_{s_{n-q}}$ that 

\begin{equation}\label{eqn:Omega_pos-cond}\Omega^{q,\,q}\geq 0 \hspace{2ex} \mbox{implies} \hspace{2ex} \Omega_{L\bar{L}}\geq 0 \hspace{2ex} \mbox{for all} \hspace{1ex} L \hspace{1ex} \mbox{with} \hspace{1ex}  |L|=q.\end{equation}

\noindent (We have used the usual notation: $L$ and $R$ stand for ordered multi-indices $L=(1\leq l_1<\dots <l_q\leq n)$, resp. $R=(1\leq r_1<\dots <r_q\leq n)$ of length $q$ and $idz_L\wedge d\bar{z}_R : = idz_{l_1}\wedge d\bar{z}_{r_1}\wedge\dots\wedge idz_{l_q}\wedge d\bar{z}_{r_q}$.)

 In the special case when $\Omega^{q,\,q} = \gamma^q$ for some positive definite smooth $(1,\,1)$-form ($=$ Hermitian metric) $\gamma$ on $X$, if we write

\begin{equation}\label{eqn:gamma-coord}\gamma = \sum\limits_{j,k=1}^n\gamma_{j\bar{k}}\,idz_j\wedge d\bar{z}_k,\end{equation}

\noindent then Sylvester's criterion ensures that $M_{L\bar{L}}(\gamma)>0$ for all multi-indices $L\subset\{1,\dots , n\}$ of any length $l\in\{1,\dots ,n\}$. (For any multi-indices $L, K\subset\{1, \dots , n\}$ of equal lengths, $M_{K\bar{L}}(\gamma)$ denotes the minor of the matrix $(\gamma_{j\bar{k}})_{j,\,k}$ corresponding to the rows with index in $K$ and the columns with index in $L$.) Clearly, $M_{L\bar{L}}(\gamma) = \Omega_{L\bar{L}}$ for all $L$ with $|L|=q$.

\begin{Lem}\label{Lem:star-ineq} Let $\alpha$, $\beta$ be arbitrary Hermitian metrics on a complex manifold $X$ with $\mbox{dim}_{\C}X=n$. 

\noindent The following pointwise inequalities hold for every $p\in\{1, \dots , n\}$ and for every smooth form $\Omega^{n-p,\,n-p}\geq 0$ of bidegree $(n-p,\,n-p)$ on $X$:

\begin{equation}\label{eqn:star-ineq}{n\choose p}\,\frac{\alpha^{n-p}\wedge\beta^p}{\alpha^n}\cdot\frac{\alpha^p\wedge\Omega^{n-p,\,n-p}}{\alpha^n}\geq\frac{\beta^p\wedge\Omega^{n-p,\,n-p}}{\alpha^n}\end{equation}

\noindent and

\begin{equation}\label{eqn:star-ineq-bis}\frac{n}{p}\,\frac{\alpha^{n-1}\wedge\beta}{\alpha^n}\cdot\frac{\alpha^p\wedge\Omega^{n-p,\,n-p}}{\alpha^n}\geq\frac{\alpha^{p-1}\wedge\beta\wedge\Omega^{n-p,\,n-p}}{\alpha^n}.\end{equation}

\end{Lem}

\noindent {\it Proof.} Let us first prove (\ref{eqn:star-ineq}). The special case when $p=1$ and $\Omega^{n-1,\,n-1}=\gamma^{n-1}$ for some $(1,\,1)$-form $\gamma>0$ was proved in [Pop14, Lemma 3.1]. We fix any point $x\in X$ and choose local coordinates $z_1, \dots , z_n$ about $x$ such that

\begin{equation}\label{eqn:coord-choice}\alpha (x) = \sum\limits_{j=1}^nidz_j\wedge d\bar{z}_j   \hspace{2ex}\mbox{and}\hspace{2ex} \beta(x)= \sum\limits_{j=1}^n\beta_j\,idz_j\wedge d\bar{z}_j.\end{equation}

\noindent Thus $\beta_j>0$ for all $j$. At $x$ we get: $\frac{\beta^p}{p!} = \sum\limits_{j_1<\dots <j_p}\beta_{j_1}\dots\beta_{j_p}\,\bigwedge\limits_{l\in\{j_1,\dots , j_p\}}(idz_l\wedge d\bar{z}_l)$, hence

\begin{equation}\label{eqn:alpha_n-p_beta-p}\frac{\alpha^{n-p}\wedge\beta^p}{\alpha^n} = \frac{1}{{n \choose p}}\,\sum\limits_{j_1<\dots <j_p}\beta_{j_1}\dots\beta_{j_p} = \frac{\beta_1\dots \beta_n}{{n \choose p}}\,\,\bigg(\sum\limits_{|K|=n-p}\frac{1}{\beta_K}\bigg) \hspace{2ex} \mbox{at} \hspace{1ex} x,\end{equation}

\noindent where $\beta_K:=\beta_{k_1}\dots\beta_{k_{n-p}}$ whenever $K=(1\leq k_1<\dots k_{n-p}\leq n)$. On the other hand, using (\ref{eqn:Omega_local-form}) with $q=n-p$, we get at $x$:

$$\frac{\alpha^p\wedge\Omega^{n-p,\,n-p}}{\alpha^n} = \frac{1}{{n \choose p}}\,\sum\limits_{|L|=n-p}\Omega_{L\bar{L}}  \hspace{2ex} \mbox{and} \hspace{2ex} \frac{\beta^p\wedge\Omega^{n-p,\,n-p}}{\alpha^n} = \frac{\beta_1\dots \beta_n}{{n \choose p}}\,\sum\limits_{|L|=n-p}\frac{\Omega_{L\bar{L}}}{\beta_L}.$$

\noindent Thus, inequality (\ref{eqn:star-ineq}) at $x$ is equivalent to:

$$\bigg(\sum\limits_{|L|=n-p}\Omega_{L\bar{L}}\bigg)\,\frac{\beta_1\dots \beta_n}{{n \choose p}}\,\bigg(\sum\limits_{|K|=n-p}\frac{1}{\beta_K}\bigg) \geq \frac{\beta_1\dots \beta_n}{{n \choose p}}\,\sum\limits_{|L|=n-p}\frac{\Omega_{L\bar{L}}}{\beta_L},$$

\noindent which clearly holds since $\Omega_{L\bar{L}}\geq 0$ and $\beta_K>0$ for all multi-indices $K, L$. 

 Let us now prove (\ref{eqn:star-ineq-bis}). With the above notation, we have at $x$:

$$\frac{\alpha^{n-1}\wedge\beta}{\alpha^n} = \frac{1}{n}\,\sum\limits_{j=1}^n\beta_j  \hspace{2ex} \mbox{and} \hspace{2ex} \frac{\alpha^{p-1}}{(p-1)!} = \sum\limits_{|J|=p-1}idz_J\wedge d\bar{z}_J,$$

\noindent and the second identity yields at $x$:

$$\frac{\alpha^{p-1}\wedge\beta}{(p-1)!} = \sum\limits_{|J|= p-1}\sum\limits_{j\in\{1, \dots , n\}\setminus J}\beta_j\, idz_j\wedge d\bar{z}_j\wedge idz_J\wedge d\bar{z}_J,$$

\noindent which, in turn, implies the following identity at $x$:

$$\frac{\alpha^{p-1}\wedge\beta\wedge\Omega^{n-p,\,n-p}}{(p-1)!\,(n-p)!} = \sum\limits_{|J|=p-1}\sum\limits_{j\in\{1, \dots , n\}\setminus J}\beta_j\, \Omega_{L_{jJ}\bar{L}_{jJ}}\,\frac{\alpha^n}{n!},$$

\noindent where we have set $\Omega_{L_{jJ}\bar{L}_{jJ}} := \Omega_{L\bar{L}}$ with $L:=\{1, \dots , n\}\setminus (\{j\}\cup J)$ ordered increasingly. Thus, $\{j\}, J$ and $L$ form a partition of $\{1,\dots , n\}$, so any two of them uniquely determine the third.

\noindent Consequently, inequality (\ref{eqn:star-ineq-bis}) at $x$ translates to 

$$\frac{n}{p}\,\frac{1}{{n \choose p}\,n}\,\bigg(\sum\limits_{|L|=n-p}\Omega_{L\bar{L}}\bigg)\,\bigg(\sum\limits_{j=1}^n\beta_j\bigg) \geq \frac{(p-1)!\,(n-p)!}{n!}\,\sum\limits_{|J|=p-1}\sum\limits_{j\in\{1, \dots , n\}\setminus J}\beta_j\, \Omega_{L_{jJ}\bar{L}_{jJ}},$$

\noindent which is clear since $\frac{n}{p}\,\frac{1}{{n \choose p}\,n} = \frac{(p-1)!\,(n-p)!}{n!}$, $\Omega_{L\bar{L}} \geq 0$ for every $L$, $\beta_j>0$ for every $j$ and the terms in the double sum on the r.h.s. of the above inequality are precisely all the products of the shape $\Omega_{L\bar{L}}\,\beta_j$ with $j\notin L$, so they form a subset of the terms on the l.h.s.  \hfill $\Box$

\vspace{3ex}

 Note that inequalities (\ref{eqn:star-ineq}) and (\ref{eqn:star-ineq-bis}) of Lemma \ref{Lem:star-ineq} allow a kind of ``simplification'' of $\alpha^n$ between the numerators and the denominators. For possible future use, we notice a simultaneous reinforcement of inequalities (\ref{eqn:star-ineq}) and (\ref{eqn:star-ineq-bis}) that has not been used in this paper. For this reason and since the proof of the general case involves rather lengthy calculations, we will only prove a special case.

\begin{Lem}\label{Lem:starstar-ineq} Let $\alpha$, $\beta$ be arbitrary Hermitian metrics on a complex manifold $X$ with $\mbox{dim}_{\C}X=n$. Let $p\in\{1, \dots , n\}$ be arbitrary.

 If $\Omega^{n-p,\,n-p}$ is proportional to $\alpha^k\wedge\beta^{n-p-k}$ for some $k\in\{0, \dots , n-p\}$, then the factor ${n\choose p}$ can be omitted from (\ref{eqn:star-ineq}). In other words, for all $p, k\in\{0, \dots , n\}$ such that $p+k\leq n$ we have:

\begin{equation}\label{eqn:starstar-ineq}\frac{\alpha^{n-p}\wedge\beta^p}{\alpha^n}\cdot\frac{\alpha^{p+k}\wedge\beta^{n-p-k}}{\alpha^n}\geq\frac{\alpha^k\wedge\beta^{n-k}}{\alpha^n} \hspace{3ex}\mbox{on}\hspace{1ex} X.\end{equation}

\end{Lem}

\noindent {\it Proof.} We will only prove here the case when $p+k=n-1$, i.e.

\begin{equation}\label{eqn:star-p_1}\frac{\alpha^{n-p}\wedge\beta^p}{\alpha^n}\cdot\frac{\beta\wedge\alpha^{n-1}}{\alpha^n}\geq\frac{\alpha^{n-p-1}\wedge\beta^{p+1}}{\alpha^n} \hspace{3ex}\mbox{on}\hspace{1ex} X,\end{equation}

\noindent which is equivalent to $(\Lambda_{\gamma}\alpha)\,(\Lambda_{\alpha}\beta)\geq n\,\Lambda_{\gamma}\beta$ when $\gamma^{n-1} = t\,\alpha^{n-p-1}\wedge\beta^p$ for some constant $t>0$. Notice that this last inequality improves by a factor $n$ in the special case when $\gamma^{n-1} = t\,\alpha^{n-p-1}\wedge\beta^p$ the general lower bound proved in [Pop14, Lemma 3.1].

 We fix an arbitrary point $x\in X$ and choose local coordinates as in (\ref{eqn:coord-choice}). Using identities analogous to those in the proof of Lemma \ref{Lem:star-ineq}, we see that (\ref{eqn:star-p_1}) translates at $x$ to \begin{equation}\label{eqn:star-p_2}\frac{n-p}{n}\,\bigg(\sum\limits_{j_1<\dots <j_p}\beta_{j_1}\dots\beta_{j_p}\bigg)\,\bigg(\sum\limits_{l=1}^n\beta_l\bigg)\geq (p+1)\,\sum\limits_{k_1<\dots <k_{p+1}}\beta_{k_1}\dots\beta_{k_{p+1}}.\end{equation}

\noindent Now, the l.h.s. of inequality (\ref{eqn:star-p_2}) equals

$$\frac{n-p}{n}\,\bigg((p+1)\,\sum\limits_{k_1<\dots <k_{p+1}}\beta_{k_1}\dots\beta_{k_{p+1}} + \sum\limits_{j_1<\dots <j_p}\beta_{j_1}\dots\beta_{j_p}\,(\beta_{j_1} + \dots + \beta_{j_p})\bigg),$$

\noindent so (\ref{eqn:star-p_2}) is equivalent to

\begin{equation}\label{eqn:star-p_3}(n-p)\,\sum\limits_{j_1<\dots <j_p}\beta_{j_1}\dots\beta_{j_p}\,(\beta_{j_1} + \dots + \beta_{j_p})\geq p(p+1)\,\sum\limits_{k_1<\dots <k_{p+1}}\beta_{k_1}\dots\beta_{k_{p+1}}.\end{equation}

\noindent We will now prove (\ref{eqn:star-p_3}). Let us fix an arbitrary ordered sequence $1\leq k_1<\dots <k_{p+1}\leq n$. For every $r,s\in\{k_1, \dots , k_{p+1}\}$ with $r<s$, we have:

\begin{equation}\label{eqn:prod-squares}2\beta_{k_1}\dots\beta_{k_{p+1}} = (2\beta_r\,\beta_s)\,\prod\limits_{l\notin\{r,s\}}\beta_l \leq \beta_r^2\,\prod\limits_{l\notin\{r,s\}}\beta_l + \beta_s^2\,\prod\limits_{l\notin\{r,s\}}\beta_l,\end{equation}

\noindent where all the products above bear on the indices $l\in\{k_1, \dots , k_{p+1}\}\setminus\{r,s\}$. Note that $\beta_r^2\,\prod_{l\notin\{r,s\}}\beta_l$ is obtained from $\beta_{k_1}\dots\beta_{k_{p+1}}$ by omitting $\beta_s$ and counting $\beta_r$ twice. Summing up these inequalities over all the ${p+1 \choose 2}$ pairs of indices $r<s$ selected from $k_1, \dots , k_{p+1}$, we get

\begin{eqnarray}\label{eqn:k's-fixed-ineq}\nonumber {p+1 \choose 2}\,2\beta_{k_1}\dots\beta_{k_{p+1}} & \leq & \beta_{k_2}\dots\beta_{k_{p+1}}(\beta_{k_2} + \dots + \beta_{k_{p+1}})\\  
\nonumber & + & \beta_{k_1}\,\beta_{k_3}\dots\beta_{k_{p+1}}(\beta_{k_1} + \beta_{k_3} + \dots + \beta_{k_{p+1}})\\
\nonumber & + & \dots\dots\dots\\
    & + & \beta_{k_1}\dots\beta_{k_p}(\beta_{k_1} + \dots + \beta_{k_p}).\end{eqnarray}

\noindent Note that for every $s\in\{1,\dots , p+1\}$, $\beta_{k_s}$ does not feature in the $s^{th}$ line on the r.h.s. of (\ref{eqn:k's-fixed-ineq}). Adding up these inequalities over all the ordered sequences $1\leq k_1<\dots <k_{p+1}\leq n$, we get the desired inequality (\ref{eqn:star-p_3}) because any ordered sequence $1\leq j_1<\dots <j_p\leq n$ occurs inside exactly $(n-p)$ ordered sequences $1\leq k_1<\dots <k_{p+1}\leq n$. Indeed, the extra index for $1\leq k_1<\dots <k_{p+1}\leq n$ can be chosen arbitrarily in $\{1,\dots , n\}\setminus\{j_1,\dots , j_p\}$, so there are $(n-p)$ choices for it.  

 This completes the proof of (\ref{eqn:star-p_3}), hence the proof of (\ref{eqn:starstar-ineq}) when $p+k=n-1$. \hfill $\Box$

\vspace{3ex}

 Again for the record, we notice that an application of Lemma \ref{Lem:starstar-ineq} is an inequality between intersection numbers of cohomology classes reminiscent of the Hovanskii-Teissier inequalities (cf. e.g. [Dem93, Proposition 5.2]). It has an interest of its own.

\begin{Prop}\label{Prop:HTineq} Let $X$ be a compact K\"ahler manifold with $\mbox{dim}_{\C}X=n$ and let $\{\alpha\}, \{\beta\}\in H^{1,\,1}_{BC}(X,\,\R)$ be {\bf nef} Bott-Chern cohomology classes. Then

\begin{equation}\label{eqn:HTineq}(\{\alpha\}^{n-p}.\{\beta\}^p)\,(\{\alpha\}^{p+k}.\{\beta\}^{n-p-k})\geq (\{\alpha\}^n)\,(\{\alpha\}^k.\{\beta\}^{n-k})\end{equation}

\noindent for all $p, k\in\{0, \dots , n\}$ such that $p+k\leq n$.

\end{Prop}

 By the density of the nef cone in the K\"ahler cone, we may assume without loss of generality that $\{\alpha\}$ and $\{\beta\}$ are K\"ahler classes in which we fix respective K\"ahler metrics $\alpha, \beta$.

\vspace{2ex}

\noindent {\it Proof 1 (deduced from a known result).\footnote{This argument was kindly pointed out to the author by S. Boucksom}} For every $j\in\{0, \dots n\}$, let

$$c_j:=\log(\{\alpha\}^j.\{\beta\}^{n-j}).$$

\noindent It is a standard result that the function $j\mapsto c_j$ is {\it concave}. Now, $k\leq n-p\leq n$ and

\begin{equation}\label{eqn:c_j1}n-p = \frac{p}{n-k}\,k + \frac{n-k-p}{n-k}\,n, \hspace{2ex}\mbox{hence, by concavity,} \hspace{2ex} c_{n-p}\geq \frac{p}{n-k}\,c_k + \frac{n-k-p}{n-k}\,c_n.\end{equation}

\noindent Similarly, $k\leq p+k\leq n$ and

 \begin{equation}\label{eqn:c_j2}p+k = \frac{n-p-k}{n-k}\,k + \frac{p}{n-k}\,n, \hspace{2ex}\mbox{hence, by concavity,} \hspace{2ex} c_{p+k}\geq \frac{n-p-k}{n-k}\,c_k + \frac{p}{n-k}\,c_n.\end{equation}

\noindent Taking the sum of (\ref{eqn:c_j1}) and (\ref{eqn:c_j2}), we get: $c_{n-p} + c_{p+k}\geq c_n + c_k,$ which is nothing but (\ref{eqn:HTineq}).  \hfill $\Box$

\vspace{3ex}

\noindent {\it Proof 2.} It uses the pointwise inequality (\ref{eqn:starstar-ineq}) via the technique introduced in [Pop14] and the {\it approximate fixed point} technique introduced in the proof of Proposition \ref{Prop:power-diff}. The arguments are essentially a repetition of some of those used above, so we will only indicate the main points.  

 First notice that the case when $k=0$ is an immediate consequence of the Hovanskii-Teissier inequalities (cf. [Dem93, Proposition 5.2]) which spell:

$$\{\alpha\}^{n-p}.\{\beta\}^p\geq (\{\alpha\}^n)^{\frac{n-p}{n}}\,(\{\beta\}^n)^{\frac{p}{n}} \hspace{2ex}\mbox{and}\hspace{2ex} \{\alpha\}^p.\{\beta\}^{n-p}\geq (\{\alpha\}^n)^{\frac{p}{n}}\,(\{\beta\}^n)^{\frac{n-p}{n}}.$$

\noindent Multiplying these two inequalities, we get (\ref{eqn:HTineq}) for $k=0$.

 For the general case of an arbitrary $k$, we consider the Monge-Amp\`ere equation:

\begin{equation}\label{eqn:M-A_alpha-k_beta-n-k}\widetilde\alpha^n = \frac{\{\alpha\}^n}{\{\alpha\}^k.\{\beta\}^{n-k}}\,\,\alpha^k\wedge\beta^{n-k}, \hspace{2ex}\mbox{or equivalently}\hspace{2ex} \mbox{det}_{\beta}\widetilde\alpha =  \frac{\{\alpha\}^n}{\{\alpha\}^k.\{\beta\}^{n-k}}\,\,\frac{\alpha^k\wedge\beta^{n-k}}{\beta^n},\end{equation}

\noindent for which the {\it approximate fixed point technique} introduced in the proof of Proposition \ref{Prop:power-diff} produces, for every $\varepsilon>0$, a K\"ahler metric $\widetilde\alpha_{\varepsilon}$ in the K\"ahler class $\{\alpha\}$ (in which we have fixed beforehand a K\"ahler metric $\omega$) such that

$$\widetilde\alpha_{\varepsilon}^n = \frac{\{\alpha\}^n}{\{\alpha\}^k.\{\beta\}^{n-k}}\,\,[(1-\varepsilon)\,\omega_{\varepsilon} + \varepsilon\omega]^k\wedge\beta^{n-k} \geq (1-\varepsilon)^k\,\frac{\{\alpha\}^n}{\{\alpha\}^k.\{\beta\}^{n-k}}\,\widetilde\alpha_{\varepsilon}^k\wedge\beta^{n-k} - O(|\eta_{\varepsilon}|),$$

\noindent for some constant $\eta_{\varepsilon}\rightarrow 0$ as $\varepsilon\rightarrow 0$. Hence

\begin{equation}\label{eqn:M-A_alpha-k_beta-n-k_alpha-eps} \mbox{det}_{\beta}\,\widetilde\alpha_{\varepsilon}\geq  (1-\varepsilon)^k\,\frac{\{\alpha\}^n}{\{\alpha\}^k.\{\beta\}^{n-k}}\,\,\frac{\widetilde\alpha_{\varepsilon}^k\wedge\beta^{n-k}}{\beta^n} - O(|\eta_{\varepsilon}|).\end{equation}

 We can now rerun the argument used several times above. For every $\varepsilon>0$, we have:

\begin{eqnarray}\label{eqn:M-A_alpha-k_beta-n-k_appl}\nonumber(\{\alpha\}^{n-p}.\{\beta\}^p)\,(\{\alpha\}^{p+k}.\{\beta\}^{n-p-k}) & = & \bigg(\int\limits_X\frac{\widetilde\alpha_{\varepsilon}^{n-p}\wedge\beta^p}{\beta^n}\,\beta^n\bigg)\,\bigg(\int\limits_X\frac{\widetilde\alpha_{\varepsilon}^{p+k}\wedge\beta^{n-p-k}}{\widetilde\alpha_{\varepsilon}^n}\,(\mbox{det}_{\beta}\widetilde\alpha_{\varepsilon})\,\beta^n\bigg)\end{eqnarray}
\begin{eqnarray}\nonumber  & \stackrel{(a)}{\geq} & \bigg[\int\limits_X\bigg(\frac{\widetilde\alpha_{\varepsilon}^{n-p}\wedge\beta^p}{\beta^n}\,\frac{\widetilde\alpha_{\varepsilon}^{p+k}\wedge\beta^{n-p-k}}{\widetilde\alpha_{\varepsilon}^n}\bigg)^{\frac{1}{2}}\,(\mbox{det}_{\beta}\,\widetilde\alpha_{\varepsilon})^{\frac{1}{2}}\,\beta^n\bigg]^2\\
\nonumber & \stackrel{(b)}{\geq} & \bigg[\int\limits_X\bigg(\frac{\widetilde\alpha_{\varepsilon}^k\wedge\beta^{n-k}}{\beta^n}\bigg)^{\frac{1}{2}}\,(\mbox{det}_{\beta}\widetilde\alpha_{\varepsilon})^{\frac{1}{2}}\,\beta^n\bigg]^2\\
\nonumber & \stackrel{(c)}{\geq} & (1-\varepsilon)^k\,\frac{\{\alpha\}^n}{\{\alpha\}^k.\{\beta\}^{n-k}}\,\, \bigg(\int\limits_X\widetilde\alpha_{\varepsilon}^k\wedge\beta^{n-k}\bigg)^2 - O(|\eta_{\varepsilon}|)\\
\nonumber & = & (1-\varepsilon)^k\,\,\{\alpha\}^n\,\,(\{\alpha\}^k.\{\beta\}^{n-k}) - O(|\eta_{\varepsilon}|).\end{eqnarray}

\noindent As usual, $(a)$ follows from the Cauchy-Schwarz inequality, $(b)$ follows from the pointwise inequality (\ref{eqn:starstar-ineq}), while $(c)$ follows from the inequality (\ref{eqn:M-A_alpha-k_beta-n-k_alpha-eps}). Letting $\varepsilon\rightarrow 0$, we get (\ref{eqn:HTineq}).  \hfill $\Box$

\vspace{6ex}

\noindent {\bf References.} \\

\noindent [AB93]\, L. Alessandrini, G. Bassanelli --- {\it Metric Properties of Manifolds Bimeromorphic to Compact K\"ahler Spaces} --- J. Diff. Geom. {\bf 37} (1993), 95-121.

\vspace{1ex}

\noindent [AB96]\, L. Alessandrini, G. Bassanelli --- {\it The class of Compact Balanced Manifolds Is Invariant under Modifications} --- Complex Analysis and Geometry (Trento, 1993), 1–17, Lecture Notes in Pure and Appl. Math., {\bf 173}, Dekker, New York, 1996.

\vspace{1ex}

\noindent [BT82]\, E. Bedford, B. Taylor --- {\it A New Capacity for Plurisubharmonic Functions.} --- Acta Math. {\bf 149} (1982), 1–40.

\vspace{1ex}

\noindent [BK07]\, Z. Blocki, S. Kolodziej --- {\it On Regularization of Plurisubharmonic Functions on Manifolds} --- Proc. AMS {\bf 135}, no. 7, (2007), 2089-2093.

\vspace{1ex}

\noindent [Bou02]\, S. Boucksom --- {\it On the Volume of a Line Bundle} --- Internat. J. Math. {\bf 13} (2002), no. 10, 1043–1063.

\vspace{1ex}

\noindent [BDPP13]\, S. Boucksom, J.-P. Demailly, M. Paun, T. Peternell --- {\it The Pseudo-effective Cone of a Compact K\"ahler Manifold and Varieties of Negative Kodaira Dimension} --- J. Alg. Geom. {\bf 22} (2013) 201-248.

\vspace{1ex}

\noindent [BEGZ10]\, S. Boucksom, P. Eyssidieux, V. Guedj, A. Zeriahi --- {\it Monge-Amp\`ere Equations in Big Cohomology Classes} --- Acta Math. {\bf 205} (2010), 199-262.

\vspace{1ex}

\noindent [BFJ09]\, S. Boucksom, C. Favre, M. Jonsson --- {\it Differentiability of Volumes of Divisors and a Problem of Teissier} --- J. Alg. Geom {\bf 18} (2009) 279-308. 

\vspace{1ex}

\noindent [CRS14]\, I. Chiose, R. Rasdeaconu, I.Suvaina --- {\it Balanced Metrics on Uniruled Manifolds} --- arXiv e-print DG 1408.4769v1.

\vspace{1ex}

\noindent [Dem85]\, J.-P. Demailly --- {\it Mesures de Monge-Amp\`ere et caract\'erisation g\'eom\'etrique des vari\'et\'es alg\'ebriques affines} --- M\'em. Soc. Math. France (N.S.) {\bf 19} (1985) 1-124.

\vspace{1ex}

\noindent [Dem92]\, J.-P. Demailly --- {\it Regularization of Closed Positive Currents and Intersection Theory} --- J. Alg. Geom., {\bf 1} (1992), 361-409.

\vspace{1ex}

\noindent [Dem93]\, J.-P. Demailly --- {\it A Numerical Criterion for Very Ample Line Bundles} --- J. Diff. Geom. {\bf 37} (1993) 323-374.

\vspace{1ex}

\noindent [Dem 97]\, J.-P. Demailly --- {\it Complex Analytic and Algebraic Geometry}---http://www-fourier.ujf-grenoble.fr/~demailly/books.html

\vspace{1ex}

\noindent [DP04]\, J.-P. Demailly, M. Paun --- {\it Numerical Characterization of the K\"ahler Cone of a Compact K\"ahler Manifold} --- Ann. of Math. {\bf 159} (2004), 1247-1274.

\vspace{1ex}

\noindent [FLM11]\, H. Fang, M. Lai, X. Ma --- {\it On a Class of Fully Nonlinear Flows in K\"ahler Geometry} --- J. reine angew. Math. {\bf 653} (2011), 189—220.

\vspace{1ex}

\noindent [Gau77]\, P. Gauduchon --- {\it Le th\'eor\`eme de l'excentricit\'e nulle} --- C.R. Acad. Sc. Paris, S\'erie A, t. {\bf 285} (1977), 387-390.

\vspace{1ex}

\noindent [GZ05]\, V. Guedj, A. Zeriahi --- {\it Intrinsic Capacities on Compact K\"ahler Manifolds} --- J. Geom. Anal. {\bf 15} (2005), no. 4, 607–639.

\vspace{1ex}

\noindent [Lam99]\, A. Lamari --- {Courants k\"ahl\'eriens et surfaces compactes} --- Ann. Inst. Fourier, Grenoble, {\bf 49}, 1 (1999), 263-285.

\vspace{1ex}

\noindent [Pop13]\, D. Popovici --- {\it Aeppli Cohomology Classes Associated with Gauduchon Metrics on Compact Complex Manifolds} --- arXiv e-print DG 1310.3685v1, to appear in Bull. Soc. Math. France.

\vspace{1ex}

\noindent [Pop14]\, D. Popovici ---{\it Sufficient Bigness Criterion for Differences of Two Nef Classes} --- Math. Ann. DOI 10.1007/s00208-015-1230-z (new title of {\it An Observation Relative to a Paper by J. Xiao} --- arXiv e-print DG 1405.2518v1).

\vspace{1ex}

\noindent [Siu74] \, Y.-T. Siu --- {\it Analyticity of Sets Associated to Lelong Numbers and the Extension of Closed Positive Currents} --- Invent. Math. {\bf 27} (1974), 53-156.

\vspace{1ex}

\noindent [Tom10] \, M. Toma --- {\it A Note on the Cone of Mobile Curves} --- C. R. Math. Acad. Sci. Paris {\bf 348} (2010), 71–73.

\vspace{1ex}

\noindent [Tos15]\, V. Tosatti --- {\it The Calabi-Yau Theorem and K\"ahler Currents} --- arXiv 1505.02124.

\vspace{1ex}

\noindent [TW10]\, V. Tosatti, B. Weinkove --- {\it The Complex Monge-Amp\`ere Equation on Compact Hermitian Manifolds} --- J. Amer. Math. Soc. {\bf 23} (2010), no. 4, 1187-1195.

\vspace{1ex}

\noindent [Xia13]\, J. Xiao --- {\it Weak Transcendental Holomorphic Morse Inequalities on Compact K\"ahler Manifolds} --- to appear in Ann. Inst. Fourier.

\vspace{1ex}

\noindent [Yau78]\, S.T. Yau --- {\it On the Ricci Curvature of a Complex K\"ahler Manifold and the Complex Monge-Amp\`ere Equation} --- Comm. Pure Appl. Math. {\bf 31} (1978) 339-411.

\vspace{3ex}

\noindent Institut de Math\'ematiques de Toulouse, Universit\'e Paul Sabatier,

\noindent 118 route de Narbonne, 31062 Toulouse, France

\noindent Email: popovici@math.univ-toulouse.fr

\end{document}